\documentclass[onefignum,onetabnum]{siamart190516}

\usepackage{makecell}
\usepackage{amsfonts}
\usepackage{amsmath}
\usepackage{amsopn}
\usepackage{amssymb}  %
\usepackage{array}
\usepackage{bold-extra}
\usepackage{braket}
\usepackage{calligra}
\usepackage{epstopdf}
\usepackage{float}
\usepackage{graphicx,epstopdf}
\usepackage{graphicx}
\usepackage{hyperref}
\usepackage{lipsum}
\usepackage{mathrsfs}
\usepackage{mathtools}
\usepackage{pgfplotstable}
\usepackage{pgfplots}
\usepackage{pgfplots}
\usepackage{pifont}
\usepackage{subcaption}
\usepackage{tikz}
\usepackage{wrapfig}
\usepackage{xspace}

\usepackage{paralist}

\usepackage{times}
\usepackage{multirow}

\usepackage{amsmath}
\usepackage{amssymb}
\usepackage{graphicx}
\usepackage{url}
\usepackage{booktabs}
\makeatletter
\g@addto@macro\normalsize{%
  \setlength\abovedisplayskip{3pt}
  \setlength\belowdisplayskip{3pt}
  \setlength\abovedisplayshortskip{3pt}
  \setlength\belowdisplayshortskip{3pt}
}
\makeatother

\usetikzlibrary{arrows.meta}
\usetikzlibrary{calc,patterns,angles,quotes}
\usetikzlibrary{external}
\usetikzlibrary{pgfplots.groupplots}
\colorlet{texcscolor}{blue!50!black}
\colorlet{texemcolor}{red!70!black}
\colorlet{texpreamble}{red!70!black}
\colorlet{codebackground}{black!25!white!25}

\newcommand{\x}{x}
\newcommand{\ix}{x}
\newcommand{\iy}{{y}}
\newcommand{\ip}{{p}}
\newcommand{\iq}{{q}}
\newcommand{\R}{\mathbb R}

\newcommand{\ltwonorm}[1]{\|#1\|_2}
\newcommand{\coloneqq}{:=}
\newcommand{\metricsup}[2]{\|#1\|_{#2}}
\newcommand{\metricdyn}[3]{d_{#2, #3}(#1)}

\newcommand{\poly}[1]{\mathcal P_{#1}}
\newcommand{\SETALG}{\Lambda}
  
\newcommand{\NINTERPOLATION}{\operatorname{{\bf Interp}}}
\newcommand{\NPOS}{\operatorname{{\bf Pos}}}
\newcommand{\NMON}{\operatorname{{\bf Mon}}}
\newcommand{\NINV}{\operatorname{{\bf Inv}}}
\newcommand{\NGRAD}{\operatorname{{\bf Grad}}}
\newcommand{\NHAMILTONIAN}{\operatorname{{\bf Ham}}}
\newcommand{\NSYM}{\operatorname{{\bf Sym}}}

\newcommand{\INTERPOLATION}[1]{\NINTERPOLATION(#1)}
\newcommand{\POS}[1]{\NPOS(#1)}
\newcommand{\MON}[1]{\NMON(#1)}
\newcommand{\INV}[1]{\NINV(#1)}
\newcommand{\GRAD}{\NGRAD}
\newcommand{\HAMILTONIAN}{\NHAMILTONIAN}
\newcommand{\SYM}[1]{\NSYM(#1)}

\newcommand{\LS}{{L_{S,\Omega}}}
\newcommand{\LiS}[1]{{L_{S_{#1},\Omega}}}

\newcommand{\CD}{{C_1^{\circ}(\Omega)}}

\newtheorem{mydef}{Definition}
\newtheorem{myexp}{Example}

\usepackage[shortlabels]{enumitem}
\pgfplotsset{compat = 1.3}
\usepackage{listings}
\definecolor{commentgreen}{RGB}{2,112,10}
\definecolor{eminence}{RGB}{108,48,130}
\definecolor{weborange}{RGB}{255,165,0}
\definecolor{frenchplum}{RGB}{129,20,83}
\lstdefinestyle{JuliaListing}{
  language=Python,
  basicstyle=\small,
  numbers=none,
  stepnumber=1,
  numbersep=10pt,
  tabsize=4,
  showspaces=false,
  showstringspaces=false,
  commentstyle=\color{commentgreen},
  keywordstyle=\color{eminence},
  otherkeywords={@variable,@objective,@constraint,@polyvar},
  frame=tb,
  mathescape=true,
  escapechar=\%,
  captionpos=b,
}
\newcommand*{\CodeComment}[2][5.0cm]{\hfill\makebox[#1][l]{\textcolor{commentgreen}{\#
      #2}}}
\renewcommand{\lstlistingname}{Source Code}% Listing -> Algorithm

\newcommand\bachir[1]{{#1}}

\title{Learning Dynamical Systems with Side
  Information$\dagger$\thanks{ The authors are with the
    department of Operations Research and Financial
    Engineering at Princeton University. \newline This work
    was partially supported by the MURI award of the AFOSR,
    the DARPA Young Faculty Award, the CAREER Award of the
    NSF, the Google Faculty Award, the Innovation Award of
    the School of Engineering and Applied Sciences at
    Princeton University, and the Sloan Fellowship. }}
\date{}

\author{
 Amir Ali Ahmadi and  Bachir El Khadir \\
 \{\texttt{aaa}, \texttt{bkhadir}\}\texttt{@princeton.edu}
}

\begin{document}

\maketitle

\begin{abstract}%
  We present a mathematical and computational framework for
  the problem of learning a dynamical system from noisy
  observations of a few trajectories and subject to
  \emph{side information}. Side information is any knowledge
  we might have about the dynamical system we would like to
  learn besides trajectory data. It is typically inferred
  from domain-specific knowledge or basic principles of a
  scientific discipline. We are interested in explicitly
  integrating side information into the learning process in
  order to compensate for scarcity of trajectory
  observations. %
  We identify six types of side information that arise
  naturally in many applications and lead to convex
  constraints in the learning problem.  First, we show that
  when our model for the unknown dynamical system is
  parameterized as a polynomial, one can impose our side
  information constraints computationally via semidefinite
  programming. We then demonstrate the added value of side
  information for learning the dynamics of basic models in
  physics and cell biology, as well as for learning and
  controlling the dynamics of a model in epidemiology.
  Finally, we study how well polynomial dynamical systems can
  approximate continuously-differentiable ones while
  satisfying side information (either exactly or
  approximately).  Our overall learning methodology combines
  ideas from convex optimization, real algebra, dynamical
  systems, and functional approximation theory, and can
  potentially lead to new synergies between these areas.
\end{abstract}


\makeatletter{\renewcommand*{\@makefnmark}{}
\footnotetext{\hspace{-.8em}$\dagger$ An $8\text{-}$page version of this
  paper \cite{learning_with_sideinfo_short} has appeared in the proceedings of the conference on Learning for Dynamics \& Control
  (L4DC), 2020.}\makeatother}

%%%%%%%%%%%%%%%%%%%%%%%%%%%%%%%%%%%%%%%%%%%%%%%%%%%%%%%%%
% Introduction
%%%%%%%%%%%%%%%%%%%%%%%%%%%%%%%%%%%%%%%%%%%%%%%%%%%%%%%%%

\section{Motivation and problem formulation} 

In several safety-critical applications, one has to learn the
behavior of an unknown dynamical system from noisy
observations of a very limited number of trajectories. For
example, to autonomously land an airplane that has just gone
through engine failure, limited time is available to learn
the modified dynamics of the plane before appropriate control
action can be taken. Similarly, when a new infectious disease
breaks out, few observations are initially available to
understand the dynamics of contagion.  In situations of this
type where data is limited, it is essential to exploit ``side
information''---e.g. physical laws or contextual
knowledge---to assist the task of learning.

%%%%%%%%%%%%%%%%%%%%%%%%%%%%%%%%%%%%%%%%%%%%%%%%%%%%%%%%%
% Problem Formulation
%%%%%%%%%%%%%%%%%%%%%%%%%%%%%%%%%%%%%%%%%%%%%%%%%%%%%%%%%

More formally, our interest in this paper is to learn a continuous-time dynamical system of the form
\begin{equation}
  \dot \ix(t) = f(\ix(t)), %\quad
  \label{eq:autonomous.dynamical.system}
\end{equation}
over a given compact set $\Omega \subset \R^n$ from noisy
observations of a limited number of its trajectories. 
Here, $\dot x(t)$ denotes the time derivative of the state $x(t) \in \mathbb R^n$ at time $t$.
We assume that the unknown vector field 
%$f: \mathbb R^n \rightarrow \mathbb R^n$
$f$ that is to be learned is continuously differentiable over
{an open set containing $\Omega$,}
% \bachir{\footnote{{The reason we require continuous
%       differentiability over an open set containing $\Omega$
%       is because this notion is more naturally defined over
%       open sets.}}}
      an assumption that is often met in
applications.
% and is known to be a sufficient
% condition for existence and uniqueness of solutions to
% \eqref{eq:autonomous.dynamical.system} (see, e.g.,
% \cite{khalil2002nonlinear}).
%
In our setting, we have access
to a training set of the form 
\begin{equation}
 \mathcal D \coloneqq \{ (\ix_i, \iy_i),\quad i = 1,\ldots,N\},\label{eq:training-set}
\end{equation}
where $\ix_i \in \Omega$ (resp.  $\iy_i \in \R^n$) is a possibly noisy measurement of the state of the dynamical system (resp. of $f(x_i)$). Typically, this training set
is obtained from observation of a few trajectories of
\eqref{eq:autonomous.dynamical.system}. The vectors $y_i$
could be either directly accessible (e.g., from sensor
measurements) or approximated from the state variables using a finite-difference scheme.

%
% By learning the dynamical system \eqref{eq:autonomous.dynamical.system} we mean finding a vector field $f_{\mathcal F}$ that best agrees with
% the unknown vector field $f$ among a particular subspace
% $\mathcal F$ of continuously-differentiable functions.

Finding a vector field $f_{\mathcal F}$ that best agrees with
the training set $\mathcal D$ among a particular class
$\mathcal F$ of continuously-differentiable functions amounts
to solving the optimization problem
\begin{equation}
  \label{eq:program_fit_vector_field} f_{\mathcal F} \in \underset{p \in
    \mathcal F}{\arg\min} \sum_{(\ix_i, \iy_i) \in \mathcal D} \ell(p(\ix_i), \iy_i),
\end{equation}
where $\ell(\cdot, \cdot)$ is some loss function that penalizes deviation of $p(x_i)$ from $y_i$. For instance, $\ell(\cdot, \cdot)$ could simply be the $\ell_2$ loss function
\[\ell_2(u, v) \coloneqq \|u - v\|_2 \quad \forall u, v \in \mathbb R^n,\]
though the computational machinery that we propose can readily handle various other convex loss functions (see \Cref{sec:learn-poly-vec}).

%
% the $\ell_1$ loss, the $\ell_\infty$ loss, and
% any loss given by an \emph{sos-convex} function (such as
% $\ell_2$ loss) (see~\cite{helton_sosconvexity_2010} for a
% definition and also \cite[Theorem
% 3.3]{lasserre2009convexity}). 
% In our applications in
% \Cref{sec:applications}, we chose to work with the $\ell_2$ loss given by
% \[\ell_2(x, y) \coloneqq \|x - y\|_2^2 \quad \forall x, y \in \mathbb R^n\]
% for simplicity. 

In addition to fitting  the training set $\mathcal D$, we desire for our learned vector field $f_{\mathcal F}$ to generalize well, i.e., to be consistent as much as possible with the behavior of the unknown vector field $f$ on all of $\Omega$. 
%
%that were not included in the training set.
% In addition to consistency on the training set $\mathcal D$, we desire for our
% learned vector field $f_{\mathcal F}$ to also generalize well
% in conditions that were not observed in the training
% data.
Indeed, the optimization problem
in~\eqref{eq:program_fit_vector_field} only dictates how the
candidate vector field should behave on the training data.
This could easily lead to overfitting, especially if the
function class $\mathcal F$ is large and the observations are
limited. Let us demonstrate this phenomenon by a simple example.

\begin{figure}[H]
  \centering
   \begin{minipage}[c]{0.4\textwidth}
     \includegraphics[scale=.33]{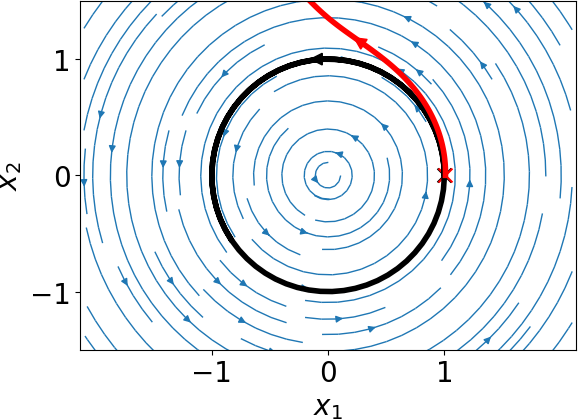}
   \end{minipage}
   \begin{minipage}[c]{0.57\textwidth}
   \caption{\label{fig:example_intro}Streamplot of the vector
    field $f$ in \eqref{eq:example_f_intro} (in blue),
    together with two sample trajectories of the vector field $h$ in \eqref{eq:example_h} with $g(x) = x$
    when started from $(1, 0)^T$ (drawn in black) and from $(1.01, 0)^T$ (drawn in red). The
    trajectories of $f$ and $h$ match exactly when
    started from $(1, 0)^T$, but get arbitrarily far from
    each other when started from $(1.01, 0)^T$.}
\end{minipage}

\end{figure}

\begin{myexp}Consider the two-dimensional vector field
  \begin{equation}
  f(x_1, x_2) \coloneqq (-x_2,x_1)^T.\label{eq:example_f_intro}  
\end{equation}
  The trajectories of
  the system $\dot \ix(t) = f(\ix(t))$ from any initial condition
  are given by circular orbits. In particular, if started
  from the initial condition $x_{\text{init}} = (1, 0)^T$, the trajectory is
  given by $\ix(t, x_{\text{init}}) =(\cos(t), \sin(t))^T$. Hence,
  for any function 
  $g:\mathbb R^2\rightarrow \mathbb R^2,$
  the vector field
  \begin{equation}\label{eq:example_h}
      h(\ix) \coloneqq f(\ix) + (x_1^2+x_2^2-1) g(\ix)
  \end{equation}
  agrees
  with $f$ on the sample trajectory $\ix(t, x_{\text{init}})$.  However,
  the behavior of the trajectories of $h$ depends on the
  arbitrary choice of the function $g$. If $g(\ix) = \ix$ for
  instance, the trajectories of $h$ starting outside of the
  unit disk diverge to infinity. See \cref{fig:example_intro}
  for an illustration.
\end{myexp}

%TODO: change cross to dot

To address the issue of insufficiency of data and to avoid overfitting, we would like to exploit the fact that in many applications, one may have contextual information about the vector field $f$ without knowing $f$ precisely. We call such contextual information \emph{side information}. Formally, every side information is
a subset $S$ of the set of all continuously-differentiable
functions that the vector field $f$ is known to belong to.  Equipped with a list of side information $S_1, \ldots, S_k$, 
our goal is to replace the optimization problem in \eqref{eq:program_fit_vector_field} with
\begin{equation}
  \label{eq:program_fit_vector_field_with_side_info} f_{\mathcal F \cap S_1 \cap \dots \cap S_k} \in \underset{p
\in \mathcal F \cap S_1 \cap \dots \cap S_k} {\arg\min}\sum_{(\ix_i,
\iy_i) \in \mathcal D} \ell(p(\ix_i),\iy_i),
\end{equation}
i.e., to find a vector field $ f_{\mathcal F \cap S_1 \cap \dots \cap S_k}  \in \mathcal F$ that is closest to $f$ on the training set $\mathcal D$ and also satisfies the side information $S_1, \ldots, S_k$ that $f$ is known to satisfy.

\subsection{Outline and contributions of the paper}
In the remainder of this paper, we build on the mathematical formalism we have introduced thus far and make problem \eqref{eq:program_fit_vector_field_with_side_info} more concrete and amenable to computation. In \Cref{sec:side-info}, we identify six notions of side information that are commonly encountered in practice and that have attractive convexity properties, therefore leading to a convex optimization formulation of problem \eqref{eq:program_fit_vector_field_with_side_info}. In \Cref{sec:learn-poly-vec}, we show that when the function class $\mathcal F$ is chosen as the set of polynomial functions of a given degree, then any combination of our six notions of side information can be enforced by semidefinite programming.
The derivation of these semidefinite programs leverages ideas from sum of squares optimization, a concept that we briefly review in the same section for the convenience of the reader.
In \Cref{sec:applications}, we demonstrate the applicability of our approach on three examples from epidemiology, classical mechanics, and cell biology. In each example, we show how common sense and contextual knowledge translate to the notions of side information we present in this paper. Furthermore, in each case, we show that by imposing side information via semidefinite programming, we can learn the behavior of the unknown dynamics from a very limited set of observations.
In our epidemiology example, we also show the benefits of our approach for a downstream task of optimal control (\Cref{sec:opt_control}).
In \Cref{sec:approx_results}, we study the question of how well trajectories of a continuously differentiable vector field that satisfies some side information can be approximated by trajectories of a polynomial vector field that satisfies the same side information either exactly or approximately. %
% we present some density results on approximation of trajectories of continuously differentiable vector fields by trajectories of polynomial vector fields. We a
% we present density results on polynomial vector fields that either exactly or approximately satisfy side information. 
We end the paper with a discussion of future research directions in \Cref{sec:conclusion}.

% \bachir{We emphasize that our aim in this paper is to propose a road map for incorporating side information in one (and not only) method of learning dynamical systems from data. We do not consider this this work to be an alternative to other learning approaches, but our hope is that it stimulates future research on incorporating side information in many other approaches to learning dynamical systems. 
% }

\bachir{We emphasize that our aim in this paper is not to propose our framework as an alternative to other learning algorithms, but rather to present a road map for incorporating side information in the problem of learning dynamical systems from data in general. We make  this road map more explicit by focusing on the common problem of fitting a polynomial vector field to data. However, our hope is that this work stimulates future research on incorporating side information in many other approaches to learning dynamical systems.}

% e connection between optimization over the set polynomial functions  and semidefinite 

% problem \eqref{eq:program_fit_vector_field_with_side_info}, with any combination of our six notions of side information, can be 

% and can
% be enforced in any combination by semidefinite programming
% (SDP). After presenting these notions in
% \Cref{sec:side-info}, we describe the SDP formulation in
% \Cref{sec:learn-poly-vec}, demonstrate the applicability of
% the approach on three examples in \Cref{sec:applications},
% and end with theoretical justification of our methodology in
% \Cref{sec:approx_results}.

% While we work
% with the least-squares loss of simplicity, it turns out that
% our SDP-based approach can readily handle other types of
% losses such as the $\ell_1$ loss, the $\ell_\infty$ loss, and
% any loss given by an \emph{sos-convex} function
% (see~\cite{helton_sosconvexity_2010} for a definition and
% also \cite[Theorem 3.3]{lasserre2009convexity}).

\subsection{Related work}

The idea of using sum of squares and semidefinite
optimization for verifying various properties of a known
dynamical system has been the focus of much research in the
control and optimization communities \cite{PhD:Parrilo,
  PabloGregRekha_BOOK, lasserre_book_2010, boyd_book_1994}. 
%TODO: cite pablo thesis
% Book pablo reckas 
% Text book by lasserre: moments positive polynomials and their applications
% Book by Boyd: Linear Matrix Inequalities in System and Control Theory
Our work borrows some of these techniques to instead impose a desired set of properties on a candidate dynamical system that is to be learned from data.

Learning dynamical systems from data is an important problem
in the field of system identification; \bachir{see e.g. \cite{brunton2019data,keesman2011system, aastrom1971system} and references therein.} 
%TODO: add in response to referee
Various classes of
vector fields have been proposed throughout the years as candidates for the function class $\mathcal F$ in \eqref{eq:program_fit_vector_field}; e.g., reproducing kernel Hilbert
spaces~\cite{sindhwani2018learning,sindhwani2018learnstabilizable,cheng2015learn}, Guassian mixture models~\cite{khansari2011learning}, and neural networks
\cite{vikas_legged_robot_2020,foster2020learning}.
%lewis2012reinforcement,
 %neural_laypunov_nips_2019}, and
%
Some recent approaches to learning dynamical systems from data
impose additional properties on the candidate vector
field. 
%---that could be thought of as
%side information constraints--- . 
These properties include
contraction~
\cite{sindhwani2018learning,khadir2019teleoperator},
stabilizability~
\cite{singh2019learning}, and stability~\cite{kolter2019learning},
and can be thought of as side information. In contrast to our work, imposing these properties requires formulation of nonconvex optimization problems, which can be hard to solve to global optimality. Furthermore, these references impose the desired properties only on sample trajectories (as opposed to the entire space where the properties are known to hold), or introduce an additional layer of nonconvexity to impose the constraints globally.

% That is because it
% would require more elaborate optimization machinery.
% Some of
% these formulations suffer from the disat
% In order
% to impose these additional properties on the candidate vector
% field, these approaches rely on a solving a nonconvex
% problem. Furthermore, 

% We also note that the problem of fitting a polynomial vector field
% to data has appeared e.g. in \cite{schaeffer2018polylasso}, % Bachir: not appeard
% though the focus there is on imposing sparsity of the
% coefficients of the vector field as opposed to side
% information. The closest work in the literature to our work
% is that of Hall on shape-constrained regression \cite[Chapter
% 8]{phd:georgina18}, where similar algebraic techniques are
% used to impose constraints such as convexity and monotonicity
% on a polynomial regressor. See also~\cite{shapeconstraints}
% for some statistical properties of these regressors and for several applications. Our work can be seen as an extension of
% this approach to a dynamical system setting.

We also note that the problem of fitting a polynomial vector field
to data has appeared e.g. in \cite{schaeffer2018polylasso},
though the focus there is on imposing sparsity of the
coefficients of the vector field as opposed to side
information. \bachir{More generally, the problem of finding sparse representations of dynamical systems form data has been studied in \cite{sindy}, where an algorithm for sparse identification of nonlinear dynamics is introduced.} The closest work in the literature to our work
is that of Hall on shape-constrained regression \cite[Chapter
8]{phd:georgina18}, where similar algebraic techniques are
used to impose constraints such as convexity and monotonicity
on a polynomial regressor. See also~\cite{shapeconstraints}
for some statistical properties of these regressors and
several applications. Our work can be seen as an extension of
this approach to a dynamical system setting.

% have been considered in
% . \textcolor{red}{\Large
%   \bf
%   $\ldots$}

% Various learning approaches have been devised for building
% mathematical models from trajectory data. For instance,
% models based on deep and reinforcement learning have been
% considered in
% e.g. \cite{lewis2012reinforcement,lenz2015deepmpc,neural_laypunov_nips_2019}. These
% models are known to require large amounts of data and are
% often expensive to learn. 

% Methods 
% \cite{hills2015algorithm,
% ahmadi2018control}

% RL, deepnets \cite{lewis2012reinforcement,lenz2015deepmpc}

% trigonometric functions, and functions parameterized by
% neural networks
% \cite{neural_laypunov_nips_2019,greydanus2019hamiltonian}.

% control without system identification with sos:
% \cite{guo2020learning}

%%%%%%%%%%%%%%%%%%%%%%%%%%%%%%%%%%%%%%%%%%%%%%%%%%%%%%%%%
% Side info
%%%%%%%%%%%%%%%%%%%%%%%%%%%%%%%%%%%%%%%%%%%%%%%%%%%%%%%%%

\section{Side information}
\label{sec:side-info}

In this section, we identify six types of side information
which we believe are useful in practice (see, e.g.,
\Cref{sec:applications}) and that lead to a convex
formulation of problem
\eqref{eq:program_fit_vector_field_with_side_info}. For example, we will
see in \Cref{sec:learn-poly-vec} that semidefinite
programming can be used to impose any list of side
information constraints of the six types below on a candidate vector field that is parameterized as a polynomial function. The set $\Omega$ that appears in these definitions is a compact
  subset of $\R^n$ over which we would like to learn an unknown vector field $f$. Throughout this paper, the notation $f \in \CD$ denotes that $f$ is continuously differentiable over an open set containing $\Omega$.

% Make left margin smaller

\setlist{leftmargin=5.5mm}

\vspace{1em}

\begin{itemize}

\setlength\itemsep{1em}
  
\item {\bf Interpolation at a finite set of points.} For
  a set of points
  \(\{(\ix_i, \iy_i) \in \Omega \times \R^n\}_{i=1}^m,\) we
  denote by $\INTERPOLATION{\{(\ix_i,y_i)\}_{i=1}^m}$\footnote{{To simplify notation, we drop the dependence of the side information on the set $\Omega$.}}
  the set of vector fields {$f \in \CD$} that satisfy $f(\ix_i)~=~y_i$
  for $i=1,\dots,m$. An important special case is the setting where the vectors $y_i$ are equal to~$0$. In this
  case, the side information is the knowledge of certain
  equilibrium points of the vector field $f$.

\item {\bf Group symmetry.}  For two given \emph{linear
    representations}\footnote{Recall that a linear
    representation of a group $G$ on the vector space
    $\mathbb R^n$ is any group homomorphism from $G$ to the
    group $GL(\mathbb R^n)$ of invertible $n \times n$
    matrices. That is, a linear representation is a map
    $\mu: G \rightarrow GL(\mathbb R^n)$ that satisfies
    $\mu(gg') = \mu(g)\mu(g') \; \forall g, g' \in
    G$. \bachir{See \cite{repr-theory-book} for more details about linear representations of groups.}} $\sigma, \rho: G \rightarrow \mathbb R^{n \times n}$
  of a finite group $G$, with
  $\sigma(g)x \in \Omega \; \forall (x, g) \in \Omega \times
  G$, we define $\SYM{G, \sigma, \rho}$ to be the set of
  {vector fields $f \in \CD$} that satisfy the symmetry
  condition
  \[f(\sigma(g)\ix) = \rho(g) f(\ix) \; \forall \ix \in \Omega,
    \; \forall g \in G.\]
  For example, let $\Omega$ be the unit ball in $\R^n$ and
  consider the group $F_2 = \{1, -1\}$, with scalar
  multiplication as the group operation.
   % of permutations of the set
%   $\{1, 2\}$, and denote its two elements by $e$ (for the
%   identity) and $\tau$ (for the permutation swapping the two
%   points $1$ and $2$).
   \bachir{If we take $\sigma_{F_2}$ to be the
   linear representation of $F_2$ defined by
   $\sigma_{F_2}(1) = -\sigma_{F_2}(-1) = I$, where $I$ denotes the $n\times n$
   identity matrix, then a vector field $f \in \CD$ is in  $\SYM{F_2,\sigma_{F_2},\sigma_{F_2}}$
    if and only if
   \[f(\sigma_{F_2}(1)\ix) = \sigma_{F_2}(1) f(\ix) \text{ and } f(\sigma_{F_2}(-1)\ix) = \sigma_{F_2}(-1) f(\ix)\; \forall \ix \in \Omega,\]
   that is 
   \[f(\ix) = - f(\ix) \; \forall \ix \in \Omega.\]
   In other words, the set
   $\SYM{F_2,\sigma_{F_2},\sigma_{F_2}}$
    is exactly the set of even vector fields in $\CD$. Similarly, if we take  $\rho_{F_2}$ to be the
   linear representation of $F_2$ defined by
   $ \rho_{F_2}(1) =
   \rho_{F_2}( -1 ) = I$, then 
   $\SYM{F_2,\sigma_{F_2},\rho_{F_2}}$ 
   is  the set of odd vector fields in $\CD$.}
   As another example, consider the group $S_n$ of all
   permutations of the set $\{1, \ldots, n\}$, with
   composition as the group operation. If we take
   $\sigma_{S_n}$ to be the map that assigns to an element
   $p \in S_n$ the permutation matrix $P$ obtained by
   shuffling the columns of the identity matrix according to
   $p$, and $\rho_{S_n}$ to be the constant map that assigns
   the identity matrix to every $p~\in~S_n$, then the set
   $\SYM{S_n, \sigma_{S_n}, \rho_{S_n}}$ is the set of
   \emph{symmetric} functions in $\CD$, i.e., functions in
   $\CD$ that are invariant under permutations of their
   arguments. We remark that a finite combination of side
   information of type $\NSYM$ can often be written
   equivalently as a single side information of type
   $\NSYM$. For example, the set
   $\SYM{F_2,\sigma_{F_2},\rho_{F_2}} \cap
   \SYM{S_n,\sigma_{S_n},\rho_{S_n}}$ of even symmetric
   functions in $\CD$ is equal to
   $\SYM{F_2 \times S_n,\sigma,\rho}$, where $F_2 \times S_n$
   is the direct product of $F_2$ and $S_n$, $\sigma$ is given
   by
   $\sigma(g, g') = \sigma_{F_2}(g)\sigma_{S_n}(g') \,
   \forall (g, g') \in F_2 \times S_n$, and $\rho$ is the
   constant map that assigns the identity matrix to every
   element in $F_2 \times S_n$.

\item {\bf Coordinate nonnegativity.}  For given sets
  $P_i, N_i \subseteq {\Omega}$, $i=1,\ldots,n$, we denote by
  $\POS{ \{(P_i, N_i)\}_{i=1}^n }$ the set of
  vector
  fields {$f \in \CD$} that satisfy
  \[f_i(\ix) \ge 0 \; \forall \ix \in P_i, \text{ and } f_i(\ix) \le 0 \; \forall \ix \in N_i,\;
    \forall i \in \{1, \dots, n\}.\]
  These
  constraints are useful when we know that certain components
  of the state vector are increasing or decreasing
  functions of time in some regions of the state
  space.\footnote{There is no loss of generality in assuming
    that each coordinate of the vector field is nonnegative
    or nonpositive on a single set since one can always
    reduce multiple sets to one by taking unions. The same
    comment applies to the side information of coordinate
    directional monotonicity that is defined next.}

\item {\bf Coordinate directional monotonicity.} For given
  sets $P_{ij}, N_{ij} \subseteq \Omega$, $i,j = 1,\ldots,n$,
  we denote by $\MON{ \{(P_{ij}, N_{ij})\}_{i, j=1}^n}$ the
  set of vector fields {$f \in \CD$} that satisfy
  \[\frac {\partial f_i}{\partial  x_j} (\ix) \ge 0 \;
    \forall \ix \in P_{ij}, \text{ and } \frac {\partial f_i}{\partial  x_j} (\ix) \le 0 \;
    \forall \ix \in N_{ij},\; \forall i, j \in \{1, \dots,
    n\}.\]
 See \cref{fig:coord_mon} for an illustration of a simple example.
  
    \begin{figure}[H]
%todo: violate montonocity cdts for which N or P is empty

 \begin{minipage}[c]{0.4\textwidth}
\centering
\begin{tikzpicture}[scale=.75]
      \draw[black] (0, 0) -- (1, 0);
      \draw[black] (1, 0) -- (2, 0);
      \draw[black, -{Latex[length=2.0mm, width=1mm]}] (2, 0) -- (5, 0) node[right]{$x_1$};
      \draw[black, -{Latex[length=2.0mm, width=1mm]}] (0, 0) -- (0, 3) node[right]{$x_2$};
      \draw[black, -{Latex[length=2.0mm, width=1mm]}] (0, 0) -- (.5, .5);
      \draw[black, -{Latex[length=2.0mm, width=1mm]}] (1.5, 0) -- (2, 1.2);
      \draw[black, -{Latex[length=2.0mm, width=1mm]}] (3, 0) -- (2.5, 1.9) node[right]{$f(x_1, x_2)$};
\end{tikzpicture}
\end{minipage}\hfill
\begin{minipage}[c]{0.57\textwidth}
  %\captionsetup{margin={1.5cm, 0cm, 0cm, 0cm}}
\caption{\label{fig:coord_mon}An example of the behavior of a vector field $f: \mathbb R^2 \rightarrow \mathbb R^2$ satisfying $\MON{ \{(P_{ij},N_{ij})\}_{i, j=1}^2}$ with $P_{21} = N_{11} = [0, 1] \times \{0\}$  (i.e., $\frac{\partial f_2}{\partial x_1}(x_1, 0)~\ge~0 \text{ and } \frac{\partial f_1}{\partial x_1}(x_1, 0) \le 0\; \forall x_1 \in [0,1]$), and with the rest of the sets $P_{ij}$ and $N_{ij}$ equal to the empty set.}
\end{minipage}
\vspace{-5mm}
\end{figure}
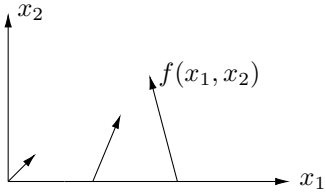

%%% Local Variables:
%%% mode: latex
%%% TeX-master: "../Learning_with_side_information_SIREV_v1.tex"
%%% fill-column: 61
%%% End:

    In the special case where $P_{ij} = \mathbb R^n$ and
    $N_{ij} = \emptyset$ for all $i,j \in \{1, \ldots, n\}$
    with $i \ne j$, and $P_{ii} = N_{ii} = \emptyset$ for all
    $i \in \{1, \ldots, n\}$, the side information is the
    knowledge of the following property of the vector field
    $f$ which appears frequently in the literature on
    monotone systems \cite{smith2008monotone}:
    \[ \forall x_{\text{init}}, \tilde x_{\text{init}} \in
      \mathbb R^n, \quad x_{\text{init}} \le \tilde
      x_{\text{init}} \implies x(t, x_{\text{init}}) \le x(t,
      \tilde x_{\text{init}}) \;\forall t \ge 0.\]
    Here, the inequalities are interpreted elementwise, and
    the notation $x(t, x_{\text{init}})$ is used  as before
    to denote the trajectory of the vector field $f$ starting
    from the initial condition~$x_{\text{init}}$.

\item {\bf Invariance of a set.} A set 
  $B \subseteq \Omega$ is  invariant under a vector field
  $f$ if any trajectory of the dynamical system
  $\dot \ix(t) = f(\ix(t))$ which starts in $B$ stays in $B$
  forever.\footnote{\bachir{In some other texts, this property is referred to as forward invariance  (to be contrasted with forward-and-backward invariance).}}  In particular, if
  $B = \{ x \in \R^n \; | \; h_j(\ix) \ge 0, \; j =
  1,,\dots,m\}$ for some differentiable
  functions $h_j: \mathbb R^n \rightarrow \mathbb R$, then
  invariance of the set~$B$ under the vector field~$f$
  implies the following constraints (see \Cref{fig:invariance} for an illustration): 
  \begin{equation}
    \label{eq:eq-invariance}
    \forall j \in \{1, \ldots, m\},\;\forall x \in B, \quad \left[h_j(x) = 0 \implies \langle f(\ix), \nabla h_j(\ix) \rangle \ge 0\right].
  \end{equation}

  Indeed, suppose for some $\tilde x \in B$ and
  for some $j \in \{1, \ldots, m\}$, we had $h_j(\tilde x) = 0$ but
  $\langle f(\tilde x), \nabla h_j(\tilde x) \rangle = \dot
  h(\tilde x) < 0$, then $h(x(t, \tilde x)) < 0$ for $t$
  small enough, implying that $x(t, \tilde x) \not \in B$ for
  $t$ small enough. It is also straightforward to verify that
  if the ``$\ge$'' in \eqref{eq:eq-invariance} were replaced
  with a ``$>$'', then the resulting condition would be
  sufficient for invariance of the set $B$ under $f$.  In
  fact, it follows from a theorem of Nagumo
  \cite{nagumo_1942,blanchini_invariance_1999} that condition
  \eqref{eq:eq-invariance} is necessary and sufficient for
  invariance of the set $B$ under $f$ if $B$ is convex,
  the functions $h_1, \ldots, h_m$ are continuously-differentiable, and for
  every point $x$ on the boundary of $B$, the vectors
  $\{\nabla h_j(x) \; |\; j \in \{1, \ldots, m\}, h_j(x) =
  0\}$ are linearly independent.\footnote{In an earlier draft
    of this paper \cite{learning_with_sideinfo_short}, we had
    incorrectly claimed that condition
    \cref{eq:eq-invariance} is equivalent to invariance of
    the set $B$, while in fact additional assumptions are
    needed for its sufficiency.}
    % referenece for cotingent cone: https://link.springer.com/content/pdf/10.1007/BF02193096.pdf
    Given sets
    $B_i = \{ x \in \R^n \; | \; h_{ij}(\ix) \ge 0, \; j =
    1,,\dots,m_i \}$, $i=1,\ldots, r$, defined by
    differentiable functions
    $h_{ij}: \mathbb R^n \rightarrow \mathbb R$, we denote by
    $\INV{\{B_i\}_{i=1}^r}$ the set of all vector fields
    {$f \in \CD$} that satisfy \eqref{eq:eq-invariance}
    for $B \in \{B_1, \ldots, B_r\}$.
   
%   \fixinv{Condition
%     \eqref{eq:eq-invariance} is also sufficient in the case
%     where $B$ is a polyhedron, an ellipsoid, or a Lorenz
%     cones. \cite{horvath_invariance_2016}}  

  % \fixinv{ We note that any vector
  %   field that satisfies the strict version of
  %   \eqref{eq:eq-invariance}, namely
  % \begin{equation}
  %   \label{eq:eq-invariance-strict}
  %   \forall i \in \{1, \ldots, m\},\;\forall x \in B, \quad [h_i(x) = 0 \implies \langle f(\ix), \nabla h_i(\ix) \rangle > 0],
  % \end{equation}    
  % leaves the set $B$ invariant.}

  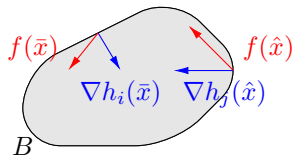
\begin{figure}[H]
%todo: violate montonocity cdts for which N or P is empty
\begin{minipage}[c]{0.41\textwidth}
\centering
\begin{tikzpicture}[scale=.2]
\draw[rounded corners=.5cm] (0,0) -- (0,5) -- (10,10) -- (15, 5) -- (10,0) -- cycle node[text=black] {$B$};
\draw[rounded corners=.5cm, fill=gray!20] (0,0) -- (0,5) -- (10,10) -- (15, 5) -- (10,0) -- cycle;
\draw[blue, -{Latex[length=2mm, width=1mm]}] (5, 7.5) -- (6.5, 5)node[below]{$\nabla h_i(\bar{x})$};
\draw[red, -{Latex[length=2mm, width=1mm]}] (5, 7.5) -- (3, 5)node[above left]{$f(\bar{x})$};
\draw[blue, -{Latex[length=2mm, width=1mm]}] (14, 5) -- (10, 5)node[below right]{$\nabla h_j(\hat x)$};
\draw[red, -{Latex[length=2mm, width=1mm]}] (14, 5)node[above right]{$f(\hat x)$} -- (11, 8);
\end{tikzpicture}
\end{minipage}
% \captionsetup{margin={1.5cm, 0cm, 0cm, 0cm}}
\begin{minipage}[c]{0.51\textwidth}
\centering
\caption{\label{fig:invariance}An example of the behavior of
  a vector field $f: \mathbb R^2 \rightarrow \mathbb R^2$
  satisfying $\INV{\{ B\}}$, where
  \(B \coloneqq \{x  \in \R^2 \; | \; h_1(x) \ge 0, \ldots, h_m(x) \ge 0\}\) is the
  set shaded in gray.}
\end{minipage}
\end{figure}

%%% Local Variables:
%%% mode: latex
%%% TeX-master: "../Learning_with_side_information_SIREV_v1.tex"
%%% fill-column: 61
%%% End:

\item {\bf Gradient and Hamiltonian systems.}  A vector field
  $f \in \CD$ is said to be a \textit{gradient} vector field
  if there exists a differentiable, scalar-valued function
  $V:\R^n \rightarrow \R$ such that
  \begin{equation}
  f(\ix) = - \nabla V(\ix)\; \forall \ix \in \Omega.\label{eq:grad}
\end{equation}
  Typically, the function $V$ is interpreted as a notion of
  potential or energy associated with the dynamical system
  $\dot x(t) = f(x(t))$. Note that the value of the function
  $V$ decreases along the trajectories of this dynamical
  system. We denote by $\GRAD$ the subset of $\CD$ consisting of gradient vector
  fields.

  A vector field $f\in \CD$ over $n$ state variables $(x_1, \ldots, x_n)$ is said to be \textit{Hamiltonian} if $n$ is even and there
  exists a differentiable scalar-valued function
  $H: \R^n \rightarrow \R$ such that
  \[f_i(p, q) = - \frac{\partial H}{\partial q_i} (p, q),
    %\text{ and }
    f_{\frac{n}2+i}(p,q) = \frac{\partial H}{\partial p_i}
    (p, q), \; {\forall (p,q) \in \Omega}, \; \forall i \in \left\{1, \ldots, \frac n2\right\},\]
  where $p \coloneqq (x_1, \ldots, x_{\frac{n}2})^T$ and
  $q \coloneqq (\ix_{\frac{n}2+1}, \ldots, \ix_n)^T$. The
  states $\ip$ and $\iq$ are usually referred to as
  \textit{generalized momentum} and \textit{generalized
    position} respectively, following terminology from
  physics. Note that a Hamiltonian system conserves the
  quantity $H$ along its trajectories. We denote by
  $\HAMILTONIAN$ the subset of $\CD$ consisting of Hamiltonian vector
    fields. For related work on learning Hamiltonian
  systems,
  see~\cite{ahmadi2018control,greydanus2019hamiltonian}.
  
\end{itemize}

\section{Learning Polynomial Vector Fields Subject to Side Information}
\label{sec:learn-poly-vec}

% To completely define the optimization problem in
% \eqref{eq:program_fit_vector_field_with_side_info}, we still
% have to specify the function class~$\mathcal F$. Among the
% possible choices are reproducing kernel Hilbert spaces
% \cite{sindhwani2018learning,sindhwani2018learnstabilizable,cheng2015learn},
% trigonometric functions, and functions parameterized by
% neural networks
% \cite{neural_laypunov_nips_2019,greydanus2019hamiltonian}.
In this paper, we take the function class $\mathcal F$ in
\eqref{eq:program_fit_vector_field_with_side_info} to be the
set of polynomial vector fields of a given degree $d$, \bachir{i.e., polynomial vector fields whose monomials have degree at most $d$}. We denote this function class by
\[{\poly d }\coloneqq \{ p: \R^n \rightarrow \R^n \; | \; p_i
  \text{ is a (scalar-valued) polynomial of degree $d$ for
    $i=1,\ldots,n$} \}.\]
Furthermore, we assume that the set
$\Omega$ over which we would like to learn the unknown dynamical system, the sets $P_i,N_i$ in the definition of
$\POS{ \{(P_i, N_i)\}_{i=1}^n }$, the sets $P_{ij},N_{ij}$
in the definition of
$\MON{ \{P_{i j}, N_{i j}\}_{i, j=1}^n}$, and the
sets $B_i$ in the definition of $\INV{\{B_i\}_{i=1}^r}$ are all closed
  semialgebraic. We recall that a \emph{closed basic
  semialgebraic} set is a set of the form
\begin{equation}\label{eq:basic.semi.algebraic}
\SETALG \coloneqq \{x\in\mathbb{R}^n|\ g_i(x)\geq 0, i=1,\ldots,m  \},
\end{equation}
where $g_1,\ldots,g_m$ are (scalar-valued) polynomial
functions, and that a \emph{closed semialgebraic} set is a
finite union of closed basic semialgebraic sets.

Our choice of working with polynomial functions to describe the vector field and the sets that appear in the side information definitions are motivated by two reasons.
%Our choices are motivated by two reasons. 
The first is that polynomial functions are expressive enough
to represent or approximate a large family of functions and
sets that appear in applications. The second reason, which
shall be made clear shortly, is that because of some
connections between real algebra and semidefinite
optimization, several side information constraints that are
commonly available in practice can be imposed on polynomial
vector fields in a numerically tractable fashion.

With our aforementioned choices, the optimization problem in
\eqref{eq:program_fit_vector_field_with_side_info} has as
decision variables the coefficients of a candidate polynomial
vector field $p: \mathbb R^n \rightarrow \mathbb R^n$.  When
the notion of side information is restricted to the six types
presented in \Cref{sec:side-info}, and under the mild
assumptions that $\Omega$ is full dimensional (i.e, that it
contains an open set), the
constraints of
\eqref{eq:program_fit_vector_field_with_side_info} are of the
following two types:
\setlist{leftmargin=10mm}
\begin{enumerate}
\item[(i)] Affine constraints in the coefficients of $p$.
\item[(ii)] Constraints of the type
  \begin{equation}\label{eq:p>=0.on.Omega}
    q(x)\geq 0\  \forall x\in \SETALG,
  \end{equation}
  where $\SETALG$ is a given closed basic semialgebraic set
  of the form \eqref{eq:basic.semi.algebraic}, and $q$ is a
  (scalar-valued) polynomial whose coefficients depend
  affinely on the coefficients of the vector field $p$.
\end{enumerate}
% Archimedean
For example, membership to
$\INTERPOLATION{\{(x_i,y_i)\}_{i=1}^m}$,
$\SYM{G,\sigma,\rho}$, $\GRAD$, or
$\HAMILTONIAN$ can be enforced by affine
constraints,\footnote{To see why membership of a polynomial
  vector field $p$ to $\GRAD$ can be enforced by affine
  constraints (a similar argument works for membership to
  $\HAMILTONIAN$), observe that if there exists a
  continuously-differentiable function
  $V: \R^n \rightarrow \R$ such that $p(x) = -\nabla V(x)$
  for all $x$ in a full-dimensional set $\Omega$, then the
  function $V$ is necessarily a polynomial of degree equal to
  the degree of $p$ plus one. Furthermore, equality between
  two polynomial functions over a full-dimensional set can be
  enforced by equating their coefficients.}
% For the latter two types of
% side infromation, this is because a gradient
% (resp. Hamiltonian) polynomial vector field is necessarily
% the gradient (resp. Hamiltonian) of a polynomial function,
% equality between two polynomials is equivalently to equality
% between their coefficients.
% TODO:fix gradient 
while membership to $\POS{ \{(P_i, N_i)\}_{i=1}^n }$,
$\INV{\{B_i\}_{i=1}^r}$, or
$\MON{ \{(P_{ij}, N_{i j})\}_{i, j=1}^n}$ can be cast as
constraints of the type
\eqref{eq:p>=0.on.Omega}. Unfortunately, imposing the latter
type of constraints is NP-hard already when $q$ is a quartic
polynomial and $\SETALG=\mathbb{R}^n$, or when $q$ is
quadratic and $\SETALG$ is a polytope (see, e.g.,
\cite{murty1985some}).

An idea pioneered to a large extent by Lasserre
\cite{lasserre_moment} and Parrilo \cite{sdprelax} has been
to write algebraic sufficient conditions for
(\ref{eq:p>=0.on.Omega}) based on the concept of sum of
squares polynomials. We say that a polynomial $h$ is a
\emph{sum of squares} (sos) if it can be written as
$h=\sum_i q_i^2$ for some polynomials~$q_i$. Observe that if
we succeed in finding sos polynomials
$\sigma_0,\sigma_1,\ldots, \sigma_m$ such that the polynomial
identity
\begin{equation}\label{eq:putinar-certificate}
q(x)=\sigma_0(x)+\sigma_1(x)g_1(x) + \ldots + \sigma_m(x)g_m(x)
\end{equation}
holds (for all $x \in \mathbb R^n$), then, clearly, the constraint in
(\ref{eq:p>=0.on.Omega}) must be satisfied. When the degrees
of the sos polynomials $\sigma_i$ are bounded above by an integer
$r$, we refer to the identity in
\eqref{eq:putinar-certificate} as the \emph{$\text{degree-}r$
  sos certificate} of the constraint in \eqref{eq:p>=0.on.Omega}. %corresponding to the constraint in
%\eqref{eq:p>=0.on.Omega}.
Conversely, the following
celebrated result in algebraic geometry
\cite{putinar1993positive} states that if $g_1,\ldots,g_m$
satisfy the so-called ``Archimedean property'' (a condition
slightly stronger than compactness of the set $\SETALG$),
then positivity of $q$ on $\Lambda$ guarantees existence of a
$\text{degree-}r$ sos certificate for some integer $r$ large
enough.

\begin{theorem}[Putinar's
  Positivstellensatz~\cite{putinar1993positive}]\label{th:putinar}
  Let 
  \[\SETALG =\{x \in \mathbb{R}^n ~|~ g_1(x)\geq 0, \ldots, g_m(x)\geq 0\}\]
  and assume that the collection of polynomials $\{g_1,\ldots,g_m\}$ satisfies the Archimedean property, i.e., there exists a positive scalar $R$ such that
  \[R^2-\sum_{i=1}^n x_i^2=s_0(x)+s_1(x) g_1(x)+\ldots+s_m(x) g_m(x),\] 
  where $s_0,\ldots,s_m$ are sos polynomials.\footnote{If $\SETALG$ is known to be contained in a ball of radius $R$, one can add the redundant constraint $R^2 - \sum_{i=1}^n x_i^2 \ge 0$ to the description of $\SETALG$, and then the Archimedean property will be automatically satisfied.} For any polynomial $q$, if  $q(x) > 0 \; \forall~x~\in~\SETALG$, then \[q(x)=\sigma_0(x)+\sigma_1(x)g_1(x)+\ldots+\sigma_m(x)g_m(x),\]
  for some sos polynomials $\sigma_0,\ldots,\sigma_m$.
\end{theorem}

%:\mathbb{R}^n\rightarrow\mathbb{R}
\bachir{The computational appeal of the sum of squares approach stems
from its connection to \emph{semidefinite programming} (SDP).}
\bachir{We recall that semidefinite programming is the problem of minimizing a linear function of a symmetric matrix over the intersection of the cone of positive semidefinite matrices with an affine subspace. Semidefinite programs can be solved to arbitrary accuracy in time that scales polynomially with their input size; see
\cite{vandenberghe1996semidefinite} for a survey of the
theory and applications of this subject.}

\bachir{To make the connection between \cref{th:putinar} and SDP more clear, we remark that the search for sos polynomials
$\sigma_0,\sigma_1,\ldots, \sigma_m$ of a given degree that
verify the polynomial identity in
\eqref{eq:putinar-certificate} can be automated via
SDP}. This is true even when some
coefficients of the polynomial $q$ are left as decision
variables. This claim is a straightforward consequence of the
following well-known fact (see, e.g.,~\cite{PhD:Parrilo}): A
polynomial $h$ of degree $2d$ is a sum of squares if and only
if there exists a symmetric matrix $Q$ which is positive
semidefinite and verifies the identity
\begin{equation}
h(x)=z(x)^TQz(x),\label{eq:sos_to_sdp}
\end{equation}
where $z(x)$ denotes the vector of all monomials in $x$
of degree less than or equal to $d$.
\bachir{The size of the matrix $Q$ is ${n+d \choose d}$, which is  polynomial in $n$ (resp. $d$) if $d$ (resp. $n$) is fixed.
Identity
\eqref{eq:sos_to_sdp} can be written in an equivalent manner
as a system of ${n+2d \choose 2d}$ linear equations involving the entries of
the matrix $Q$ and the coefficients of the polynomial $h$.  These equations come from equating the
coefficients of the polynomials appearing on the left and
right hand sides of \eqref{eq:sos_to_sdp}. 
The problem of finding a positive semidefintie matrix $Q$ whose entries satisfy these linear equations is a semidefinite program.}
For implementation purposes, there exist modeling languages,
such as YALMIP~\cite{yalmip}, SOSTOOLS~\cite{sostools}, or
SumOfSquares.jl~\cite{sumofsquares_jl}, that accept sos
constraints on polynomials directly and do the conversion to
a semidefinite program in the background. See e.g.
\cite{laurent_survey_2009,PabloGregRekha_BOOK,hall_sos_engineering_2019}
for more background on sum of squares techniques.

To end up with a semidefinite programming formulation of
problem \eqref{eq:program_fit_vector_field_with_side_info},
we also need to take the loss function $\ell$ that appears in
the objective function to be semidefinite representable
(i.e., we need its epigraph to be the projection of the
feasible set of a semidefinite program; \bachir{see \cite[Chapter 3.2]{bental} for a discussion on semidefinite representability and several examples}).  Luckily, many
common loss functions in machine learning are semidefinite
representable. Examples of such loss functions include
\begin{inparaenum}[(i)]
    \item any $\ell_p$ norm for a rational number $p \ge 1$,
      or for $p = \infty$,
    \item any convex piece-wise linear function,
    \item any $\emph{sos-convex}$ polynomial (see
      e.g.~\cite{helton_sosconvexity_2010} for a definition), and
    \item any positive integer power of the previous three
      function classes.
\end{inparaenum}

\section{Illustrative Experiments}
\label{sec:applications}
In this section, we present numerical experiments from four
application domains to illustrate our methodology. The first
three applications are learning experiments and the last one
involves an optimal control component. In all of our
experiments, we use the SDP-based approach explained in
\Cref{sec:learn-poly-vec} to tackle problem
\eqref{eq:program_fit_vector_field_with_side_info} and take our loss
function $\ell(\cdot, \cdot)$ in
\eqref{eq:program_fit_vector_field_with_side_info} to be
$\ell(u, v) = \|u - v\|_2^2 \quad \forall u,v\in\mathbb  R^n.$
The added value of side information for
learning dynamical systems from data will be demonstrated in
these experiments.

\subsection{Diffusion of a contagious disease}
\label{sec:diffusion_contagious_disease}

\begin{figure}
  \begin{minipage}[c]{0.4\textwidth}
  \includegraphics[width=.8\textwidth, trim=0 10 0 0,
  clip]{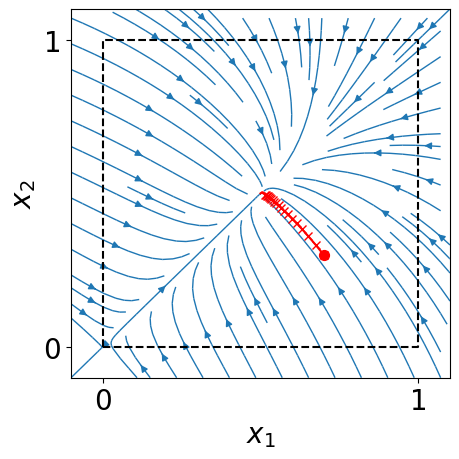}
\end{minipage}
 \begin{minipage}[c]{0.5\textwidth}
   \caption{\label{fig:vf_diffusion_disease} Streamplot of
     the vector field in \eqref{eq:diffusion_disease}. We
     consider this vector field to be the ground truth and
     unknown to us. We would like to learn it over
     $[0, 1]^2$ from noisy snapshots of a
     single trajectory starting from $(0.7, 0.3)^T$ (plotted
     inred).}
  \end{minipage}
\end{figure}

The following dynamical system has appeared in the epidemiology literature (see, e.g., \cite{anderson1992infectious}) as a model for the spread of a sexually transmitted disease in a heterosexual population:
\begin{equation}  \label{eq:diffusion_disease}
    \dot \ix(t) = f(\ix(t)), \text{where} \; \ix(t) \in \R^2 \text{ and }
  f(\ix) = \left(
    \begin{array}{l}
      -a_1   x_1 + b_1 (1-x_1) x_2\\
      -a_2   x_2 + b_2 (1-x_2) x_1
    \end{array}\right).
\end{equation}
Here, the quantity $x_1(t)$ (resp. $x_2(t)$) represents the
fraction of infected males (resp. females) in the
population. The parameter $a_i$ (resp. $b_i$)
denotes the recovery rate (resp. the infection rate) in the male population when $i=1$,
and in the female population when $i=2$. We take
\begin{equation}
    (a_1, b_1, a_2, b_2) = (0.05 , 0.1, 0.05, 0.1) 
\label{eq:disease_params}
\end{equation}
and plot the resulting vector field $f$ in
\cref{fig:vf_diffusion_disease}. We suppose that this vector
field is \emph{unknown} to us, and our goal is to learn it
over $\Omega \coloneqq [0, 1]^2$ from a few noisy snapshots
of a single trajectory. More specifically, we have access to
the training data set
\begin{equation}\label{eq:training_data_epi}
\mathcal D \coloneqq \left\{ \left(\ix(t_i, x_{\text{init}}),
      f(\ix(t_i, x_{\text{init}}))+
      10^{-4} \begin{pmatrix}\varepsilon_{i,1}\\\varepsilon_{i,2}\end{pmatrix}\right)
  \right\}_{i=1}^{20},
  \end{equation}
where $\ix(t, x_{\text{init}})$ is the
trajectory of the system
\eqref{eq:diffusion_disease} starting from the initial
condition 
\[\ix_{\text{init}} = (0.7, 0.3)^T,\]
the scalars
$t_i \coloneqq i$ represent a uniform subdivision of the
time interval $[0, 20]$, and the scalars
$\varepsilon_{1,1}, \varepsilon_{1,2}, \ldots, \varepsilon_{20,1}, \varepsilon_{20,2}$ are independently sampled from the
standard normal distribution.

Following our approach in \Cref{sec:learn-poly-vec}, we
parameterize our candidate vector
field \[p:~\R^2~\rightarrow~\R^2\] as a polynomial
function. We choose the degree of this polynomial to be $d =
3$. The degree $d$ is taken to be larger than $2$ because we
do not want to assume knowledge of the degree of the true
vector field in \eqref{eq:diffusion_disease}. This makes the
task of learning more difficult; see the end of this
subsection where we also learn a vector field of degree $2$
for comparison.

%  by taking the degree to be $3$, and we will
% show that more side information is needed to prevent
% overfitting to noise in that case. The solution to problem
% \eqref{eq:disease-fitting-polynomial-lsq} corresponding to
% $d=3$ is plotted in
% \cref{fig:vf_diffusion_disease_3_no_side_info}. Unsurprisingly,
% this solutions overfits to the observed trajectory but
% behaves very differently from $f$ on the rest of the unit
% box. We resume our experiements by imposing the following
% additional side information constrats.

% In this experiment, we pretend that the true vector field $f$
% is unknown to us. 
%  In this experiment, we pretend that $f$ is
% unknown to us and consider an over-parameterized model of the
% true dynamics by taking $d=3$.
In absence of any side
information, one could solve the least-squares
problem
\begin{equation}
  \min_{p \in \mathcal P_d} \sum_{(\ix_i, \iy_i) \in \mathcal D} \| p(\ix_i) - \iy_i\|_2^2
  \label{eq:disease-fitting-polynomial-lsq}
\end{equation}
to find a polynomial of degree $d$ that best agrees with the
training data. For this experiment only, and for educational purposes, we include a template code using the library
SumOfSquares.jl~\cite{sumofsquares_jl} of the Julia
programming language and demonstrate how the code changes as we impose side information constraints. We initiate our template with the following code that solves optimization problem \cref{eq:disease-fitting-polynomial-lsq}. %(i.e., when there is no side information) is as follows:

% We present the code responsible for defining
% and solving this optimization problem using the library
% SumOfSquares.jl~\cite{sumofsquares_jl} of the Julia
% programming language. %but there are other choices such as
% % YALMIP~\cite{yalmip} and SOSTOOLS~\cite{sostools}.

\renewcommand{\lstlistingname}{Template Code}
\vspace{1em}
\begin{lstlisting}[style=JuliaListing,title={Julia template code for learning dynamical systems with side information.}]
%\CodeComment[13cm]{Input: vectors $X_1$, $X_2$, $Y_1$, $Y_2$ $\in \mathbb R^{20}$ representing the training set in  \eqref{eq:training_data_epi},  }% 
%\CodeComment[13cm]{\quad \quad \; where $X_1 = \{ x_1(t_i, x_{\text{init}})\}_{i=1}^{20}$, $X_2 = \{ x_2(t_i, x_{\text{init}})\}_{i=1}^{20}$,     }% 
%\CodeComment[13cm]{\quad \quad \; $Y_1 = \{ f_1(x(t_i, x_{\text{init}})) + 10^{-4}\varepsilon_{i,1}\}_{i=1}^{20}$, and $Y_2 = \{ f_2(x(t_i, x_{\text{init}})) + 10^{-4}\varepsilon_{i,2}\}_{i=1}^{20}$}%

model = SOSModel(solver)% \CodeComment{solver could be any SDP solver,}% 
% \CodeComment{e.g., Mosek \cite{mosek}, SDPT3 \cite{SDPT3}, CSDP \cite{csdp}}% 
@polyvar $x_1$ $x_2$% \CodeComment{Define state variables $x_1$, $x_2$}%
d = 3                  % \CodeComment{Construct vector of monomials}% 
z = monomials([$x_1$, $x_2$], 0:d)% \CodeComment{in ($x_1$, $x_2$) up to degree $d$}%
@variable(model, $p_1$, Poly(z))% \CodeComment{Declare a polynomial vector field}%
@variable(model, $p_2$, Poly(z))% \CodeComment{whose coefficients are decision variables}%
error_vec = [p[1].($X_1$, $X_2$) - $Y_1$;%\CodeComment{Vector of individual terms appearing}%
             p[2].($X_1$, $X_2$) - $Y_2$]%\CodeComment{in the objective of \cref{eq:disease-fitting-polynomial-lsq}}%
@objective model Min error_vec$'$ * error_vec
# Side information constraints go here
# ...
optimize!(model)%\CodeComment{Solve the optimization problem}%

\end{lstlisting}

% The solution to problem
% \eqref{eq:disease-fitting-polynomial-lsq} returned by the solver MOSEK \cite{mosek} is plotted in
% \cref{fig:vf_diffusion_disease_2_no_side_info} and is given
% by\footnote{The empty set in the subscript of $p_{\emptyset,\text{deg }2}$ signifies the fact that the list of side information is empty.}% the quadratic polynomial
% \[ p_{\emptyset,\text{deg }2}(x_1, x_2) = \begin{pmatrix}
%   -1.28x_1^2 - 2.30x_1x_2 - 0.99x_2^2 + 2.33x_1 + 2.18x_2 - 1.11\\
%   10.29x_1^2 + 17.13x_1x_2 + 7.74x_2^2 - 18.80x_1 - 16.25x_2 + 8.73
%  \end{pmatrix}.
% \]

\begin{figure}
  \centering
%   \begin{subfigure}[t]{0.27\textwidth}
%     \includegraphics[width=.9\textwidth, trim=0 10 0 0, clip]{imgs/disease/disease_true_vf.png}
%     \captionsetup{width=.8\linewidth}
%     \hspace{-1cm}\caption{\label{fig:disease_true_vf_2}The true vector field (unknown to us)}
%   \end{subfigure}
  
  \advance\leftskip-1.5cm
  \begin{subfigure}[t]{0.27\textwidth}
    \includegraphics[width=.9\textwidth, trim=0 10 0 0, clip]{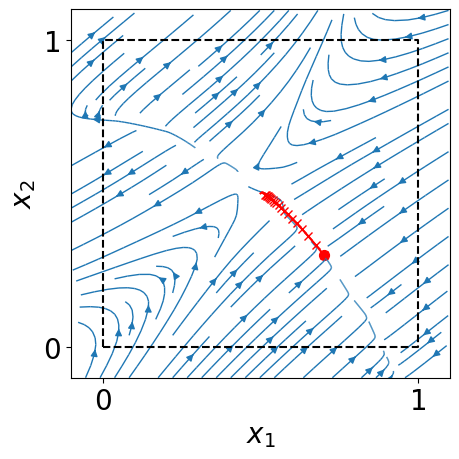}
    \caption{\label{fig:vf_diffusion_disease_3_no_side_info}No side information}
  \end{subfigure}~
  \begin{subfigure}[t]{0.27\textwidth}
    \includegraphics[width=.9\textwidth,
    trim=0 10 0 0, clip]{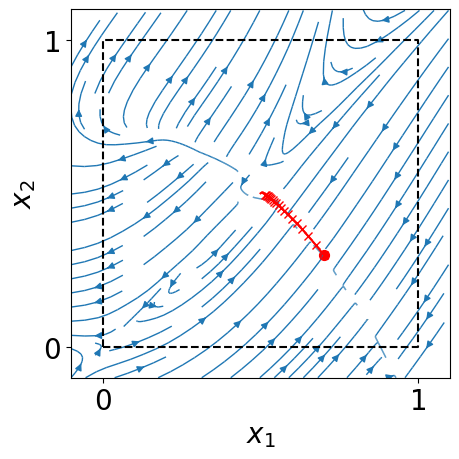}
  \caption{\label{fig:vf_diffusion_disease_3_eq}$\NINTERPOLATION$}
\end{subfigure}~
\begin{subfigure}[t]{0.27\textwidth}  
  \includegraphics[width=.9\textwidth,
  trim=0 10 0 0,
  clip]{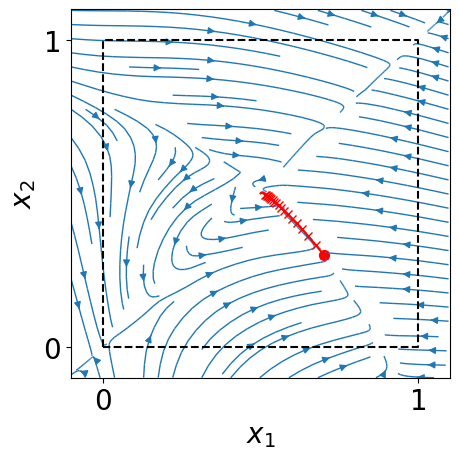}
\caption{\label{fig:vf_diffusion_disease_3_eq_inv}$\NINTERPOLATION \cap \NINV$}
\end{subfigure}~
\begin{subfigure}[t]{0.27\textwidth}  
  \includegraphics[width=.9\textwidth,
  trim=0 10 0 0,
  clip]{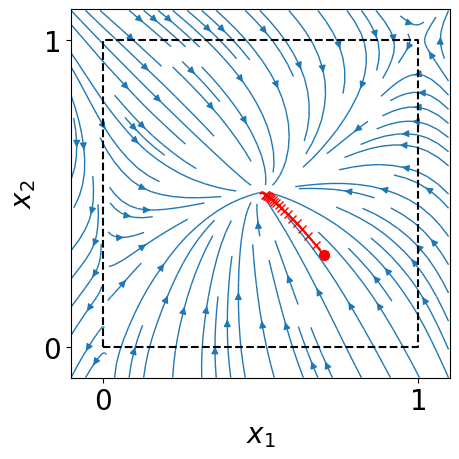}
\caption{\label{fig:vf_diffusion_disease_3_eq_inv_mon}$\NINTERPOLATION \cap \NINV \cap \NMON$}
\end{subfigure}

 \caption{ 
   \label{fig:disease-lsq-deg-3}
   % (\cref{fig:disease_true_vf_2})  Streamplot of the vector field in \eqref{eq:diffusion_disease} along with 
   \bachir{
   Streamplots} of polynomial vector fields of degree $3$
   that are optimal to \eqref{eq:disease-fitting-polynomial-lsq} with different side information constraints appended to it. %(\cref{fig:vf_diffusion_disease_3_no_side_info,fig:vf_diffusion_disease_3_eq,fig:vf_diffusion_disease_3_eq_inv,fig:vf_diffusion_disease_3_eq_inv_mon}). 
   %In each case, the trajectory of the learned
   %vector field starting from $(0.7, 0.3)^T$ is also
   %plotted.
 }
 \vspace{-5mm}
\end{figure}

The solution to problem
\eqref{eq:disease-fitting-polynomial-lsq} returned by the solver MOSEK \cite{mosek} is plotted in
\cref{fig:vf_diffusion_disease_3_no_side_info}.
Observe that
while the learned vector field replicates the behavior of the
true vector field $f$ on the observed trajectory, it differs
significantly from $f$ on the rest of the unit square. To remedy
this problem, we leverage the following list of side information that
is available from the context without knowing the exact
structure of $f$.

\vspace{1em}
\begin{itemize}
  \setlength\itemsep{1em}
\item {\bf Equilibrium point at the origin
    ($\NINTERPOLATION$).} Naturally, if no
  male or female is infected, there would be no contagion and the number of infected individuals will remain at zero. This side information
  corresponds to our vector field $p$ having an equilibrium
  point at the origin, i.e., $p(0, 0) = 0$.  Note from \bachir{\cref{fig:vf_diffusion_disease,fig:vf_diffusion_disease_3_no_side_info}} that the true vector field $f$ in \eqref{eq:diffusion_disease} satisfies this constraint, but the vector field learned by solving the least-squares problem in \eqref{eq:disease-fitting-polynomial-lsq} does not. We can impose
  this linear constraint by simply adding the following lines of code to our template:
  %Note from
  %\cref{fig:vf_diffusion_disease_3_no_side_info} that the
  %least-squares solution does not satisfy this side
  %information.
  %The code for this constraint is given below.
  \begin{lstlisting}[style=JuliaListing]
    @constraint model %$p_1(0, 0)$% == 0
    @constraint model %$p_2(0, 0)$% == 0
  \end{lstlisting}
  The vector field resulting from solving this new problem is plotted in
  \cref{fig:vf_diffusion_disease_3_eq}. %Interestingly, its behavior is already quite close to that of the true vector field in \cref{fig:vf_diffusion_disease}.

\item {\bf Invariance of the box $[0, 1]^2$ ($\NINV$).} The
  state variables $(x_1, x_2)$ of the dynamics in
  \eqref{eq:diffusion_disease} represent fractions of
  infected individuals and as such, the vector $\ix(t)$
  should be contained in the box $[0, 1]^2$ at all times
  $t \ge 0$. Note that this property is violated by the
  vector fields learned in
  \cref{fig:vf_diffusion_disease_3_no_side_info,fig:vf_diffusion_disease_3_eq}.
  Mathematically, the invariance of the unit box is
  equivalent to the four (univariate) polynomial
  nonnegativity constraints
  \[p_2(x_1, 0) \ge 0, p_2(x_1, 1) \le 0 \;\forall x_1 \in [0, 1],\]
  \[ p_1(0, x_2) \ge 0, p_1(1, x_2) \le 0 \; \forall x_2 \in
    [0, 1].\]
  These constraints imply that the vector field points inwards on the
  four edges of the unit box. We replace each one of these
  four constraints with the corresponding $\text{degree-}2$
  sos certificate of the type in
  \eqref{eq:putinar-certificate}.
  For instance, we replace the constraint
  \[p_2(x_1, 0) \ge 0 \; \forall x_1 \in [0, 1]\] with
  linear and semidefinite constraints obtained from equating the coefficients
  of the two sides of the polynomial identity
  \begin{equation}
  p_2(x_1, 0) = x_1s_0(x_1) + (1-x_1) s_1(x_1),\label{eq:poly_identity_putinar_coord_nonneg} 
\end{equation}
  and requiring that the newly-introduced (univariate) polynomials $s_0$ and $s_1$  be quadratic and sos. Obviously, the algebraic identity \eqref{eq:poly_identity_putinar_coord_nonneg} is sufficient
  for nonnegativity of $p_2(x_1, 0)$ over $[0, 1]$; in this
  case, it also happens to be necessary
  \cite{univariateposinterval}. The code in Julia for imposing the $\text{degree-}2$ sos certificate in \eqref{eq:poly_identity_putinar_coord_nonneg} is as follows:
  \begin{lstlisting}[style=JuliaListing]
    @variable(model, %$s_1$%, Poly([%$1, x_1, x_1^2$%]))% \CodeComment[4.0cm]{Declare decision polynomial $s_1$}%
    @variable(model, %$s_2$%, Poly([%$1, x_1, x_1^2$%]))% \CodeComment[4.0cm]{Declare decision polynomial $s_2$}%
    @constraint(model, %$s_1$%, in SOSCone())%\CodeComment[4.0cm]{Enforce $s_1$ to be sos}%
    @constraint(model, %$s_2$%, in SOSCone())%\CodeComment[4.0cm]{Enforce $s_2$ to be sos}%

    polynomial_identity = %\CodeComment[4.0cm]{Enfroce polynomial}%
    % $\quad p_2((x_1, x_2) \Rightarrow  (0, x_2))$% - % $x_1$%*%$s_1$%-(1-%$x_1$%)*%$s_2$% %\CodeComment[4.0cm]{identity in \cref{eq:poly_identity_putinar_coord_nonneg}}%

    @constraint(model,coefficients(polynomial_identity).== 0)
  \end{lstlisting}

  The output of the semidefinite program which
  imposes the invariance of the unit box and the equilibrium point
  at the origin is plotted in
  \cref{fig:vf_diffusion_disease_3_eq_inv}.
  
\item {\bf Coordinate directional monotonicity ($\NMON$).}
  Naturally, one would expect that if the fraction of infected males rises in
  the population, the rate of infection of females should
  increase. Mathematically, this observation is equivalent to the constraint that
  \[\frac {\partial p_2}{\partial x_1}(\ix) \ge 0 \;
  \forall \ix \in [0, 1]^2.\]
  Similarly, swapping the roles played by males and
  females leads to the constraint
  \[\frac {\partial p_1}{\partial x_2}(\ix) \ge 0 \; \forall
  \ix \in [0, 1]^2.\]
 Note that this property is violated by the vector fields learned in \cref{fig:vf_diffusion_disease_3_no_side_info,fig:vf_diffusion_disease_3_eq,fig:vf_diffusion_disease_3_eq_inv}.
 Just as in the previous bullet point, we replace each of the above two nonnegativity constraints with its corresponding $\text{degree-}2$ sos certificate. To do this, we represent the closed basic semialgebraic set $[0, 1]^2$ with the polynomial inequalities
 \[x_1 \ge 0, x_2 \ge 0, 1-x_1 \ge 0, 1-x_2 \ge 
 0.\]
%   Since $[0, 1]^2$ is a closed basic semialgebraic set,
%   in the same spirit as the previous bullet point, we
%   can replace each of the above constraints with its
%   corresponding sos certificate. We choose the sos certificate to be $2$. 
The Julia code for imposing this side information is similar
to the one of the previous bullet point and therefore
omitted.  \cref{fig:vf_diffusion_disease_3_eq_inv_mon}
demonstrates the vector field learned by our semidefinite
program when all side information constraints discussed thus
far are imposed.
  \end{itemize}
% \vspace{-\topsep}

Note from
\cref{fig:vf_diffusion_disease_3_no_side_info,fig:vf_diffusion_disease_3_eq,fig:vf_diffusion_disease_3_eq_inv,fig:vf_diffusion_disease_3_eq_inv_mon}
that as we add more side information, the learned vector
field respects more and more properties of the true vector
field $f$. In particular, the learned vector field in
\cref{fig:vf_diffusion_disease_3_eq_inv_mon} is quite similar
qualitatively to the true vector field 
\bachir{in \cref{fig:vf_diffusion_disease}}
even though only noisy snapshots of a single trajectory are used for
learning. 

\begin{figure}
  
  \centering
  
%     \begin{subfigure}[t]{0.27\textwidth}
%     \includegraphics[width=.9\textwidth, trim=0 10 0 0, clip]{imgs/disease/disease_true_vf.png}
%     \captionsetup{width=.8\linewidth}
%     \hspace{-1cm}\caption{The true vector field (unknown to us)}
%   \end{subfigure}
  
  \advance\leftskip-1.5cm
  \begin{subfigure}[t]{.27\textwidth}
    \includegraphics[width=.9\textwidth, trim=0 10 0 0, clip]{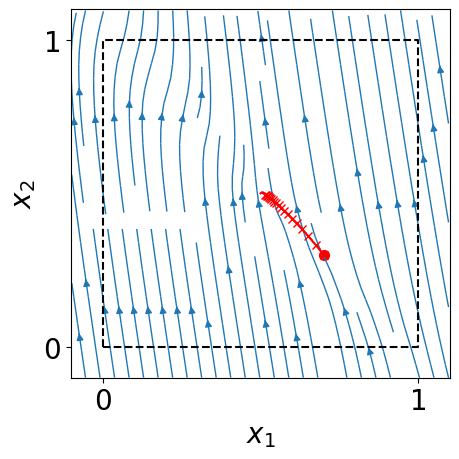}
    \caption{\label{fig:vf_diffusion_disease_2_no_side_info}No side information}
  \end{subfigure}~
  \begin{subfigure}[t]{.27\textwidth}
    \includegraphics[width=.9\textwidth,
    trim=0 10 0 0, clip]{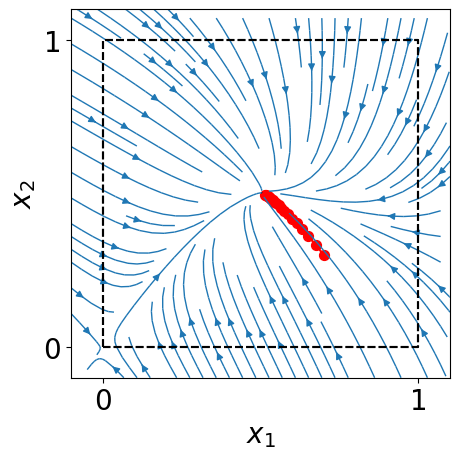}
  \caption{\label{fig:vf_diffusion_disease_2_eq}$\NINTERPOLATION$}
\end{subfigure}~
\begin{subfigure}[t]{.27\textwidth}  
  \includegraphics[width=.9\textwidth,
  trim=0 10 0 0,
  clip]{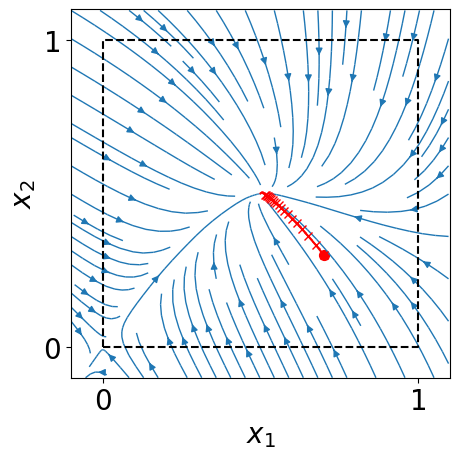}
\caption{\label{fig:vf_diffusion_disease_2_eq_inv}$\NINTERPOLATION \cap \NINV$}
\end{subfigure}~
\begin{subfigure}[t]{.27\textwidth}  
  \includegraphics[width=.9\textwidth,
  trim=0 10 0 0,
  clip]{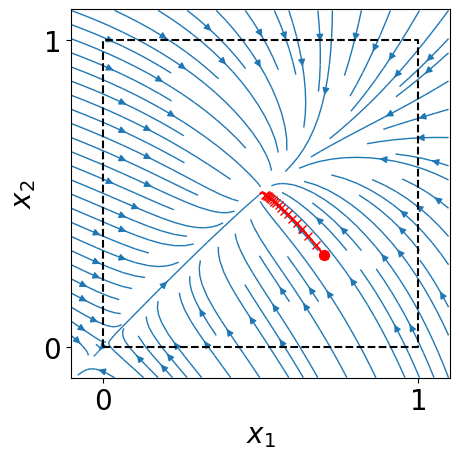}
\caption{\label{fig:vf_diffusion_disease_2_eq_inv_mon}$\NINTERPOLATION \cap \NINV \cap \NMON$}
\end{subfigure}

 \caption{
   \label{fig:disease-lsq-deg-2} 
   %Streamplot of the vector field in \eqref{eq:diffusion_disease} (\cref{fig:disease_true_vf_2}) along with
   %
   \bachir{Streamplots} of polynomial vector fields of degree $2$
   that are optimal to \eqref{eq:disease-fitting-polynomial-lsq} with different side information constraints appended to it. %(\cref{fig:vf_diffusion_disease_2_no_side_info,fig:vf_diffusion_disease_2_eq,fig:vf_diffusion_disease_2_eq_inv,fig:vf_diffusion_disease_2_eq_inv_mon}).
 }
\end{figure}

It is interesting to observe what would happen if we try to learn a $\text{degree-}2$ vector field from the same training set using the list of side information discussed in this subsection. The outcome of this experiment is plotted in \cref{fig:disease-lsq-deg-2}.
Note that
with the equilibrium-at-the-origin side information, the
behavior of the learned vector field of degree $2$
is already quite close to that of the true vector field.
\cref{fig:vf_diffusion_disease_2_eq_inv_mon} shows that when we impose all side information, the learned vector field is almost indistinguishable from the true vector field (even though, once again, only noisy snapshots of a single trajectory are used for
learning). The vector field plotted in \cref{fig:vf_diffusion_disease_2_eq_inv_mon} is given by 
\[\small
  p_{{\NINTERPOLATION \cap \NINV \cap \NMON, \text{deg}\;2}}(x_1, x_2) = \begin{pmatrix}
0.038x_1^2 - 0.100 x_1x_2 - 0.009x_2^2 - 0.084x_1 + 0.119x_2\\
 -0.101x_1x_2 + 0.003x_2^2 + 0.101x_1 - 0.052x_2   
 \end{pmatrix},
\]
which is indeed very close to the vector field in \eqref{eq:diffusion_disease}. \bachir{For experiments with other degrees and noise levels, see the appendix.}

% \begin{figure}[H]
%   \centering
%   \begin{subfigure}[t]{.3\textwidth}
%     \includegraphics[width=.9\textwidth, trim=0 10 0 0,
%     clip]{imgs/disease/disease_true_vf.png}
%     \caption{\label{fig:vf_diffusion_disease_true_vf_comparison}True and unknown vector field in
%       \eqref{eq:diffusion_disease}.}
%     \end{subfigure}\quad
%     \begin{subfigure}[t]{.3\textwidth}
%     \includegraphics[width=.9\textwidth, trim=0 10 0 0, clip]{imgs/disease/deg_2_eq_invariance_monotone.png}
%     \caption{\label{fig:vf_diffusion_disease_2_final}$\NINTERPOLATION \cap \NINV \cap \NMON$,\\degree $2$.}
%     \end{subfigure}
%     \begin{subfigure}[t]{.3\textwidth}
%     \includegraphics[width=.9\textwidth, trim=0 10 0 0, clip]{imgs/disease/deg_3_eq_invariance_monotone.png}
%     \caption{$\NINTERPOLATION \cap \NINV \cap \NMON$,\\degree $3$.}
%     \end{subfigure}
%  \caption{
%   \label{fig:disease-final-comparison} Comparison of the
%   streamplots of the true (unknown) vector field $f$
%   (\cref{fig:vf_diffusion_disease_true_vf_comparison}) and the
%   polynomial vector fields of degree $2$
%   (\cref{fig:vf_diffusion_disease_2_eq_inv_mon}) and $3$
%   (\cref{fig:vf_diffusion_disease_3_eq_inv_mon}) returned by
%   our semidefinite programs with all side information constraints added.}
% \end{figure}

\subsection{Dynamics of the simple pendulum}
\label{sec:simple-pendulum}

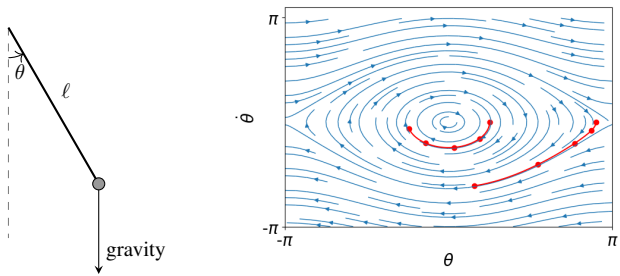
\begin{figure}
  \centering
  \begin{tikzpicture}[scale=.8, transform shape]
    % save length of g-vector and theta to macros
    \pgfmathsetmacro{\Gvec}{1.5}
    \pgfmathsetmacro{\myAngle}{30}
    % calculate lengths of vector components
    \pgfmathsetmacro{\Gcos}{\Gvec*cos(\myAngle)}
    \pgfmathsetmacro{\Gsin}{\Gvec*sin(\myAngle)}

    \coordinate (centro) at (0,0);
    \draw[dashed,gray,-] (centro) -- ++ (0,-3.5) node (mary)
    [black,below]{$ $};
    \draw[thick] (centro) -- ++(270+\myAngle:3) coordinate (bob) node
    [black, midway, above right]{$\ell$};
    \pic [draw, ->, "$\theta$", angle eccentricity=1.5] {angle = mary--centro--bob};
    % \draw [blue,-stealth] (bob) -- ($(bob)!\Gcos cm!(centro)$) node
    % [midway, below left]{$\vec T$};
    % \draw [-stealth] (bob) -- ($(bob)!-\Gcos cm!(centro)$)
    %   coordinate (gcos)
    %   node[midway,above right] {$a\cos\theta$};
    % \draw [-stealth] (bob) -- ($(bob)!\Gsin cm!90:(centro)$)
    %   coordinate (gsin)
    %   node[midway,above left] {$a\sin\theta$};
    \draw [-stealth] (bob) -- ++(0,-\Gvec)
      coordinate (g)
      node[near end,right] {gravity};
    % \pic [draw, ->, "$\theta$", angle eccentricity=1.5] {angle = g--bob--gcos};
    \filldraw [fill=black!40,draw=black] (bob) circle[radius=0.1];
\end{tikzpicture}\quad\quad
  \includegraphics[width=.4\textwidth]{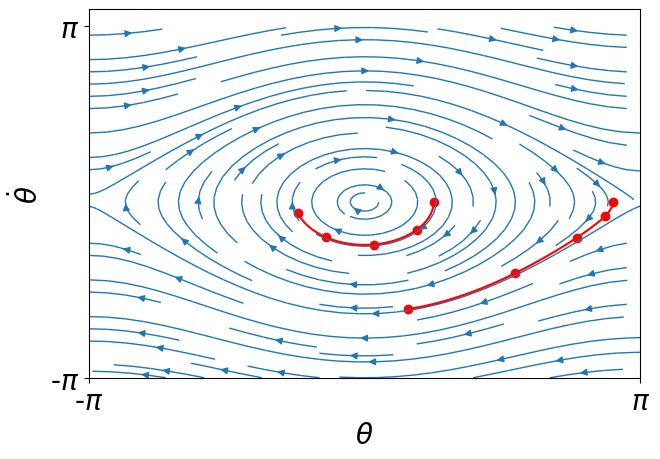}%\hspace*{5in}
  \caption{\label{fig:fig-simple-pendulum}The simple pendulum
    and the streamplot of its vector field. We would like to
    learn this vector field over $[-\pi, \pi]^2$ from 10
    noisy snapshots coming from two trajectories. }
  \vspace{-5mm}
\end{figure}

In this subsection, we consider the dynamics of the simple pendulum,
i.e., a mass $m$ hanging from a massless rod of length
$\ell$ (see \cref{fig:fig-simple-pendulum}). The state variables
of this system are given by $x = (\theta, \dot \theta)$,
where $\theta$ is the angle that the rod makes with the
vertical axis and $\dot \theta$ is the time derivative of
this angle. By convention, the angle $\theta \in [-\pi, \pi]$
is positive when the mass is to the right of the vertical
axis, and negative otherwise.  By applying Newton's second
law of motion, the equation
\[\ddot \theta(t) = -\frac g \ell \sin (\theta(t))\]
for the dynamics of the pendulum
can be derived, where $g$ here is the acceleration due to
gravity. This is a one-dimensional second-order system that
we convert to a first-order system as follows:
\begin{equation}
  \dot x(t) = f(x(t)) \quad \text{where } \quad x(t) \coloneqq \begin{pmatrix} \theta(t)\\\dot \theta(t)\end{pmatrix} \text{ and }  f(\theta, \dot \theta) \coloneqq \begin{pmatrix}\dot \theta\\-\frac g\ell \sin \theta\end{pmatrix}.\label{eq:ode-pendulum}
\end{equation}

We take the vector field in \eqref{eq:ode-pendulum} with
$g=\ell=1$ to be the ground truth. We observe from this
vector field a noisy version of two trajectories
$x(t, x_{\text{init}}^1)$ and $x(t, x_{\text{init}}^2)$
sampled at times $t_i = \frac{3i}5$, where
$i\in\{0, \ldots, 4\}$, with
$x_{\text{init}}^1 = (\frac \pi4, 0)^T$ and
$x_{\text{init}}^2 = (\frac{9\pi}{10}, 0)^T$ (see
\cref{fig:fig-simple-pendulum}). More precisely, we assume
that our training set (with a slightly different
representation of its elements) is given by
\begin{equation}
  %\footnotesize
  \mathcal D \coloneqq
  \bigcup_{k=1}^2
  \left\{ \begin{pmatrix}\theta(t_i, x_{\text{init}}^k)\\ \dot
    \theta(t_i, x_{\text{init}}^k)\\ \ddot \theta(t_i,
    x_{\text{init}}^k)\end{pmatrix} + 10^{-2} \varepsilon_{i,k} \right\}_{i=0}^4,
  \label{eq:data-pendulum}
\end{equation}
where the $\varepsilon_{i,k}$ (for $k=1,2$ and $i=0,\ldots,4$)
are independent $3 \times 1$ standard normal vectors.

We are interested in learning the vector field $f$ over the
set $\Omega \coloneqq [-\pi, \pi]^2$ from the training data in
\eqref{eq:data-pendulum}. We parameterize our candidate vector field $p: \mathbb R^2 \rightarrow \mathbb R^2$
as a $\text{degree-}5$ polynomial. Note that
$p_1(\theta, \dot \theta) = \dot \theta$, just from the
meaning of our state variables. The only unknown is therefore the polynomial
$p_2(\theta, \dot \theta)$. 

In absence of side information, one can solve a least-squares
problem that finds a polynomial of degree $5$ that best
agrees with the training data. As it can be seen in
\cref{fig:pendulum_no_side_info}, the resulting vector field
is very far from the true vector field and is unable to even
replicate the observed trajectories when started from the
same two initial conditions. To learn a better model, we
describe a list of side information which could be derived
from contextual knowledge without knowing the true vector
field $f$.

\newcommand{\pendulumsubfigure}[2]{
  \begin{subfigure}[t]{.3\textwidth}
    \includegraphics[width=\textwidth, trim=0 10 0 0,
    clip]{imgs/pendulum/#2}
    \caption{#1}
  \end{subfigure}
}

\begin{figure}
  \centering
%   \pendulumsubfigure{\label{fig:pendulum_true_vf_2}The true vector field (unkown to us)}{pendulum_true_vf_with_data.png}
  \pendulumsubfigure{\label{fig:pendulum_no_side_info}No side
    information}{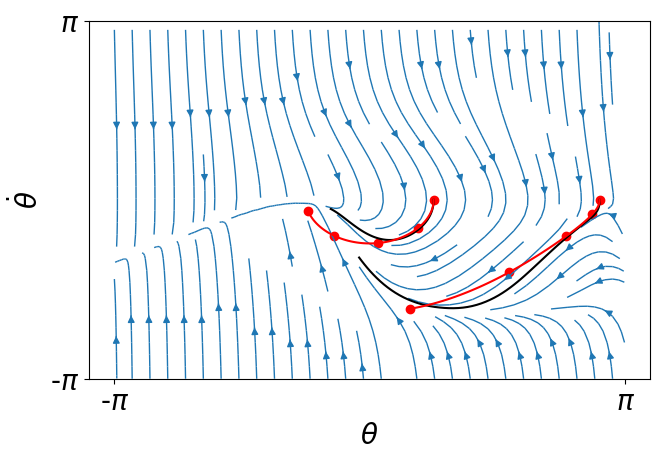}
%   \pendulumsubfigure{no side
%     information}{deg_5_sym_false_sign_false_hamiltonian_false.png}
  \pendulumsubfigure{\label{fig:pendulum_sym}
    $\NSYM$}{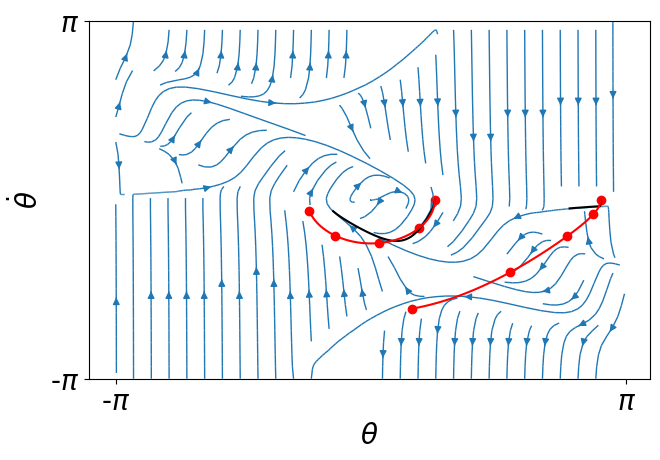}
  \pendulumsubfigure{\label{fig:pendulum_sym_pos}
    $\NSYM \cap
    \NPOS$}{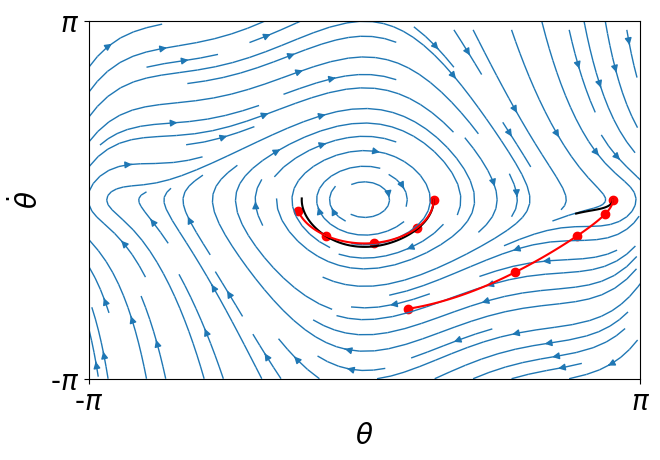}
  \pendulumsubfigure{\label{fig:simple-pendulum-learned-vf-all-sideinfo}
    $\NSYM \cap \NPOS \cap
    \HAMILTONIAN$}{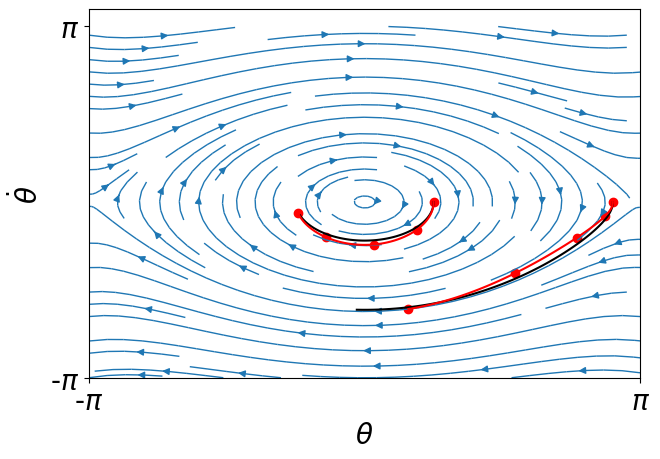}
  
  \caption{\label{fig:simple-pendulum-learned-vf} %Streamplot of the vector field in \eqref{eq:ode-pendulum}
    %(\cref{fig:pendulum_true_vf_2}) along with
    \bachir{Streamplots} of
    polynomial vector fields of degree $5$ that best agree
    with the data (in the least-squares sense) and obey an
    increasing number of side information constraints.
    %(\cref{fig:pendulum_no_side_info,fig:pendulum_sym,fig:pendulum_sym_pos,fig:simple-pendulum-learned-vf-all-sideinfo}).
    In
    each case, the trajectories of the learned vector field
    starting from the same two initial conditions as the
    trajectories observed in the training set are plotted in
    black.}

\end{figure}

\newcommand{\trajectoryplot}[2]{
\begin{figure}[H]\centering
\begin{tikzpicture}
\begin{axis}
  \addplot table [x=t, y=theta, col sep=comma, dashed] {data/pendulum/#1};
\end{axis}
\end{tikzpicture}
\caption{#2}
\end{figure}
}

\begin{figure}
  \scriptsize
  \centering
\begin{tikzpicture}
  \begin{axis}[align=center,  enlargelimits=false, xlabel=$t$, xlabel near ticks, ylabel=$\theta(t)$, xlabel shift=-10pt, xtick={0,10}, yticklabels={$-\frac{\pi}4$, $0$, $\frac{\pi}{4}$}, ytick={-0.7854, 0, 0.7854}, ymin=-1, ymax=1.4, width=.4\textwidth, ylabel shift=-10pt, height=3cm]
    
    \addplot[restrict x to domain=0:3.5, mark=*, red, mark repeat={8}] table [x=t, y=theta1, col sep=comma]
    {data/pendulum/trajectory_simple_pendulum_true_vf.csv};
    \addplot[smooth, restrict x to domain=0:10,  red, dashed] table [x=t, y=theta1, col sep=comma]
    {data/pendulum/trajectory_simple_pendulum_true_vf.csv};
    \addplot[smooth] table [x=t, y=theta, restrict x to domain=0:11,
    col sep=comma]
    {data/pendulum/trajectory_simple_pendulum_vf_no_side_info.csv};
  \end{axis}
\end{tikzpicture}
\begin{tikzpicture}
  \begin{axis}[legend pos=outer north east, enlargelimits=false, xlabel=$t$, xlabel near ticks, ylabel={}, xlabel shift=-10pt,xtick={0,10}, yticklabels={}, ytick={-0.7854, 0, 0.7854}, ymin=-1, ymax=1.4, width=.4\textwidth, height=3cm]
    \addlegendentry{Training Data}
    \addplot[restrict x to domain=0:3.5, mark=*, red, mark repeat={8}] table [x=t, y=theta1, col sep=comma]
    {data/pendulum/trajectory_simple_pendulum_true_vf.csv};

    \addlegendentry{Ground truth}
    \addplot[smooth, restrict x to domain=0:10,  red, dashed] table [x=t, y=theta1, col sep=comma]
    {data/pendulum/trajectory_simple_pendulum_true_vf.csv};

    \addlegendentry{Learned vector field}
    \addplot[smooth] table [x=t, y=theta, restrict x to domain=0:11,
    col sep=comma]
    {data/pendulum/trajectory_simple_pendulum_vf_all_side_info.csv};

  \end{axis}
\end{tikzpicture}
\caption{\label{fig:trajectory-pendulum}Comparison of the trajectory
  of the simple pendulum in \eqref{eq:ode-pendulum} starting from
  $(\frac \pi4, 0)^T$ (dotted) with the trajectory from the same
  initial condition of the polynomial vector field of degree 5 that
  best agrees with the data (in the least-squares solution) in the
  absence of side information (left), and subject to side information
  constraints $\NSYM \cap \NPOS \cap \NHAMILTONIAN$ (right).}
\vspace{-5mm}
\end{figure}
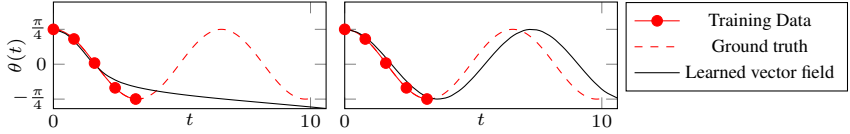

%%% Local Variables:
%%% mode: latex
%%% TeX-master: "../../Learning_with_side_information_SIREV_v1"
%%% End:

\vspace{1em}
\begin{itemize}  \setlength\itemsep{1em}
\item {\bf Sign symmetry ($\NSYM$).} The pendulum obviously behaves symmetrically with respect to the vertical axis (plotted with a dotted line in 
  \cref{fig:fig-simple-pendulum}). We therefore require our candidate
  vector field $p$ to satisfy the same symmetry condition
  \[p(- \theta, - \dot \theta) = - p(\theta, \dot \theta)
    \quad \forall (\theta, \dot \theta) \in \Omega.\]
  Note that this is an affine constraint in the coefficients
  of the polynomial $p$, and that the true vector field $f$
  in \eqref{eq:ode-pendulum} satisfies this constraint.
  
\item {\bf Coordinate nonnegativity ($\NPOS$).} 
% The only
%   external force applied on the pendulum system is that of
%   gravity; see \cref{fig:fig-simple-pendulum}. 
%  
  We know that the force of gravity
  pulls the pendulum's mass down and pushes the angle $\theta$
  towards~$0$. This means that the angular velocity
  $\dot \theta$ decreases when $\theta$ is positive and
  increases when $\theta$ is negative.  Mathematically, we
  must have
  \[p_2(\theta, \dot \theta) \le 0 \; \forall (\theta, \dot
    \theta) \in [0, \pi] \times [-\pi, \pi] \text{ and } p_2(
    \theta, \dot \theta) \ge 0 \; \forall (\theta, \dot
    \theta) \in [-\pi, 0]\times [-\pi, \pi].\] We replace
  each one of these constraints with their corresponding
  $\text{degree-}4$ sos certificate
  (see the definition following equation~\eqref{eq:putinar-certificate}).  (Note that, because
  of the previous symmetry side information, we actually only
  need to impose one of these two constraints.)
\item{\bf Hamiltonian ($\NHAMILTONIAN$).} 
%The system in \eqref{eq:ode-pendulum} is {Hamiltonian}. Indeed, 
    In the
  simple pendulum model, there is no dissipation of energy
  (through friction for example), so the total energy 
  \begin{equation}
    E(\theta, \dot \theta) = \frac {1}2 ml^2{\dot
      \theta}^2 + mgl  (1 -
    \cos(\theta))\label{eq:total-energy-pendulum}
  \end{equation}
  is conserved. The two terms appearing in this equation can
  be interpreted physically as the kinetic and the potential
  energy of the system. Furthermore, the total energy $E$ satisfies
    \[\dot \theta(t) = \frac1{m\ell^2} \frac{\partial E}{\partial
\dot \theta}(\theta(t), \dot \theta(t)),\text{ and } \ddot
\theta(t) = -\frac1{m\ell^2} \frac{\partial E}{\partial \theta}(\theta(t), \dot
\theta(t)).\]
  The simple pendulum system is therefore a Hamiltonian system, with the associated Hamiltonian function $\frac{E}{m \ell^2}$.
  %In physics, dynamical systems that satisfie
    %
  Note that neither the vector field in
  \eqref{eq:ode-pendulum} describing the dynamics of the
  simple pendulum nor the associated Hamiltonian are polynomial functions.
  In our learning procedure, we use only the fact that the
  system is Hamiltonian, i.e., that there exists a function
  $H$ such that
\begin{equation}\label{eq:pendulum_ham}
  p_1(\theta, \dot \theta) = \frac{\partial H}{\partial
\dot \theta}(\theta, \dot \theta),\text{ and } p_2(\theta, \dot
\theta) = -\frac{\partial H}{\partial \theta}(\theta, \dot
\theta),
\end{equation}
but not the exact form of this Hamiltonian. Since we are
parameterizing the candidate vector field $p$ as a
$\text{degree-}5$ polynomial, the function $H$ must be a
(scalar-valued) polynomial of degree $6$. The Hamiltonian
structure can thus be imposed by adding affine constraints on
the coefficients of $p$, or by directly learning $H$ and
obtaining $p$ from \eqref{eq:pendulum_ham}.
\end{itemize}
\vspace{1em}

Observe from \cref{fig:simple-pendulum-learned-vf} that as
more side information is added, the behavior of the learned
vector field gets closer to the truth. In particular, the
solution returned by our semidefinite program in
\cref{fig:simple-pendulum-learned-vf-all-sideinfo} is almost
identical to the true dynamics in
\cref{fig:fig-simple-pendulum} even though it is obtained
only from $10$ noisy samples on two
trajectories. \cref{fig:trajectory-pendulum} shows that even
if we start from an initial condition from which we have made
partial trajectory observations, using side information can
lead to better predictions for the future of the trajectory.

\subsection{Growth of cancerous tumor cells}

\begin{figure}[H]
\begin{minipage}{.45\textwidth}
  \includegraphics[scale=.38]{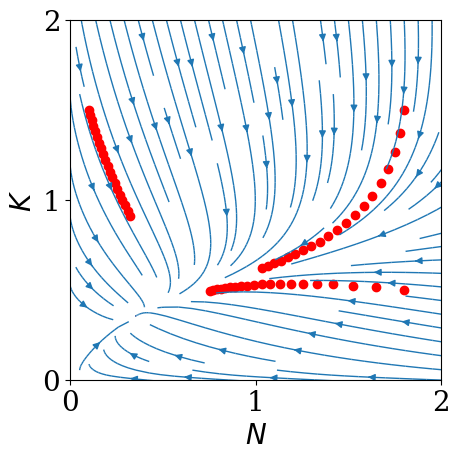}
\end{minipage}
\begin{minipage}{.5\textwidth}
  \caption{\label{fig:two_compartment_model_vf}Streamplot of
    the vector field in \cref{eq:f_g_cancer_vf} describing
    the time evolution of the volume $N$ of a cancerous tumor
    and the host's carrying capacity $K$  (in \textit{cubic
      centimeters}). We consider this vector field to be the
    ground truth and unkown to us. We try to learn it over
    $[0, 2]^2$ from noisy measurements of three trajectories
    (plotted in red).}
\end{minipage}
\vspace{-5mm}
\end{figure}

In this subsection, we consider a model governing the time
evolution of the volume $N$ of a cancerous tumor inside a
human body \cite{tumor_growth_ode_2010}. Cancerous tumors
depend for their development on availability of the so-called
{\it Endothelial} cells,
%
%a
%special type of cells that are responsible for %the formation
%of new blood vessels.
the supply of which is characterized by a quantity called the {\it carrying capacity} $K$.
Intuitively, $K$, which has the same unit as $N$, is proportional to the physical and energetic resources available for cell growth.

Two common modeling assumptions in cancer cell biology are that (i) the growth rate of the tumor decreases as the tumor grows, and (ii) that the volume of the tumor increases (resp. decreases) if it is below (resp. above) the  carrying capacity. We follow the dynamics proposed in \cite{tumor_growth_ode_2010}, which in contrast to prior works in that literature, also models the time evolution of the carrying capacity.
% }, we use the notation $N(t)$ for the volume of the tumor at time $t$ and $K(t)$ for the carrying capacity 
% with a time-varying quantity, $K(t)$,, and has the same unit as $N(t)$. 
The dynamics reads
\begin{equation}\label{eq:dynamics-tumor}
\begin{pmatrix}\dot N(t)\\\dot K(t)\end{pmatrix} = f(N(t), K(t)),
\end{equation}
where $f:\mathbb R^2 \rightarrow \mathbb R^2$ is given by
\begin{equation}
  f_1(N, K) \coloneqq \frac{\mu N}{\nu} \left[1 - \left(\frac
      NK\right)^\nu \right], \quad f_2(N, K) \coloneqq \omega N -
  \gamma N^{\frac23} K.\label{eq:f_g_cancer_vf}
\end{equation}
Here, $\mu, \nu, \gamma$ and $\omega$ are positive
parameters. The dynamics of $N$ is motivated by the so-called
generalized logistic differential equation; the term
$\omega N$ models the influence of the tumor on the
Endothelial cells via short-range stimulators, and the term
$- \gamma N^{\frac 23} K$ models this same influence via
long-range inhibitors. See \cite{tumor_growth_ode_2010} for
more details.

Having an accurate model of the growth dyanmics of cancerous
tumors is crucial for the follow-up task of designing
treatment plans, e.g., via radiation therapy. We hope that in
the future leveraging side information can lead to learning
more accurate models directly from patient data (as opposed
to postulating an exact functional form such as
\eqref{eq:f_g_cancer_vf}). For the moment however, we take
\eqref{eq:f_g_cancer_vf} with the following parameters to be
the ground truth:
\[(\mu, \nu, \gamma, \omega) = \frac1{10}(1, 5,  1, 2).\]
See
\cref{fig:two_compartment_model_vf} for a streamplot of the corresponding vector
field.

%
% We take the vector field $f$ in \cref{eq:f_g_cancer_vf} to be
% the (unknown) ground truth, and we observe a noisy version of
% five trajectories plotted in red in
% \cref{fig:two_compartment_model_vf}. 
%
We consider the task of
learning the vector field $f$ over the compact set
$\Omega \coloneqq [0, 2]^2$ from noisy snapshots of three
trajectories. Each trajectory was started from a random
initial conditions
$x_{\text{init}}^{k} \coloneqq (N_{\text{init}}^{k}, K_{\text{init}}^{k})$
(with $k=1,2,3$) inside $\Omega$ and sampled at times
$t_i = i / {20}$, with $i=0,\ldots,19$ (see
\cref{fig:two_compartment_model_vf}). More precisely, we have
access to the following training data:
\begin{equation}
  \scriptsize
    \mathcal D \coloneqq 
  \left\{ \left( (N(t_i, x_{\text{init}}^{k}), K(t_i, x_{\text{init}}^{k})), (\dot N(t_i, x_{\text{init}}^{k}) + 10^{-4}\varepsilon_{i,k,1}, \dot K(t_i, x_{\text{init}}^{k})+ 10^{-4}\varepsilon_{i,k,2})\right) \right\}_{0 \le i < 20, 1 \le k\le3},
  \label{eq:data-cancer}
\end{equation}
where the $\varepsilon_{i,k,l}$ (for $i=0,\ldots,19$, $k=1,\ldots,3$, and $l=1,2$) are independent standard normal variables.
We parameterize our candidate vector field $p$ as a
$\text{degree-}5$ polynomial. Without any side infromation, fitting this candidate vector
field to the data in \cref{eq:data-cancer} via a
least-squares problem leads to the vector field plotted
in \cref{fig:cancer-no-side-info}. Once again, the vector field
obtained in this way is very far from the true vector field.

\newcommand{\cancersubfigure}[3]{
  \begin{subfigure}[t]{#3}
    \includegraphics[width=3cm, trim=0 10 0 0,
    clip]{imgs/cancer/#2}
    \caption{#1}
  \end{subfigure}
}
\begin{figure}
  \centering
  
%   \cancersubfigure{\label{fig:cancer-true-vf-2}The true
%     vector field (unknown to us)}{cancer_true_vf.png}{.3\textwidth}
%   \quad 
  \cancersubfigure{\label{fig:cancer-no-side-info}No side information}{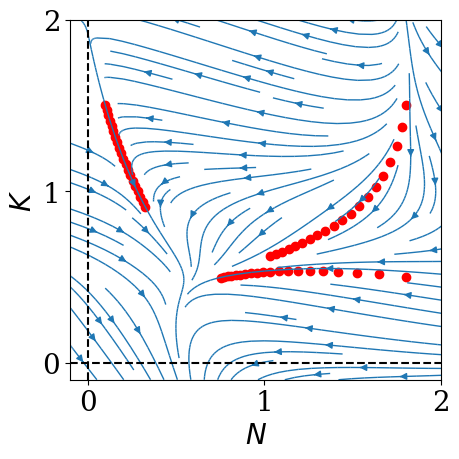}{.3\textwidth}
  \quad
  \cancersubfigure{\label{fig:cancer-posmon}$\NPOS\text{-}\NMON$}{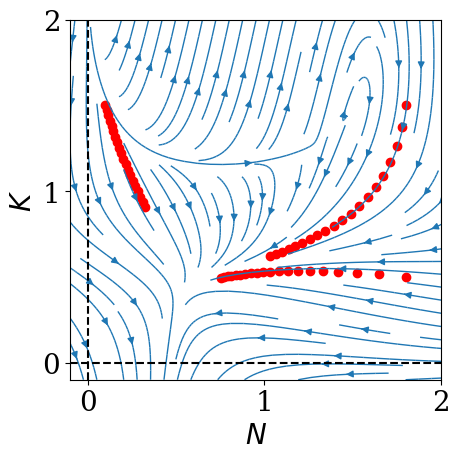}{.3\textwidth}

  \cancersubfigure{\label{fig:cancer-posmon-inv}$\NPOS\text{-}\NMON \cap \NINV$}{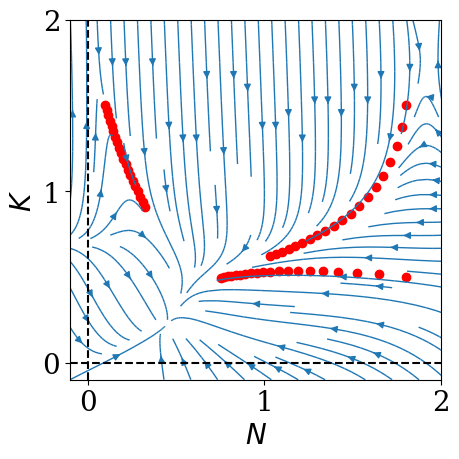}{.3\textwidth}
  \quad
  \cancersubfigure{\label{fig:cancer-posmon-inv-pos}$\NPOS\text{-}\NMON \cap \NINV \cap  \NPOS$}{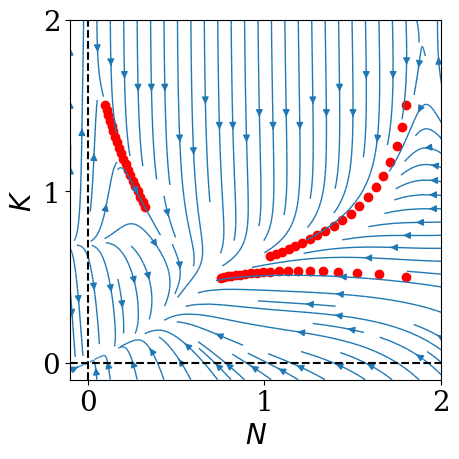}{.3\textwidth}
  \quad
  \cancersubfigure{\label{fig:cancer-learned-vf-all-sideinfo}$\NPOS\text{-}\NMON
    \cap \NINV \cap \NPOS$\\$\cap \NINTERPOLATION $}{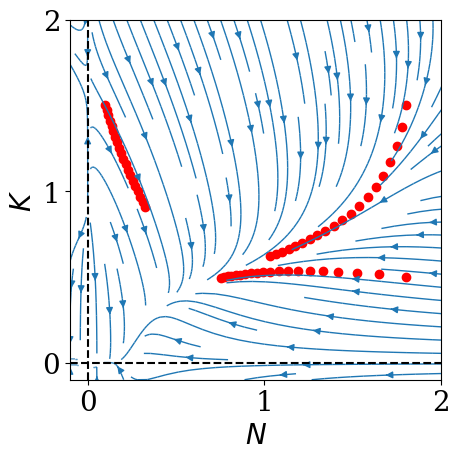}{.3\textwidth}
  % \cancersubfigure{$\NINTERPOLATION$}{cancer_learned_vf_eq.png}{.32\textwidth}
  % %
  % \cancersubfigure{$\NINTERPOLATION \cap \NINV$}{cancer_learned_vf_eq_nonneg.png}{.32\textwidth}
  % %
  % \cancersubfigure{$\NINTERPOLATION \cap \NINV \cap \NPOS$}{cancer_learned_vf_eq_nonneg_Ndot_K_ge_N.png}{.32\textwidth}
  % %
  % \quad
  % \cancersubfigure{\label{fig:cancer-learned-vf-all-sideinfo}$\NINTERPOLATION \cap \NINV \cap \NPOS \cap \NPOS\text{-}\NMON$}{cancer_learned_vf_eq_nonneg_Ndot_K_ge_N_Ndot_rate.png}{.4\textwidth}
 
 \caption{\label{fig:cancer-learned-vf}
   %Streamplot of the vector field in \eqref{eq:f_g_cancer_vf} (\cref{fig:cancer-true-vf-2}) along with  
   \bachir{Streamplots} of polynomial vector fields of degree $5$
   that best agree with the data (in the least-squares sense) and obey an increasing number of side
   information constraints. %(\cref{fig:cancer-no-side-info,fig:cancer-posmon,fig:cancer-posmon-inv,fig:cancer-posmon-inv-pos,fig:cancer-learned-vf-all-sideinfo}).
   }
\vspace{-5mm}
\end{figure}

To do a better job at learning, we impose the side
information constraints listed below that come from expert knowledge in the tumor growth literature (see, e.g.,
\cite{tumor_growth_ode_2010,reference_one_tumor,reference_two_tumor}):

\vspace{1em}
\begin{itemize}  \setlength\itemsep{1em}
  \item {\bf A mix between coordinate nonnegativity and
    coordinate directional monotonicity
    ($\NPOS\text{-}\NMON$).}
  %
  % A known fact about tumor cells is that their growth rate,
  % $\frac{\dot N}{N}$, decelerate as they grow larger, i.e., %
  % \[ \frac{\partial}{\partial N} \left(\frac{\dot N}{N}\right) \le
  %   0.\] %
  As stated in \cite{tumor_growth_ode_2010},
  ``one of the few near-universal observations about solid
  tumors is that almost all decelerate, i.e., reduce their
  specific growth rate $\frac{\dot N}{N}$, as they grow
  larger.''

  Based on this contextual knowledge, our candidate vector field $p$ should satisfy
  \[ \frac1{N} \frac{\partial p_1}{\partial N}(N, K) -
      \frac1{N^2} p_1(N, K) \le 0 \quad \forall N \in (0, 2], \; \forall  K \in [0, 2].\] 
  Since the state variable $N$ is nonnegative at all times,
  we can clear denominators to obtain the constraint
     \[N \frac{{\partial} p_1}{{\partial}N}(N, K) -
      p_1(N, K) \le 0 \; \forall (N, K) \in [0, 2]^2.\]    
This is a polynomial nonnegativity constraint over a closed basic semialgebraic set.
\item {\bf Invariance of the nonnegative orthant ($\NINV$).} The
  state variables $N$ and $K$ quantify volumes, and
  as such, should be nonnegative at all times. This
  corresponds to the nonnegativity constraints
  \[p_2(N, 0) \ge 0\;\forall N \in [0, 2], \quad p_1(0, K)\ge0\; \forall K \in [0, 2].\]
  \item{\bf Coordinate nonnegativity ($\NPOS$).}   As mentioned before, the rate of change $\dot N$ in
    the tumor volume is nonnegative when the
    carrying capacity $K$ exceeds $N$, and nonpositive
    otherwise~\cite{tumor_growth_ode_2010}. Mathematically, we must have
    \[ p_1(N, K) \ge 0 \quad \forall N \in [0, 2],\; \forall K \in [N, 2],\]
    \[p_1(N, K) \le 0 \quad \forall N \in [0, 2],\; \forall K \in [0, N].\]
\item {\bf Equilibrium point at the origin
    ($\NINTERPOLATION$).} The tumor does not grow if the
  volume of cancerous cells and the carrying capacity are
  both zero. This corresponds to the constraint
  $p(0, 0) = 0$.  
\end{itemize}
\vspace{1em}

We observe from \cref{fig:cancer-learned-vf} that as
more side information is added, the behavior of the learned
vector field gets closer and closer to the ground truth. In particular, the
solution returned by our semidefinite program in
\cref{fig:cancer-learned-vf-all-sideinfo} is very close to the true dynamics in \cref{fig:two_compartment_model_vf}.

\subsection{\bachir{Learning the Lorenz system}}
\label{sec:lorenz}

\begin{figure}
  \centering
  \begin{minipage}{.3\textwidth}
    \includegraphics[scale=.3]{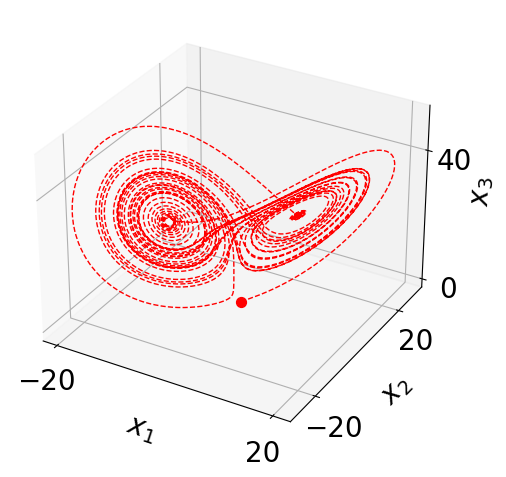}
  \end{minipage}
  \begin{minipage}{.39\textwidth}
    \caption{\label{fig:traj_lorentz}\bachir{Sample trajectory of the system in \eqref{eq:lorentz} starting from $x(0) = (1, 1, 1)^T$.}}
\end{minipage}
\end{figure}

\bachir{
In this subsection, we consider the classical Lorenz system (see, e.g., \cite{lorenzbook}) that is known for the chaotic properties of its solutions: %(i.e., slight changes in initial conditions lead to large changes in trajectories).
%It models the movement of a two dimensional fluid (e.g., air) generated by the difference in temperature between its two sides. More formally, the Lorenz system is given by
\begin{equation}  \label{eq:lorentz}
    \dot \ix(t) = f(\ix(t)), \text{where} \; \ix(t) \in \R^3 \text{ and }
  f(\ix) = \left(
    \begin{array}{l}
    \sigma (x_2 - x_1)\\
      x_1 (\rho - x_3) - x_2\\
      x_1 x_2 - \beta x_3
    \end{array}\right).
\end{equation}
% Here, $x_1$ represents the convective intensity, and $x_2$ (resp. $x_3$) represents the horizontal (resp. vertical) temperature variation. The parameters $\sigma$, $\rho$, and $\gamma$ are positive physical constants that relate to the {\it  Prandtl} number,  the {\it Rayleigh number}, and the physical dimensions of the fluid respectively.  
Here, we work with the commonly used parameter values
\[\rho = 2, \sigma = 10, \text{ and } \beta = \frac 8 3.\]
We consider the task of learning the Lorenz system in \eqref{eq:lorentz} from 20 noisy snapshots of $N_t$ trajectories. More precisely, our data is $\mathcal D^1 \cup \ldots \cup \mathcal D^{N_t}$, where for $k=1, \ldots, N_t$, 
\begin{equation}\label{eq:training_data_lorenz}
\mathcal D^k \coloneqq \left\{ \left(\ix(t_i, x_{\text{init}}^k),
      f(\ix(t_i, x_{\text{init}}^k))+
      10^{-1} \begin{pmatrix}\varepsilon^k_{i,1}\\\varepsilon^k_{i,2}\\\varepsilon^k_{i,3}\end{pmatrix}\right)
  \right\}_{i=1}^{20}.
  \end{equation}
 Here, $x_{\text{init}}^k$ is generated uniformly at random from the box $[0, 10]^3$, and the scalars $\varepsilon^k_{i,l}$ (for $i=1,\ldots,20$, and $l=1,\ldots,3$) are independent standard normal variables.
 We parameterize our candidate  vector field $p: \mathbb R^3 \rightarrow \mathbb R^3$ as a polynomial of degree $3$.
 The degree is taken to be larger than $2$ because we do not want to assume knowledge of the degree of the true vector field in \cref{eq:lorentz}. We assume access 
to the following list of side information:
\begin{itemize}
    \item 
 {\bf Equilibrium point at the origin
    ($\NINTERPOLATION$).}  This corresponds to the constraint $p(0)~=~0$.  
\item {\bf Coordinate nonnegativity ($\NPOS$).}
The state variable $x_1$ should increase when $x_2~\ge~x_1$ and decrease otherwise. In other words, 
\[p_1(\ix) \ge 0 \; \forall \ix \in S, \quad p_1(\ix) \le 0 \; \forall \ix \in S', \]
where $S \coloneqq \{ x \in \mathbb R^3 \; | \; x_2 - x_1 \ge 0 \}$ and $S' \coloneqq \{ x \in \mathbb R^3 \; | \; x_2 - x_1 \le 0 \}.$
\item {\bf Monotonicity ($\NMON$).} The rate of change of any state variable decreases as the variable gets larger. In other words, 
\[\frac{\partial p_i}{\partial x_i}(x) \le 0 \quad \forall x \in \mathbb R^3, \quad i = 1, 2, 3.\]
\end{itemize}
\vspace{1em}
Since a 3-dimensional streamplot is too clutered to be insightful, we instead report in \cref{tbl:lorentz} the test error of the polynomial vector field of degree $d=3$ learned from the side information and the data described above.  Here, the test error of the vector field $p$ is measured as
\[ \frac1{1000} \sum_{x \in D} \|p(x) - f(x)\|,\]
where $D$ is a uniform discretization of the box $[-10, 10]^3$ with 1000 samples. As one can observe, the test error is reduced most of the time as more side information constraints are incorporated. For experiments with other degrees, noise levels, and a different notion of test error, see the appendix.
}

\begin{table}
\begin{tabular}{lcccc}
\toprule
     Side information: &      $\emptyset$ &       $\NINTERPOLATION$  &         $\NINTERPOLATION \cap \NPOS$ &        $\NINTERPOLATION \cap \NPOS \cap \NMON$ 
   \\
\midrule
$N_t=1$  &   1.47e+03 &   516 &   664 &  54.9 \\
       $N_t=2$ &        4.61 &  2.51 &  2.25 & 0.963  \\
       $N_t=3$ &     1.38 & 0.943 & 0.787 & 0.352 \\
\bottomrule
\end{tabular}
\caption{\label{tbl:lorentz} \bachir{Test error of cubic vector fields learned from side information and noisy snapshots of $N_t$ trajectories of the Lorenz system.}}
\end{table}

\subsection{Following learning with optimal control}
\label{sec:opt_control}
In this subsection, we revisit the contagion dynamics
\eqref{eq:diffusion_disease}
%that was studied in
%\Cref{sec:diffusion_contagious_disease}, 
and study the effect of side information on policy decisions to contain an outbreak. Suppose that by an initial screening of a random subset of the population, it is estimated that a fraction $0.5$ (resp. $0.4$) of males (resp. females) are infected with the disease. We would like to contain the outbreak by performing daily widespread testing of the population. We introduce two
control decision variables $u_1$ and $u_2$, representing respectively the fraction of the population of males and females that are tested per unit of time. We suppose that testing slows down the spread of the disease (due e.g. to appropriate action that can be taken on the positive cases) as follows:
\begin{equation}  \label{eq:diffusion_disease_control}
  \dot \ix(t) = f(\ix(t))
  - \begin{pmatrix}u_1x_1\\u_2x_2\end{pmatrix} ,
\end{equation}
 where $f(x(t))$ is the unknown dynamics of the spread of the disease in the absence of any control.
 
%  vector field in \eqref{eq:diffusion_disease} with parameters in \eqref{eq:disease_params}. $u_1, u_2 \in [0, 1]$, $x(0) = (0.5, 0.4)$, and
 
We suppose that the monetary cost of testing a fraction $u_1$ of males and
$u_2$ of females is given by $\alpha (u_1+u_2)$ for some known positive scalar $\alpha$. Our goal is to minimize the sum
\begin{equation}
  \label{eq:cost_control}
  c(u_1, u_2) \coloneqq x_1(T, \hat x_{\text{init}} ) + x_2(T, \hat x_{\text{init}} ) + \alpha (u_1+u_2),
\end{equation}
of the total number $x_1(T, \hat x_{\text{init}} ) + x_2(T, \hat x_{\text{init}})$ of infected individuals at
the end of a desired time period~$T$, and the monetary cost
$\alpha (u_1+u_2)$ of our control law. Here, $x_1(t,\hat x_{\text{init}} )$ and $x_2(t, \hat x_{\text{init}} )$ evolve according to \eqref{eq:diffusion_disease_control} when started from the point $\hat x_{\text{init}} = (0.5, 0.4)^T$.
In our experiments, we take $T = 20$, $\alpha = 0.4$, and $f$ to be the vector field in \eqref{eq:diffusion_disease} with parameters in \eqref{eq:disease_params}.

\begin{figure}
  \centering
  \begin{minipage}{.6\textwidth}
    \includegraphics[scale=.085]{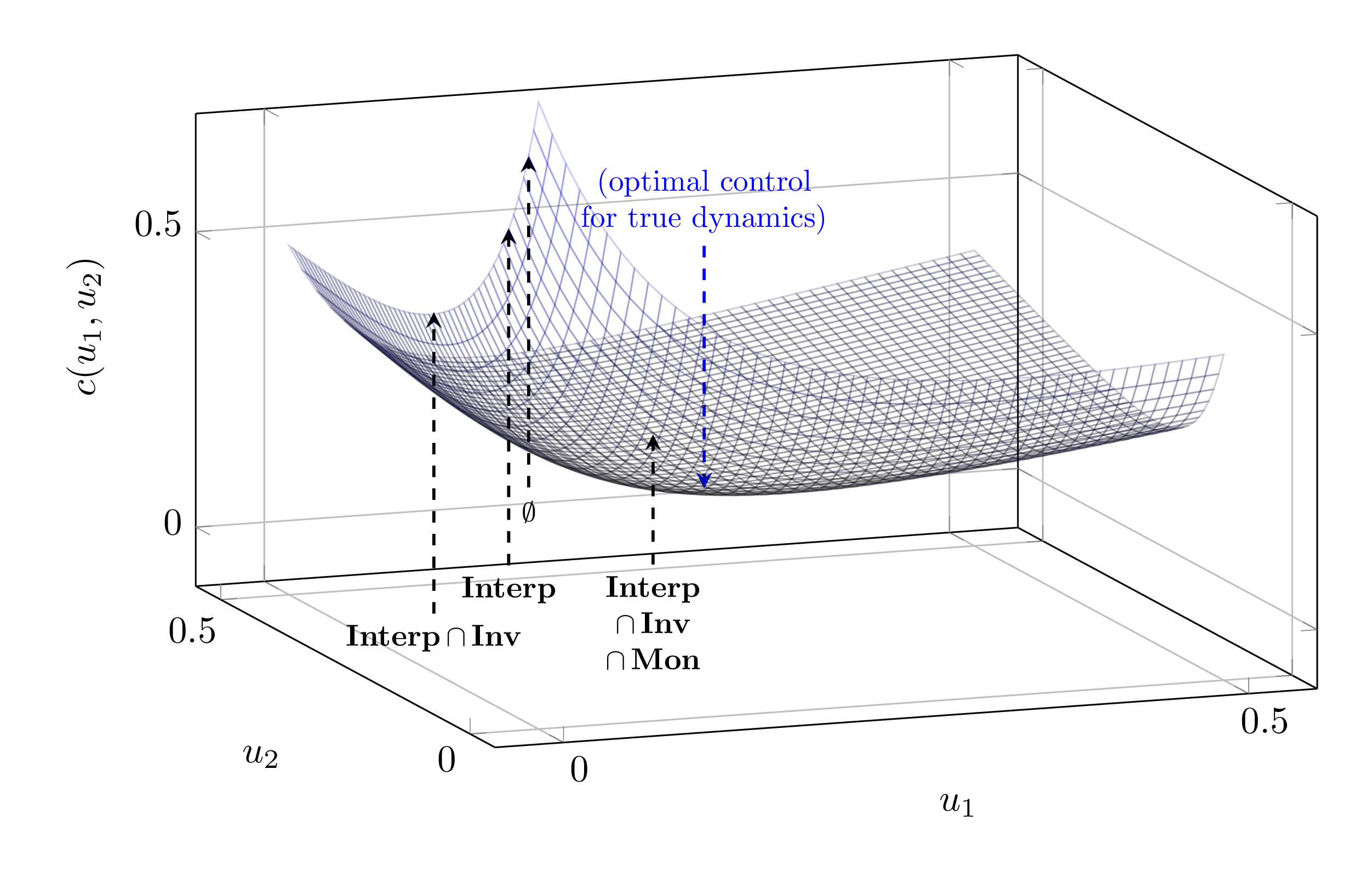}
  \end{minipage}
  \begin{minipage}{.39\textwidth}
    \caption{\label{fig:opt_control}The graph of the function
      $c(u_1, u_2)$ in \eqref{eq:cost_control} with $T~=~20$,
      $\alpha = 0.4$, and $f$ as in
      \eqref{eq:diffusion_disease} with parameters in
      \eqref{eq:disease_params}. The minimizer of the
      function $c(u_1, u_2)$, which corresponds to the
      optimal control law, is indicated with a blue arrow.
      The control laws that are optimal for dynamics learned
      from a single trajectory of $f$ with different side
      information constraints are indicated with black
      arrows.}
\end{minipage}
\end{figure}

Given access to the vector field $f$, one could simply design
an optimal control law $(u_1^*, u_2^*)$ that minimizes the
cost function $c(u_1, u_2)$ in \eqref{eq:cost_control} by
gridding the control space $[0, 1]^2$, and computing
$c(u_1, u_2)$ for every point of the grid. Indeed, for a
given $(u_1, u_2)$, $c(u_1, u_2)$ can be computed by
simulating the dynamics in
\eqref{eq:diffusion_disease_control} from the initial
condition $\hat x_{\text{init}}$. The optimal control law obtained
by following this strategy is depicted in
\Cref{fig:opt_control} together with the graph of the
function $c(u_1, u_2)$.

We now consider the same setup as \Cref{sec:diffusion_contagious_disease}, where we do not know the vector field $f$, but have observed $20$ noisy samples of a single trajectory of it starting from $(0.7, 0.3)^T$ (c.f. \eqref{eq:training_data_epi}). We design our control laws instead for the four polynomial vector fields of degree $3$ that were learned in \Cref{sec:diffusion_contagious_disease} with no side information, with $\NINTERPOLATION$, with $\NINTERPOLATION \cap \NINV$, and with $\NINTERPOLATION \cap \NINV \cap \NMON$. This is done by following the procedure described in the previous paragraph, but using the learned vector field instead of $f$. The corresponding four control laws are depicted in \Cref{fig:opt_control} with black arrows. We emphasize that while the control laws are computed from the learned vector fields, their associated costs in \Cref{fig:opt_control} are computed by applying them to the true vector field. It is interesting to observe that adding side information constraints during the learning phase leads to the design of better control laws.

From \Cref{tbl:opt_control}, we see that if we had access to the true vector field, an optimal control law would lead to the eradication of the disease by time $T$. The first four rows of this table demonstrate the fraction of infected males and females at time $T$ when control laws that are optimal for dynamics learned with different side information constraints are applied to the true vector field. It is interesting to note that with no side information, a large fraction of the population remains infected, whereas control laws computed with three side information constraints are able to eradicate the disease almost completely.

% In this section, We assume that the true vector
% field $f$ is not know, and we use the polynomial vector
% fields of degree $3$ that were learned from training data and
% side information constraints in
% \Cref{sec:diffusion_contagious_disease} to design the
% control $(u_1, u_2)$. More specifically, the control law
% corresponding to a learned vector field $p$ is given by
% \[\arg\min_{u_1, u_2} x_1'(T) + x_2'(T) + \alpha (u_1+u_2),\]
% where $x_1'(T)$, $x_2'(T)$ are the values of the state
% variables of the dynamical system
% \begin{equation}  \label{eq:diffusion_disease_control_p}
%   \dot \ix'(t) = p(\ix'(t)) - \begin{pmatrix}u_1x_1'\\u_2x_2'\end{pmatrix}
% \end{equation}
% at time $T$ when started from $x_{\text{init}}' = (0.7, 0.3)$.
% We report the results in \cref{fig:opt_control}. As
% expected, as we consider more side information constraints,
% the performance of the learned control law improves.

\begin{table}
\centering
\begin{tabular}{c|c|c}
Side information & $x_1(T, \hat x_{\text{init}})$ & $x_2(T, \hat x_{\text{init}})$ \\
\hline
None & 0.45 &  0.40\\
$\NINTERPOLATION$ & 0.41 &  0.29\\
$\NINTERPOLATION \cap \NINV$ & 0.31 &  0.12\\
$\NINTERPOLATION \cap \NINV \cap \NMON$ & 0.01 &  0.01\\
\hline
True vector field & 0.00 &  0.00\\
\end{tabular}
\caption{\label{tbl:opt_control} The first four rows indicate
  the fraction of infected males and females at the end of
  the period $T$ when a control law, optimal for dynamics
  learned from a single trajectory with different side
  information constraints, is applied to the true vector
  field. The last row indicates the fraction of infected
  males and females at time $T$ resulting from applying the
  optimal control law computed with access to the true
  dynamics.}
\vspace{-5mm}
\end{table}

% % \begin{figure}
% %     \centering
%     \includegraphics[scale=.5]{imgs/tables/disease_table_deg_2.png}
%   \includegraphics[scale=.5]{imgs/tables/disease_table_deg_3.png}
%   \includegraphics[scale=.5]{imgs/tables/disease_table_deg_4.png}
%     {\bf Box error and trajectory error of vector fields learned with and without side information. (Disease Example)}
% %     \label{fig:tbl_disease}
% % \end{figure}

% \newpage
% \begin{figure}
%     \centering
%     \includegraphics[scale=.5]{imgs/tables/lorentz_table_deg_2.png}
%     \includegraphics[scale=.5]{imgs/tables/lorentz_table_deg_3.png}
%     \caption{Box error and trajectory error of vector fields learned with and without side information. (Lorenz)}
%     \label{fig:tbl_lorenz}
% \end{figure}

\section{Approximation Results}
\label{sec:approx_results}

In this section, we present some density results for polynomial vector fields that obey side information. This provides some theoretical justification for our choice of parameterizing our candidate vector fields as  polynomial functions.

More precisely, we are interested in the following question:
Given a continuously-differentiable vector field $f$
satisfying a list of side information constraints from
\Cref{sec:side-info}, is there a polynomial vector field that
is ``close'' to $f$ and satisfies the same list of side
information constraints? For the purpose of learning
dynamical systems, arguably the most relevant notion of
``closeness'' between two vector fields is one that measures
how differently their corresponding trajectories can behave
when started from the same initial condition.  More formally,
we fix a compact set $\Omega \subset \R^n$ and a time horizon
$T$, and we define the following notion of distance between
any two vector fields {$f,g \in \CD$}:
\begin{equation}\label{eq:metric_dyn}
  \metricdyn{f, g}{\Omega}{T} \coloneqq \sup_{ (t, \ix_{\text{init}}) \in \mathcal S}
  \max\left\{\ltwonorm{\ix_f(t, \ix_{\text{init}}) - x_g(t,
      \ix_{\text{init}})}, \ltwonorm{\dot \ix_f(t, \ix_{\text{init}}) - \dot x_g(t,
      \ix_{\text{init}})}\right\},
\end{equation} 
where $\ix_f(t, \ix_{\text{init}})$ (resp.
$x_g(t, \ix_{\text{init}})$) is the trajectory starting from
$\ix_{\text{init}} \in \Omega$ and following the dynamics of
$f$ (resp. $g$), and
\begin{equation}\label{eq:S}
    \mathcal S \coloneqq  \{ (t, \ix_{\text{init}})  \in  [0, T] \times  \Omega\; | \;  \ix_f(s, \ix_{\text{init}}), x_g(s, \ix_{\text{init}}) \in \Omega \; \forall s \in [0, t]\}.
\end{equation}
The reason why in the definition of
$d_{\Omega,T}$, we take the supremum  over $\mathcal S$ instead of over
$[0, T] \times \Omega$ is to ensure that the trajectories
that appear in \eqref{eq:metric_dyn} are well defined.

In \Cref{sec:one_side_info}, we show that under some
assumptions that are often met in practice, polynomial vector
fields can be made arbitrarily close to any
continuously-differentiable vector field $f$ (in the sense of
\eqref{eq:metric_dyn}), even if they are required to satisfy
one side information constraint that $f$ is known to
satisfy. In \Cref{sec:multiple_side_info}, we drop our
assumptions and generalize this approximation result to any
list of side information constraints at the price of allowing
an arbitrarily small error in the satisfaction of these
constraints. Furthermore, we show that the approximate
satisfaction of side information can be certified by a sum of
squares proof.

%, then polynomial vector fields can approximate any continuously differentiable vector field $f$ arbitrarily well while satisfying any list of side information constraints that $f$ is known to satisfy.

%in the this statement is no longer true when we have more than one side information constraint

\subsection{Approximating a vector field while (exactly) satisfying
  one side information constraint}
\label{sec:one_side_info}
The following theorem is the main result of this section. We will need the following definition for a subcase of this theorem: Given a collection of sets $A_1, \ldots, A_r$, we define $\mathcal G(A_1, \ldots, A_r)$ to be the graph on $r$ vertices labeled by the sets $A_1, \ldots, A_r$, where two vertices $A_i$ and $A_j$ are connected if $A_i \cap A_j \ne \emptyset$.

\begin{theorem}\label{thm:approximation_one_side_info}
  For any compact set $\Omega \subset \R^n$, time horizon
  $T > 0$, desired accuracy $\varepsilon > 0$, and vector
  field {$f \in \CD$ which} satisfies one of the following
  side information constraints (see \Cref{sec:side-info}):
  
  \begin{enumerate}[(i),leftmargin=.6cm,rightmargin=.9cm]
  \item $\INTERPOLATION{\{(\ix_i,y_i)\}_{i=1}^m}$, where $x_1, \ldots, x_m \in \Omega$,
  \item $\SYM{G,\sigma,\rho}$,
  \item $\POS{ \{(P_i, N_i)\}_{i=1}^n }$, where for each $i \in \{1, \ldots, n\}$,  $P_i \cap N_i = \emptyset$,
  \item $\MON{ \{(P_{ij}, N_{ij}) \}_{i, j=1}^n}$, where for
    each $i, j \in\{1, \ldots, n\}$, the sets $P_{ij}$ and
    $N_{ij}$ belong to different connected components of the
    graph $\mathcal G(P_{i1},N_{i1}, \ldots, P_{in},N_{in})$,
  \item $\INV {\{B_i\}_{i=1}^r}$, where the sets $B_i$ are
    pairwise nonintersecting, and defined as
    $B_i \coloneqq \{ x \in \R^n \; | \; h_{ij}(\ix) \ge
    0,\;j =1,\dots,m_i\}$ for some concave
    continuously-differentiable functions $h_{ij}$
    that satisfy
    \[ \forall i \in \{1, \ldots, r\}, \; \exists x^i \in B_i \text{ such that }   
    %\forall j\in \{1, \ldots, m\}
     h_{ij}(x^i) > 0  \text{ for } j=1,\ldots,m_i,
    \]
  \item $\GRAD$,
  \item[(vi')] $\HAMILTONIAN$,
  \end{enumerate}\vspace{.25cm}
  there exists a polynomial vector field
  $p: \mathbb R^n \rightarrow \R^n$ that satisfies the same side information constraint as $f$ and has
  $\metricdyn{f, p}{\Omega}{T} \le \varepsilon.$
\end{theorem}
% The reader can check that in the numerical experiments of \Cref{sec:applications}, whenever we imposed a side information constraint of type (i)-(vi'), the assumptions of \cref{thm:approximation_one_side_info} were satisfied.
%
Before we present the proof, we recall the classical
Stone-Weierstrass approximation theorem. Note that while the
theorem is stated here for scalar-valued functions, it
readily extends to vector-valued ones.

\begin{theorem}(see, e.g.,
  \cite{stone_weierstrass_1948})\label{thm:weirstrass} For
  any compact set $\Omega \subset \mathbb R^n$, scalar
  $\varepsilon > 0$, and continuous function
  $f:\Omega \rightarrow \R$, there exists a polynomial
  $p:\R^n\rightarrow \R$ such that
\[\max_{x \in \Omega} | f(x) - p(x)| \le \varepsilon.\]
% and
% \[\max_{x \in \Omega, i \in \{1, \ldots, n\}, j \in \{1, \ldots, n\}} \left| \frac{\partial f_i}{\partial x_j}(x) - \frac{\partial p_i}{\partial x_j}(x)\right| \le \varepsilon.\]
\end{theorem}
 For two vector fields $f,g: \R^n \rightarrow \R^n$ and a set $\Omega \subseteq \R^n$, let us define
\[\metricsup{f-g}{\Omega} \coloneqq \max_{\ix \in \Omega}
  \ltwonorm{f(\ix)-g(\ix)}.\]
The following proposition relates this quantity to the notion
of distance $\metricdyn{f, g}{\Omega}{T}$ defined in
\eqref{eq:metric_dyn}.  We recall that for a scalar
$L \ge 0$, a vector field $f$ is said to be
$L\text{-Lipschitz}$ over $\Omega$ if
\[\ltwonorm{f(x) - f(y)} \le L \ltwonorm{x - y} \quad \forall x, y \in \Omega.\]
Note that any vector field $f \in \CD$ is
$L\text{-Lipschitz}$ over a compact set $\Omega$ for {some
  nonnegative scalar $L$.}
\begin{proposition}\label{prop:equivalence-metric-vf}For any
  compact set
  $\Omega \subset \mathbb R^n$, \bachir{any finite} time horizon $T > 0$,  
  and any
  two vector fields
  {$f,g \in \CD$}, we have
\[\metricsup{f-g}{\Omega} \le \metricdyn{f, g}{\Omega}{T} \le
  \max\{T e^{LT}, 1+LT e^{LT}\} \metricsup{f-g}{\Omega},\]
where $L \ge 0$ is any scalar for which either $f$ or $g$ is
$L\text{-Lipschitz}$ over $\Omega$.
\end{proposition}
\bachir{We note that the dependence of the right-most term  on $T$ cannot be avoided. (For example, for $n=1$, $\Omega = [0, 1]$, $f(x) = 0$, $g_\varepsilon(x) = \varepsilon x$, we have $\metricsup{f-g_\varepsilon}{\Omega} = \varepsilon$, but $\metricdyn{f, g_\varepsilon}{\Omega}{T} \ge \frac12$ for all $T \ge \frac{\log 2}{\varepsilon}$.)} To present the proof of this proposition, we need to recall the Gr\"{o}nwall-Bellman inequality.

\begin{lemma}[Gr\"{o}nwall-Bellman inequality
  \cite{bellman1943stability,khalil2002nonlinear}]
\label{lem:gronwall}
Let $I = [a, b]$ denote a nonempty interval on the real
line. Let $u$, $\alpha$, $\beta: I \rightarrow \R$ be
continuous functions satisfying
  \[u(t) \le \alpha(t) + \int_a^t \beta(s) u(s) \; {\rm d}s
\quad \forall t \in I.\] If $\alpha$ is nondecreasing and
$\beta$ is nonnegative, then
  \[u(t) \le \alpha(t) e^{\int_a^t \beta(s) \; {\rm d}s} \quad \forall t \in I.\]
\end{lemma}

\begin{proof}[Proof of \cref{prop:equivalence-metric-vf}]
  We fix a compact set $\Omega \subset \mathbb R^n$, a time
  horizon $T > 0$, and vector fields {$f,g \in \CD$}, with
  $f$ being $L\text{-Lipschitz}$ over $\Omega$ for some
  scalar $L \ge 0$.
  To see that the first inequality holds, note that for any
  $x_{\text{init}} \in \Omega$,
  $\dot \ix_f(0, \ix_{\text{init}}) = f(\ix_{\text{init}})$
  and
  $\dot x_g(0, \ix_{\text{init}}) =
  g(\ix_{\text{init}})$. Therefore,
  $ \metricsup{f - g}{\Omega} \le \metricdyn{f,
    g}{\Omega}{T}$.

  For the second inequality, fix $(t, x_{\text{init}}) \in \mathcal S$, where $\mathcal S$ is defined in \eqref{eq:S}. Let us first bound $\|x_f(t, x_{\text{init}}) -
x_g(t, x_{\text{init}})\|_2$. By definition of $x_f$ and $x_g$, we have
  \begin{align*} x_f(t,  x_{\text{init}}) - x_g(t, x_{\text{init}}) &= \int_0^t f(x_f(s,  x_{\text{init}})) - g(x_g(s,  x_{\text{init}}))
{\rm d}s \\&= \int_0^t f(x_f(s,  x_{\text{init}})) - f(x_g(s,  x_{\text{init}}))\; {\rm d}s 
\\&+
\int_0^t f(x_g(s,  x_{\text{init}})) - g(x_g(s,  x_{\text{init}})) \; {\rm d}s.
  \end{align*}
  Using the triangular inequality, we get
  \begin{equation}\label{eq:tirang_ineq}
  \begin{split} \|x_f(t, x_{\text{init}}) - x_g(t, x_{\text{init}})\|_2 & \le
\int_0^t \|f(x_f(s,  x_{\text{init}})) - f(x_g(s,  x_{\text{init}}))\|_2\; {\rm d}s 
\\&+ \int_0^t
\|f(x_g(s,  x_{\text{init}})) - g(x_g(s,  x_{\text{init}}))\|_2 \; {\rm d}s.
  \end{split}
  \end{equation}
Because
the function $f$ is $L\text{-Lipschitz}$ over $\Omega$, we have
\[\|f(x_f(s,  x_{\text{init}})) -
f(x_g(s,  x_{\text{init}}))\|_2 \le L \|x_f(s,  x_{\text{init}}) - x_g(s,  x_{\text{init}})\|_2 \quad \forall s \in [0, t].\]
Furthermore, we know that for all
$s \in [0, t]$, $x_g(s,  x_{\text{init}}) \in \Omega$, and therefore 
\[\|f(x_g(s,  x_{\text{init}})) -
g(x_g(s,  x_{\text{init}}))\|_2 \le \metricsup{f-g}{\Omega}.\]
We can hence further bound the left hand side of
\eqref{eq:tirang_ineq} as
  \begin{align*} \|x_f(t, x_{\text{init}}) - x_g(t, x_{\text{init}})\|_2 & \le  L \int_0^t \|x_f(s, x_{\text{init}}) - x_g(s, x_{\text{init}})\|_2 \; {\rm
d}s + t
\metricsup{f-g}{\Omega}.
\end{align*}
  By \Cref{lem:gronwall}, we get
  \begin{equation}\label{eq:ineq1} \|x_f(t, x_{\text{init}}) - x_g(t, x_{\text{init}})\|_2 \le t e^{Lt} \metricsup{f-g}{\Omega}.\end{equation}

Next, we bound the
quantity$\|\dot x_f(t, x_{\text{init}}) -
\dot x_g(t, x_{\text{init}})\|_2$, which can be expressed in
terms of the vector fields $f$ and $g$ as
$ \|f(x_f(t, x_{\text{init}})) - g(x_g(t, x_{\text{init}}))\|_2$.
We have
\begin{equation}\label{eq:ineq2}
\begin{split}
    \|f(x_f(t, x_{\text{init}})) - g(x_g(t, x_{\text{init}}))\|_2 
    &\le
  \|f(x_f(t, x_{\text{init}})) - f(x_g(t, x_{\text{init}}))\|_2 
  \\&+
  \|f(x_g(t, x_{\text{init}})) - g(x_g(t, x_{\text{init}}))\|_2\\
  &\le
  L \|x_f(t, x_{\text{init}}) - x_g(t, x_{\text{init}})\|_2
  + \metricsup{f - g}{\Omega} \\
  &\le
  (1
  +
  Lte^{Lt}) \metricsup{f - g}{\Omega},
\end{split}
\end{equation}
where the first inequality follows from the triangular inequality, the second from the definition of $\metricsup{\cdot}{\Omega}$ and the fact that $f$ is $L\text{-Lipschitz}$ over $\Omega$, and the third one from \eqref{eq:ineq1}.

% and both vectors
% $x_f(t, x_{\text{init}})$ and $x_g(t, x_{\text{init}})$ are in
% $\Omega$, we get
% $\|g(x_f(t, x_{\text{init}})) - g(x_g(t, x_{\text{init}}, x_{\text{init}}))\|_2 \le
% L\|x_f(t,  x_{\text{init}}) - x_g(t,
% x_{\text{init}})\|_2$. Therefore
% \[\|f(x_f(t, x_{\text{init}})) - g(x_g(t, x_{\text{init}}, x_{\text{init}}))\|_2 \le
% \metricsup{f-g}{\Omega} + L\|x_g(t, x_{\text{init}}) - x(t, x_{\text{init}})\|_2.\]

  Putting \eqref{eq:ineq1} and \eqref{eq:ineq2} together, and using the  fact that $t \le T$, we have
  \begin{align*}
       \max &\left\{ \|x_f(t, x_{\text{init}}) - x_g(t, x_{\text{init}})\|_2,  \|f(x_f(t, x_{\text{init}})) - g(x_g(t, x_{\text{init}}))\|_2\right\}
       \\&\quad\quad\quad\quad \le 
  \max\{te^{Lt}, 1
  +
  Lte^{Lt}\}  \metricsup{f - g}{\Omega}
  \\&\quad\quad\quad\quad \le 
  \max\{Te^{LT}, 1
  +
  LTe^{LT}\}  \metricsup{f - g}{\Omega}.
  \end{align*}

  Taking the supremum over $(t,x_{\text{init}}) \in \mathcal S$, we get
  \[\metricdyn{f, g}{\Omega}{T} \le \max\{T e^{LT}, 1+LT e^{LT}\}
    \metricsup{f-g}{\Omega}.\]
\end{proof}

\begin{proof}[Proof of \Cref{thm:approximation_one_side_info}]
  Let us fix a compact set $\Omega \subset \R^n$, a time horizon
  $T > 0$, and a desired accuracy $\varepsilon > 0$. Let
  {$f \in \CD$} be a vector field
  that satisfies any one of the side
  information constraints stated in the theorem. Note that $f$ is $L\text{-Lipschitz}$ over $\Omega$ for some $L \ge 0$. 
  We claim that for any $\delta > 0$, there exists
  a polynomial vector field $p:\R^n \rightarrow \R^n$ that
  satisfies the same side information constraint as $f$ and the inequality
\[\metricsup{f-p}{\Omega} \le \delta.\]
By \Cref{prop:equivalence-metric-vf}, if we take \begin{equation}\label{eq:delta}
\delta = \varepsilon / \max\{T e^{LT}, 1+LT e^{LT}\},
\end{equation}
this shows that there exists a polynomial vector field
$p:\R^n \rightarrow \R^n$ that satisfies the same side
information as $f$ and the inequality 
\[\metricdyn{f,p}{\Omega}{T} \le \varepsilon.\]

We now give a case-by-case proof of our claim above depending on which side
information $f$ satisfies. Throughout the rest of the proof, the constant $\delta$ is fixed as in \eqref{eq:delta}.

{\bf \labelitemi\; Case (i):} Suppose
$f \in \INTERPOLATION{\{(\ix_i,y_i)\}_{i=1}^m}$, where
$x_1, \ldots, x_m \in \Omega$. Without loss of generality, we
assume that the points $x_i$ are all different, or else we
can discard the redundant constraints.  Let $\delta'$ be a
positive constant that will be fixed later. By
\cref{thm:weirstrass}, there exists a polynomial vector field
$q$ that satisfies $\|f - q\|_\Omega \le \delta'$.
% Let
% $v \in \mathbb R^m$ be the vector whose $i^{\text{th}}$
% component is given by
% $v_i = y_i - q(x_i)$, and note that
% $\|v\|_2 \le \delta'$.
%
%
We claim that there exists a
polynomial ${\tilde q}$ of degree $m-1$ such that
\begin{equation}\label{eq:lin_system_poly}
  (q+{\tilde q})(x_i) = y_i \quad i=1,\ldots,m,
\end{equation}
and $\|{\tilde q}\|_\Omega \le C \delta'$, where $C$ is a
constant depending only on the points $x_i$ and the set $\Omega$. Indeed, \eqref{eq:lin_system_poly} can be viewed as a linear system of equations where the unknowns are the coefficients of $\tilde q$ in some basis. For example, if we let   $N={n+m-1 \choose n}$  and ${\tilde q}_{\text{coeff}}\in \R^{N \times n}$ be the matrix whose $j\text{-th}$ column  is the vector of coefficients of $\tilde q_j$ in the standard monomial basis, then \eqref{eq:lin_system_poly} can be written as
% N = num of monomials in n vars of degree m-1 {n+m-1 \choose n}
\begin{equation}\label{eq:AqDelta}A \; {\tilde q}_{\text{coeff}} = \Delta,\end{equation}
where $\Delta \in \R^{m \times n}$ is the matrix whose $i\text{-th}$ row is given by $y_i^T - q(x_i)^T$, and $A \in \R^{m \times N}$  is the
matrix whose $i\text{-th}$ row is
the vector of all standard monomials in $n$ variables and of degree up to $m-1$ evaluated at the point $x_i$. 
One can verify that the rows of
the matrix $A$ are linearly independent (see, e.g., 
\cite[Corollary 4.4]{comon_symmetric_2008}), and so the  matrix $AA^T$ is invertible. If we let $A^+ =A^T(AA^T)^{-1}$, then 
${\tilde q}_{\text{coeff}} = A^+\Delta$ is a solution to \eqref{eq:AqDelta}. Since the matrix $A^+$ only depends on the points $x_i$, and since all the entries of the matrix $\Delta$ are bounded in absolute value by $\delta'$, there exists a constant $c$ such that all the entries of the matrix ${\tilde q}_{\text{coeff}}$ are bounded  in absolute value by  $c \delta'$. Since the set $\Omega$ is compact, there exists a constant $C$ depending only on the points $x_i$ and the set $\Omega$ such that  $\|{\tilde q}\|_\Omega \le C \delta'$.

Finally, by taking $p \coloneqq q+ {\tilde q}$, we get
\(p \in \INTERPOLATION{\{(\ix_i,y_i)\}_{i=1}^m},\)
and
\[ \metricsup{f - p }{\Omega} \le \delta' (1 + C).\]
We take $\delta' = \frac{\delta}{1 + C}$ to conclude the proof for this case.

{\bf \labelitemi\; Case (ii):} Suppose $f \in \SYM{G, \sigma,
  \rho}$, where $G$ is a finite group.  Let
$\delta'$ be a positive constant that will be fixed
later. By \cref{thm:weirstrass}, there exists a polynomial
vector field $p$ that satisfies
$\|f - p\|_\Omega \le \delta'$. Let $p^G: \R^n \rightarrow \R^n$ be the polynomial defined as
\[p^G(x) \coloneqq \frac{1}{|G|} \sum_{g \in G}
  \rho(g^{-1})p(\sigma(g)x) \quad \forall x \in \mathbb
  R^n,\] where $|G|$ is the size of the group $G$.  We claim
that $p^G \in \SYM{G, \sigma, \rho}$.  Indeed, for any
$h \in G$ and $x \in \Omega$, using the fact that $\sigma$ is
a group homomorphism, we get
\[p^G(\sigma(h)x) = \frac{1}{|G|} \sum_{g \in G} \rho(g^{-1})p(\sigma(gh)x).\]
By doing the change of variables $g' = gh$ in the sum above,
and using the fact that $\rho$ is a group homomorphism, we get
\begin{align*}
  p^G(\sigma(h)x)
  &= \frac{1}{|G|} \sum_{g' \in G} \rho(hg'^{-1})p(\sigma(g')x)
  \\&= \frac{1}{|G|} \sum_{g' \in G} \rho(h)\rho(g'^{-1})p(\sigma(g')x).
  \\&= \rho(h)p^G(x).
\end{align*}

We now claim that by by taking $\delta' = \delta
\left(\frac{1}{|G|} \sum_{g \in G}
  \|\rho(g^{-1})\|\right)^{-1}$, where $\|\cdot\|$ denotes the operator norm of its matrix argument, we get $\|f-p^G\|_\Omega \le \delta$. Indeed,
\begin{align*}
  f(x)-p^G(x) &= \frac{1}{|G|} \sum_{g \in G} \left(f(x) - \rho(g^{-1})p(\sigma(g)x)\right)
  \\&= \frac{1}{|G|} \sum_{g \in G} \rho(g^{-1}) \left(\rho(g)f(x) - p(\sigma(g)x)\right)
  \\&= \frac{1}{|G|} \sum_{g \in G} \rho(g^{-1}) \left(f(\sigma(g)x) - p(\sigma(g)x)\right),
\end{align*}
where in the last equation, we used the fact that $f \in  \SYM{G, \sigma, \rho}$.
Therefore,
\begin{align*}
  \|f(x)-p^G(x)\|
  &\le \frac{1}{|G|} \sum_{g \in G} \|\rho(g^{-1})\| \|f(\sigma(g)x) - p(\sigma(g)x)\|_2
  \\&\le \left(\frac{1}{|G|} \sum_{g \in G} \|\rho(g^{-1})\|\right) \delta' = \delta.
\end{align*}

{\bf \labelitemi\; Case (iii):} If $f \in \POS{ \{(P_i, N_i)\}_{i=1}^n }$, where for each $i \in \{1, \ldots, n\}$, the sets $P_i$ and $N_i$ are subsets of $\Omega$ and satisfy $P_i \cap N_i = \emptyset$.
For $i=1,\ldots,n$, let $d_i$ denote the distance between the sets $P_i$ and $N_i$:
\[d_i \coloneqq \min_{x \in P_i, x' \in N_i} \| x - x'\|_2.\]
Since $P_i$ and $N_i$ are compact sets with empty
intersection, the scalar $d_i$ is positive. Fix $\gamma$ to
be any positive scalar smaller than
$\min_{i=1,\ldots,n} d_i$. For $i=1,\ldots,n$, let
\[P_i^\gamma \coloneqq \{ x + \frac \gamma 2 z\; | \; x\in P_i,  z \in \mathbb R^n, \text{
    and } \|z\|_2 \le 1 \},\]
\[N_i^\gamma \coloneqq \{ x +\frac \gamma 2 z\; | \;
  x\in N_i, z \in \mathbb R^n, \text{ and } \|z\|_2
  \le 1 \}.\]
With our choice of $\gamma$,
\(P_i^\gamma \cap N_i^\gamma = \emptyset \text{ for } i=1,\ldots,n.\)
Let $\psi: \mathbb R^n \rightarrow \mathbb R^n$ be the
piecewise-constant function defined as
\[\psi_i(x) = \left\{\begin{array}{cc}
                       1 & \text{if } x \in P_i^\gamma\\
                       -1 & \text{if } x \in N_i^\gamma\\
                       0 & \text{otherwise}
                       \end{array}\right.
                     \quad \text{for } i=1,\ldots,n,\] 
and $\phi^\gamma: \R^n \rightarrow \R$ be the ``bump-like''
function that is equal to $e^{-\frac1{1-\|z\|^2}}$ when
$\|z\|_2 \le \frac \gamma 2$ and $0$ elsewhere. Let $\psi^{\text{conv}}: \mathbb R^n \rightarrow \mathbb R^n$ be the normalized convolution of $\psi$
with $\phi^\gamma$ , i.e.,
\[\psi^{\text{conv}}(x) \coloneqq  \frac1{\int_{z \in \mathbb R^n} \phi^\gamma( z) {\rm d}z} \int_{z \in \mathbb R^n} \psi(x + z) \phi^\gamma(z) {\rm d}z.\]
Note that $\psi^{\text{conv}}$ is a continuous
function as it is the convolution of a piecewise-constant function $\psi$
with a continuous function  $\phi^\gamma$. Moreover, for each $i \in \{1,\ldots,n\}$, $\psi^{\text{conv}}_i$ satisfies
\[\psi^{\text{conv}}_i(x) = 1 \; \forall x \in P_i,\;
\psi^{\text{conv}}_i(x) = -1 \; \forall x \in N_i,\; \text{and }
|\psi^{\text{conv}}_i(x)| \le 1 \; \forall x \in \Omega.\]

Now let $f^\delta \in \CD$ be the vector
field defined component-wise by
\[f^\delta_{i}(x)= f_i(x) + \frac{\delta}{2\sqrt{n}}
  \psi_i^{\text{conv}}(x) \quad i=1,\ldots,n.\]
Note that for $i=1,\ldots,n$, the function $f^\delta_{i}$ is continuous,
bounded below by $\frac{\delta}{2\sqrt{n}}$ on $P_i$, and bounded above by
$\frac{-\delta}{2\sqrt{n}}$ on $N_i$. Moreover
$\|f - f^\delta\|_\Omega \le \frac \delta 2$. 
\Cref{thm:weirstrass} guarantees the existence of a
polynomial vector field $p$ such that
$\|f^\delta - p\|_\Omega \le \frac{\delta}{2\sqrt{n}}$. In particular,
$p \in \POS{ \{(P_i, N_i)\}_{i=1}^n }$ and satisfies
$\|f - p\|_\Omega \le \delta$.

{\bf \labelitemi\; Case (iv):} If
$f \in \MON{ \{(P_{ij}, N_{ij}) \}_{i, j=1}^n}$, where for each $i, j \in\{1, \ldots, n\}$, the
    sets $P_{ij}$ and $N_{ij}$ are subsets of $\Omega$ and belong to different connected components of the graph $\mathcal G(P_{i1},N_{i1}, \ldots, P_{in},N_{in})$. (Recall the definition of this graph from the first paragraph of \Cref{sec:one_side_info}.)

    Consider an index $i \in \{1, \ldots, n\}$, and let
    $C_{i1}, \ldots, C_{ir_i}$ be the connected components of
    the graph
    $\mathcal G(P_{i1},N_{i1}, \ldots, P_{in},N_{in})$. Let $U_{il}$ be the union of the sets
    in component $C_{il}$. Since the sets
    $U_{i1}, \ldots, U_{ir_i}$ are compact and pairwise
    non-intersecting, the minimum distance
    \(d_i \coloneqq \min_{x \in U_{il}, x' \in U_{il'}, l\ne
      \l'} \|x - x'\|_2\) between any two of them is
    positive. Fix $\gamma_i$ to be a positive scalar smaller
    than $d_i$, and for each 
    $l \in \{1, \ldots, r_i\}$, let
  \[U_{il}^{\gamma_i} \coloneqq \{ x + \frac {\gamma_i} 2 z\; | \; x\in U_{il},  z \in \mathbb R^n, \text{
    and } \|z\|_2 \le 1 \}.\]
Define $\psi_i: \mathbb R^n \rightarrow \mathbb R$ to be the
piecewise-linear function defined as
% \text{ for some } $l \in \{1, \ldots, r_i\}
\[
  \psi_i(x) = \left\{\begin{array}{cc}
                       \sum_{j: P_{ij} \in C_l} x_j -
                       \sum_{j: N_{ij} \in C_l} x_j &
                                                      \text{if } x \in U_{il}^{\gamma_i} \text{ for some } l \in \{1, \ldots, r_i\} \\
                       0 & \text{otherwise.}
                   \end{array}\right.
               \]
Let $\psi_i^{\text{conv}}: \R^n \rightarrow \R$  be the normalized convolution of $\psi_i$
with the ``bump-like'' function $\phi^{\gamma_i}: \R^n
\rightarrow \R$ that is equal to $e^{-\frac1{1-\|z\|^2}}$
when $\|z\|_2 \le \frac{\gamma_i}2$ and $0$ elsewhere; that is
\[\psi_i^{\text{conv}}(x) \coloneqq  \frac1{\int_{z \in \mathbb R^n} \phi^{\gamma_i}( z) {\rm d}z} \int_{z \in \mathbb R^n} \psi_i(x + z) \phi^{\gamma_i}(z) {\rm d}z.\]
The function
$\psi_i^{\text{conv}}$
 is continuously differentiable (because $\phi^{\gamma_i}$ is continuously differentiable) and satisfies 
\[\frac{\partial\psi_i^{\text{conv}}}{\partial x_j}(x) = 1 \; \forall x \in P_{ij},\;
\frac{\partial\psi_i^{\text{conv}}}{\partial x_j}(x) = -1 \; \forall x \in N_{ij},\]
\[|\psi^{\text{conv}}_i(x)| \le \sup_{x \in \Omega} |\psi_i(x)| \quad
  \forall x \in \Omega.\]

Now, let $\psi^{\text{conv}} \coloneqq (\psi^{\text{conv}}_1, \ldots, \psi^{\text{conv}}_n)^T$ and
\(f^{\delta'}(x)\coloneqq f(x) + \delta'
  \psi^{\text{conv}}(x)\) for a constant $\delta' > 0$ that will be
    fixed later.
Note that $\|f^{\delta'} - f\|_\Omega \le
\delta' \|\psi^{\text{conv}}\|_\Omega$, and for each pair of
indices $i, j \in \{1, \ldots, n\}$,
\[\frac{\partial f^{\delta'}_i}{\partial x_j}(x) \ge \delta' \; \forall x \in P_{ij},\;
\frac{\partial f^{\delta'}_i}{\partial x_j}(x) \le
-\delta' \; \forall x \in N_{ij}.\]
% \[
% |f^{\delta'}_i(x)| \le \sup_{x \in \Omega} \|x\|_2 \quad
%   \forall x \in \Omega \setminus (P_{ij} \cap N_{ij}).\]

A generalization of the Stone-Weierstrass approximation result stated in
\Cref{thm:weirstrass} to continuously-differentiable
functions (see, e.g.,  \cite{compact_weirstrass_1976})
guarantees the existence of a polynomial $p: \R^n \rightarrow \R^n$ such
that
\[\| f^{\delta'} - p \|_\Omega \le \delta',\quad
\sup_{x \in \Omega} \left| \frac{\partial f_i^{\delta'}}{\partial x_j}(x) - \frac{\partial
    p_i}{\partial x_j}(x) \right| \le \delta'/2 \quad
\forall i, j \in \{1, \ldots, n\}.\]

In particular, $p \in \MON{ \{P_{ij}, N_{ij} \}_{i, j=1}^n}$ and
satisfies $\|f - p\|_\Omega \le \delta'(1+
\|\psi^{\text{conv}}\|_{\Omega})$.
We conclude the proof by taking $\delta' = \frac{\delta}{1+
\|\psi^{\text{conv}}\|_{\Omega}}$.

{\bf \labelitemi\; Case (v):} If
$f \in \INV {\{B_i\}_{i=1}^r}$, where the sets $B_i$ are
subsets of $\Omega$, pairwise nonintersecting, and defined as
$B_i \coloneqq \{ x \in \R^n \; | \; h_{ij}(\ix) \ge 0,\;j
=1,\dots,m_i\}$ for some continuously-differentiable concave
functions $h_{ij}: \R^n \rightarrow \R$ that satisfy
\begin{equation}
     \forall i \in \{1, \ldots, r\}, \; \exists x^i \in B_i \text{ such that }   
     h_{ij}(x^i) > 0  \text{ for } j=1,\ldots,m_i.\label{eq:inv_cdt_thm}
   \end{equation}

%   %
%   Since the set $B_i$ are compact and pairwise
%   nonintersecting, the scalar
%   %
%   \[\gamma \coloneqq \min_{x \in B_i, x' \in B_{i'}, i\ne i'}
%      \|x - x'\|_2\]
%   %
%   is positive. For each $i \in \{1, \ldots,
%   r\}$, let $\psi^i: \R^n \rightarrow \R$ denote the
%   indicator function of the set $B_i$, and let
%   \[B_{i}^{\gamma} \coloneqq \{ x + \frac {\gamma} 2 z\; |
%      \; x\in B_{i}^{\gamma_i}, z \in \mathbb R^n, \text{ and
%      } \|z\|_2 \le 1 \}.\]

%   For each $i \in \{1, \ldots, r\}$, let
%   $\psi_i^{\gamma}: \R^n \rightarrow \R$ be the
%   normalized convolution of $\psi_i$ with the ``bump-like''
%   function $\phi^{\gamma}: \R^n \rightarrow \R$ that is
%   equal to $e^{-\frac1{1-\|z\|^2}}$ when
%   $\|z\|_2 \le \frac{\gamma}2$ and $0$ elsewhere:
% \[\psi_i^{\gamma}(x) \coloneqq  \frac1{\int_{z \in \mathbb R^n} \phi^{\gamma}( z) {\rm d}z} \int_{z \in \mathbb R^n} \psi_i(x + z) \phi^{\gamma}(z) {\rm d}z.\]
   By the same argument as that for {\bf Case (iii)}, for each $i \in \{1, \ldots, r\}$, there exists a continuous function $\psi_i^{\text{conv}}:
   \R^n \rightarrow \R$ that satisfies
   \[ \psi_i^{\text{conv}}(x)
     = 1 \; \forall x \in B_i, \psi_i^{\text{conv}}(x) = 0 \; \forall x \in
     \cup_{i' \ne i} B_{i'}, \; |\psi_i^{\text{conv}}(x)| \le 1 \; \forall x \in \Omega.\]
   Let $\delta' \coloneqq \frac{\delta}{2 r (1+\max_{x,x'\in\Omega} \|x - x'\|_2)}$, and for $i =1\ldots,r$,  let $x^i \in B_i$ be any point satisfying $h_{ij}(x^i) > 0  \text{ for } j=1,\ldots,m_i$. Consider the continuous vector field 
   \[f^{\delta'}(x) \coloneqq f(x) - \delta' \sum_{i=1}^r
     \psi_i^{\text{conv}}(x) (x-x^i).\]
   For every $x
   \in \Omega$, the triangular inequality gives
   \begin{align*}
     \|f(x) - f^{\delta'}(x)\|_2 &\le \delta' \sum_{i=1}^r \|x - x^i\|_2
     \\&\le r \delta' \max_{x' \in \Omega} \|x - x'\|_2 = \frac \delta2,
   \end{align*}
   and so $\|f - f^{\delta'}\|_\Omega \le \delta/2$. Furthermore, for each $i \in \{1, \ldots, r\}$, for each
   $j~\in~\{1, \ldots, m_i\}$, and for each $x \in B_i$
   satisfying $h_{ij}(x) = 0$,
   \begin{equation}\label{eq:bound_f_delta}
   \begin{aligned}
     \langle f^{\delta'}(x), \nabla h_{ij}(x) \rangle
     &= \langle f(x), \nabla h_{ij}(x) \rangle  - \delta' \sum_{k=1}^r \phi_k^{\text{conv}}(x) \langle x - x^k , \nabla h_{ij}(x) \rangle
     \\&= \langle f(x), \nabla h_{ij}(x) \rangle - \delta'  \langle x - x^i , \nabla h_{ij}(x) \rangle
     \\&\ge - \delta' \langle x - x^i ,\nabla h_{ij}(x) \rangle
     \\&\ge \delta' (h_{ij}(x^i) - h_{ij}(x))
    \\&= \delta' h_{ij}(x^i),
   \end{aligned}
   \end{equation}
   %
%   \[ \langle x^i - x ,\nabla h_{ij}(x) \rangle
%      \ge h_{ij}(x^i) - h_{ij}(x)\]
     %
   where the second equality follows from the definition of $\psi_i^{\text{conv}}$ and the fact that $x \in B_i$, the first inequality from the fact that$f
   \in \INV{\{B_i\}_{i=1}^r}$, the second inequality from concavity of the function $h_{ij}$, and the last equality from the fact that $h_{ij}(x) = 0$.

   For a constant $\delta'' > 0$ that will be fixed later,
   \Cref{thm:weirstrass} guarantees the existence of a
   polynomial vector field $p$ such that
   $\|f^{\delta'} - p\|_\Omega \le \delta''$. By triangular inequality we have
   $\|f - p\|_\Omega \le \delta/2 + \delta''$. Furthermore, for each $i \in \{1, \ldots, r\}$, for each
   $j \in \{1, \ldots, m_i\}$, and for each $x \in B_i$
   satisfying $h_{ij}(x) = 0$, we have
\begin{align*}
  \langle p(x), \nabla h_{ij}(x) \rangle
  &= \langle p(x) - f^{\delta'}(x), \nabla h_{ij}(x) \rangle  + \langle f^{\delta'}(x), \nabla h_{ij}(x) \rangle
  \\&\ge - \delta'' \|\nabla h_{ij}(x)\|_2  + \delta' h_{ij}(x^i)
   \end{align*}
   due to \eqref{eq:bound_f_delta} and the Cauchy-Schwarz
   inequality. Let
%   \[\delta' = \frac{\delta}{r \max_{x,x'\in\Omega} \|x - x'\|_2}\]
%   and
   \[\delta'' \coloneqq \min \left\{ \frac \delta 2, \min_{i \in \{1,\ldots,r\}} \min_{j \in \{1, \ldots, m_i\}, x \in B_i} \delta' \frac{ h_{ij}(x^i)}{\|\nabla h_{ij}(x)\|_2}\right\},\]
   and note that $\delta'' > 0$ as we needed before because $h_{ij}(x^i) > 0$ for each $i \in \{1,\ldots,r\}$ and $j \in \{1, \ldots, m_i\}$. With this choice of $\delta''$,
   we get that $p \in \INV {\{B_i\}_{i=1}^r}$ and $\|f -
     p\|_\Omega \le \delta$.

     {\bf \labelitemi\; Case (vi):} If $f \in \GRAD$. In this
     case, there exists a continuously-differentiable function
     $V: \R^n \rightarrow \R$ such that
     $f(\ix) = -\nabla V(\ix)$. A generalization of the
     Stone-Weierstrass theorem to continuously-differentiable
     functions (see, e.g., \cite{compact_weirstrass_1976})
     guarantees the existence of a polynomial
     $W: \R^n \rightarrow \R$ such that
     \[\max_{x \in \Omega} \left\| \nabla V(x) - \nabla W(x) \right\|_2 \le \delta.\]
     Letting $p(x) = -\nabla W(x)$, we get that $p \in \GRAD$ and $\|f -
     p\|_\Omega \le \delta$.

     {\bf \labelitemi\; Case (vi'):} If $f \in \HAMILTONIAN$. The
     proof for this case is analogous to {\bf Case (vi)}.
\end{proof}

% \vspace{-\topsep}
\subsection{Approximating a vector field while approximately satisfying multiple side information constraints}
\label{sec:multiple_side_info}
It is natural to ask whether \cref{thm:approximation_one_side_info} can be
generalized to allow for polynomial approximation of vector
fields satisfying \emph{multiple}  side information constraints. It turns out that our proof idea of ``smoothing by convolution" %technique from the proof of that theorem 
can be used to show that the answer is positive if  the
following three conditions hold: (i) the side information
constraints are of type $\NINTERPOLATION,  \NPOS, \NMON,
\NINV$, or $\NSYM$, (ii) each side information constraint
satisfies the assumptions of
\Cref{thm:approximation_one_side_info}, and (iii) the regions
of the space where the first four types of side information constraints are imposed are pairwise nonintersecting. 
% The intuition is that under condition (iii), our ``smoothing by convolution" technique turns multiple vector fields that each satisfy  a particular side information constraint into a single smooth vector field that simultaneously satisfies all of them.
In absence of condition (iii),  the answer is no longer positive as the next example shows.

\begin{myexp}
  Consider the univariate vector field
  $f: \R \rightarrow \mathbb R$ given by

  \[ f(x) \coloneqq \begin{cases} 
      0 & x \ge 0 \\
      -e^{-\frac1{x^2}} & x < 0.
    \end{cases}
  \]
  This vector field is continuously differentiable over $\R$ and
  satisfies the following combination of side information
  constraints:
  \begin{equation}
  \INTERPOLATION{\{(0, 0), (1, 0)\}} \text{ and }
  \MON{ \{([-1, 1], \emptyset)\}}.\label{eq:side_info_f}
\end{equation}
In other words, $f$ is nondecreasing on the interval
$[-1, 1]$ and satisfies $f(0) = f(1) = 0$. Yet, the only
polynomial vector field that satisfies the constraints in
\cref{eq:side_info_f} is the identically zero polynomial.  As
a result, the vector field $f$ cannot be approximated
arbitrarily well over $[-1, 1]$ by polynomial vector fields
that satisfy the side information constraints in
\cref{eq:side_info_f}.
\end{myexp}

\begin{table}
  \small
  \centering
        \begin{tabular}{ |c|c| }
          \hline
          {\bf Side information $S$} & {\bf Functional $\LS(f)$}  \\
          \hline
          \makecell{$\INTERPOLATION{\{(\ix_i,y_i)\}}_{i=1}^m)$\\
          with $x_i \in \Omega$\\for $i=1\ldots,m$} & $\displaystyle \max_{i=1,\ldots,m}\|f(x_i) - y_i\|_2$\\\hline
          \makecell{$\SYM{G,\sigma,\rho}$\\} & $ \displaystyle  \max_{g \in G} \max_{\substack{i=1,\ldots,n\\x\in\Omega}} |f_i(\sigma(g)x) - (\rho(g)f(x))_i|$\\\hline
          \makecell{$\POS{ \{(P_i, N_i)\}_{i=1}^n }$\\
          with $P_i, N_i \subseteq \Omega$\\
          for $i=1,\ldots,n$}& $\displaystyle \max_{i=1,\ldots,n} \max\left\{0, \max_{x \in P_i} -f_i(x), \max_{x \in N_i} f_i(x)\right\} $\\\hline
          \makecell{$\MON{ \{(P_{ij}, N_{ij}) \}_{i, j=1}^n}$\\
          with $P_{ij}, N_{ij} \subseteq \Omega$\\
          for $i,j=1,\ldots,n$} & $ \displaystyle\max_{i,j=1,\ldots,n} \max\left\{0, \max_{x \in P_{ij}} -\frac{\partial f_i}{\partial x_j}(x), \max_{x \in N_{ij}} \frac{\partial f_i}{\partial x_j}(x)\right\}$\\\hline
          \makecell{$\INV {\{B_i\}_{i=1}^r}$ where\\ $B_i \coloneqq \{x \; | \; h_{ij}(x) \ge 0$\\ $\; \forall j \in \{1,\ldots,m_i\}\} \subseteq \Omega$\\ for $i=1,\ldots,r$}  & $\displaystyle  \max_{i=1\ldots,r} \max_{\substack{\ix \in B_i\\ j\in\{1,\ldots,m_i\}\\ h_{ij}(\ix) = 0}} \max\left\{0, - \langle f(\ix), \nabla h_{ij}(\ix) \rangle\right\}$\\\hline
          $\GRAD$ &  \makecell{$\displaystyle\inf_{\substack{V:\R^n \rightarrow \mathbb R}} \;\max_{\substack{i=1,\ldots,n\\x\in\Omega}} \left| f_i(x) + \frac{\partial V}{\partial x_i}(x)\right|$}\\
          \hline
          $\HAMILTONIAN$ & \makecell{$ \displaystyle \inf_{\substack{H:\R^n \rightarrow \mathbb R}}
          \max_{\substack{(p,q) \in \Omega, \\i=1\ldots,n/2}} \max \left\{\left| f_i(p, q) + \frac{\partial H}{\partial q_i}(p,q)\right|,\right.$\\ $\quad\quad\quad\quad\quad\quad\quad\left.\left| f_{i+n/2}(p, q)  -\frac{\partial H}{\partial p_i}(p,q)\right|\right\}$ }\\
          \hline
      \end{tabular}
      
  \caption{\label{tab:linear_L_S}For each side information
    $S$, the functional $\LS: \CD \rightarrow \R$ quantifies
    how close a vector field $f \in \CD$ is to satisfying
    $S$.}
  \vspace{-5mm}
\end{table}
To overcome difficulties associated with such examples,  we introduce
% We mitigate this shortcoming of polynomial vector fields by
% introducing 
the notion of approximate satisfiability of side information
over a compact set $\Omega \subset \R^n$. Before we give a
formal definition of this notion, for each side information
constraint $S$, we present in \Cref{tab:linear_L_S} a
functional $\LS: \CD \rightarrow \R$ that measures how close
a vector field in $\CD$ is to satisfying the side information
$S$. One can verify that the functional $\LS$ has the
following two properties: (i)~for any vector field
$f \in \CD$,
\begin{equation}\label{eq:L_S_0}
 \LS(f) = 0 \text{ if and only if } f \text{ satisfies } S,
\end{equation}
and (ii) for any $\delta > 0$, there exists
$\gamma > 0$, such that for any two vector fields $f,\hat f \in \CD$,
\begin{equation}\label{eq:L_S_continuous}
  \|f - {\hat f}\|_\Omega \le \gamma \text{ and } \max_{\substack{x \in \Omega,\\ i,j=1,\ldots,n}} \left|\frac{\partial f_i}{\partial x_j}(x) - \frac{\partial {\hat f}_i}{\partial x_j}(x)\right| \le \gamma
  \implies |\LS(f) - \LS({\hat f})| \le \delta.
\end{equation}
Indeed, take e.g. $S=\INV {\{B_i\}_{i=1}^r}$, where
$B_i \coloneqq \{x \in \R^n\; | \; h_{ij}(x) \ge 0,\;
j=1,\ldots,m_i\}$. It is clear from condition
\cref{eq:eq-invariance} that $\LS(f) = 0$ if and only if
$f \in \INV {\{B_i\}_{i=1}^r}$.  To verify the second
property, let $\delta > 0$ be given. If we take
\[\gamma = \delta \left\{  \max_{\substack{x \in \Omega,\\ i=1\ldots,r\\j=1\ldots,m_i}}\|\nabla h_{ij}(x)\|\right\}^{-1},\]
it is easy to see that for any two vector fields
$f, \hat f \in \CD$ satisfying
\(\|f - {\hat f}\|_\Omega \le \gamma\), we must have
$|L_S(f) - L_S(\hat f)| \le \delta$. Indeed, let
$i \in \{1, \ldots, r\}$ and $x \in B_i$ be such that
$h_{ij}(x) = 0$ for some $j \in \{1,\ldots,m_i\}$. Then, the
Cauchy-Schwarz inequality and our choice of $\gamma$ give
\[|\langle f(x), \nabla h_{ij}(x)\rangle - \langle \hat f(x),
  \nabla h_{ij}(x)\rangle| \le \|f - \hat f\|_\Omega
  \|\nabla h_{ij}(x)\| \le \delta.\] The desired result follows by taking the maximum over $i$, $j$, and $x$.

% \vspace{-\topsep}
\begin{mydef}[$\delta\text{-satisfiability}$]Let $\Omega \subset \R^n$ be a compact set and consider any side information $S$ presented in
  \Cref{tab:linear_L_S} together with its corresponding functional $\LS$. For a scalar   $\delta >0$,  we say that a vector field $f \in \CD$ $\delta\text{-satisfies}$ $S$ if
  $\LS(f) \le \delta$.
\end{mydef}
% \vspace{-\topsep}

% \vspace{-\topsep}
% \begin{myexp}
%   A vector field $f$ $\delta\text{-satisfies}$ the side
%   information $\INTERPOLATION{\{(\ix_i,y_i)\}_{i=1}^m}$ if
%   \(\|f(x_i) - y_i\|_2 \le \delta\) for \(i=1,\ldots,m\), and
%   $\delta\text{-satisfies}$ the side information
%   $\POS{ \{(P_i, N_i)\}_{i=1}^n }$ if
%   $f_i(\ix) \ge -\delta \; \forall \ix \in P_i$ and
%   $f_i(\ix) \le \delta \; \forall \ix \in N_i$ for
%   \(i=1,\ldots,n\).
% \end{myexp}
% \vspace{-\topsep}

From a practical standpoint,
for small values of $\delta$,
it is reasonable to substitute the requirement of exact satisfiability of side information for  
$\delta\text{-satisfiability}$. This is especially true since most optimization solvers return an
approximate numerical solution anyway. The following theorem shows that
polynomial vector fields can approximate any continuously-differentiable vector field
$f$ and satisfy the same side information as $f$ up to an
arbitrarily small error tolerance $\delta$.
It also shows that in the context of learning a vector field from trajectory data, one can always impose  $\delta\text{-satisfiability}$ on a candidate polynomial vector field via semidefinite programming.

% \vspace{-\topsep}
\begin{theorem}\label{thm:poly_approx}
  For any compact set $\Omega \subset \R^n$, time horizon
  $T > 0$, desired approximation accuracy $\varepsilon > 0$, desired side information satisfiability accuracy
  $\delta > 0$, and for any vector field $f \in \CD$ that satisfies
  any combination of the side information constraints from the
  first column of \Cref{tab:linear_L_S}, there exists a
  polynomial vector field $p: \mathbb R^n \rightarrow \R^n$
  that $\delta\text{-satisfies}$ the same combination of side
  information as $f$ and has
  $\metricdyn{f, p}{\Omega}{T} \le \varepsilon$. 
  
  Moreover, if the set $\Omega$, the sets $P_i,N_i$ in the
  definition of $\POS{ \{(P_i, N_i)\}_{i=1}^n }$, the sets
  $P_{ij},N_{ij}$ in the definition of
  $\MON{ \{P_{i j}, N_{i j}\}_{i, j=1}^n}$, and the sets
  $B_i$ in the definition of $\INV{\{B_i\}_{i=1}^r}$) are all
  closed basic semialgebraic and their defining polynomials
  satisfy the Archimedian property, then
  $\delta\text{-satisfiability}$ of all side information
  constraints by the polynomial vector field $p$ has a sum
  of squares certificate.
  % comes with a $\text{degree-}r$ sos
  % certificate for some integer $r$ large enough whenever the
  % sets used to define these side information constraints
  % (i.e., the sets $P_i$, $N_i$, $P_{ij}$, $N_{ij}$, and
  % $B_i$) are closed semialgebraic sets and are given by a
  % list of polynomials that satisfy the Archimedian property.
\end{theorem}
\begin{proof}

Let $f \in \CD$  satisfy a list $S_1, \ldots, S_k$ of side information constraints from the first column of \Cref{tab:linear_L_S}, and let the scalars $T$, $\varepsilon,\delta>0$ be fixed.
A generalization of the Stone-Weierstrass approximation theorem to continuously-differentiable
functions (see, e.g.,  \cite{compact_weirstrass_1976})
guarantees that for any $\gamma > 0$, there exists a polynomial $p^\gamma: \R^n \rightarrow \R^n$ such
that
\begin{equation}
  \|f - p^\gamma\|_\Omega \le \gamma \text{ and } \max_{\substack{x \in
    \Omega\\ i,j=1,\ldots,n}} \left|\frac{\partial
      f_i}{\partial x_j}(x) - \frac{\partial p^\gamma_i}{\partial
      x_j}(x)\right| \le \gamma.\label{eq:p_gamma}  
\end{equation}
For the rest of this paragraph, for any $\gamma > 0$, we fix
an (arbitrary) choice for the polynomial~$p^\gamma$. 
%For the rest of this paragraph, we think of $p^\gamma$ as any
%fixed polynomial that satisfies~\eqref{eq:p_gamma}.
Since for each $i \in \{1,\ldots,k\}$, the functional
$\LiS{i}$ satisfies \eqref{eq:L_S_continuous}, there exists a
scalar $\gamma_i>0$ for which
$\LiS{i}(p^{\gamma}) \le \delta/2$ for any
$\gamma \in (0, \gamma_i]$. If we let
\[\bar \gamma \coloneqq \min \{ \varepsilon / \max\{T e^{LT}, 1+LT e^{LT}\}, \gamma_1, \ldots, \gamma_k\}, \]
 where $L > 0$ is any scalar for which $f$ is $L\text{-Lipschitz}$ over $\Omega$, then the polynomial $p \coloneqq p^{\bar \gamma}$ $\delta/2\text{-satisfies}$ $S_1,\ldots,S_k$ (and hence $\delta\text{-satisfies}$ $S_1,\ldots,S_k$), and because of \Cref{prop:equivalence-metric-vf}, also satisfies $\metricdyn{f, p}{\Omega}{T} \le \varepsilon$.

% Let us now assume that the set $\Omega$, the sets
% $P_i,N_i$ in the definition of
% $\POS{ \{(P_i, N_i)\}_{i=1}^n }$, the sets $P_{ij},N_{ij}$ in
% the definition of $\MON{ \{P_{i j}, N_{i j}\}_{i, j=1}^n}$,
% and the sets $B_i$ in the definition of
% $\INV{\{B_i\}_{i=1}^r}$) are all closed basic semialgebraic
% and are given by a list of polynomials that satisfy the
% Archimedian property, and let us show that
% $\delta\text{-satisfiability}$ of these side information
% constraints by a polynomial vector field $p$ can be certified
% via an sos certificate.  The following result, which is a
% direct consequence of Putinar's approximation theorem in
% \cref{th:putinar}, will be useful to us: A constraint of the
% type
% \begin{equation}
% h_1(x) \ge 0 \ldots, h_s(x) \ge 0 \Rightarrow q(x) > 0
%   \quad \forall x \in \SETALG,\label{eq:valid_sdp_const}
% \end{equation}
% where $q, h_1, \ldots, h_s:\R^n \rightarrow \R$ are
% polynomial functions whose coefficients are linear combinations of the coefficients
% $p$ (and possibly, on other decision variables),
% $\SETALG = \{g_1 \ge 0, \ldots, g_r \ge 0\}$ is a (fixed)
% closed semialgebraic set, and $\{g_1,\ldots,g_r\}$ is a list
% of polynomial functions that satisfy the Archimedian
% property, admits an sos certificate.

To prove the second claim of the theorem, 
observe that for each $\ell \in \{1, \ldots, k\}$,  the fact that $p$ $\delta/2\text{-satisfies}$ $S_{\ell}$ implies the following inequalities:\footnote{We exclude the case $S_\ell = \INTERPOLATION{\{(\ix_i,y_i)\}_{i=1}^m}$ because verifying $\delta\text{-satisfiability}$ is trivial there, and the case $S_\ell = \HAMILTONIAN$ because the argument for it is similar to that of $S_\ell = \GRAD$.}

\begin{itemize}

% \item The condition that $p$ $\delta\text{-satisfies}$ $\INTERPOLATION{\{(\ix_i,y_i)\}_{i=1}^m}$, where
% $x_1, \ldots, x_m \in \Omega$ is equivalent to 
% \[\|p(x_i) - y_i\|_2 \le \delta \text{ for } i=1,\ldots,m,\]
% and can be directly checked.

\item If $S_\ell = \SYM{G,\sigma,\rho}$, 
\[ p_i(\sigma(g)x) - (\rho(g)p(x))_i + \delta > 0 \text{ and }
  (\rho(g)p(x))_i - p_i(\sigma(g)x) + \delta > 0 \; \forall x
  \in \Omega,\]
% for $G \in \{G_1, \ldots, G_r\}$,
for $g\in G$ and $i=1,\ldots,n$;

\item If $S_\ell = \POS{ \{(P_i, N_i)\}_{i=1}^n }$, 
\[ p_i(x) + \delta > 0\; \forall x
      \in P_{i} \text{ and } -p_i(x) + \delta > 0 \; \forall x
      \in N_{i},\]
for $i=1,\ldots,n$;

  \item If $S_\ell =  \MON{ \{(P_{ij}, N_{ij}) \}_{i, j=1}^n},$
\[\frac{\partial p_i}{\partial x_j}(x) + \delta >0 \; \forall x
  \in P_{ij}, \text{ and } -\frac{\partial p_i}{\partial x_j}(x) + \delta > 0\; \forall
  x \in N_{ij},\]
for $i,j=1,\ldots,n$;
\item If $S_\ell = \INV {\{B_i\}_{i=1}^r}$,
\[
 \langle p(\ix), \nabla h_{ij}(\ix)
  \rangle + \delta > 0\; \forall x \in
    B_i \cap \{x \in \R^n \; | \; h_{ij}(\ix) = 0\},\]
for $i=1\ldots,r$, $j = 1,\ldots,m_i$;

\item If $S_\ell = \GRAD$,
\[ p_i(x) + \frac{\partial V}{\partial x_i}(x) + \delta > 0 \text{ and }  -p_i(x) - \frac{\partial V}{\partial x_i}(x) + \delta > 0 \quad \forall x \in \Omega,\]
for $i=1\ldots,n$, where $V:\R^n \rightarrow \R$ is a polynomial function. (The fact that $V$ can be taken to be a polynomial function follows from an other application of the generalization of Stone-Weierstrass approximation theorem that was used at the beginning of this proof.)
\end{itemize}

%          

%
% The polynomial vector field $p$ $\delta\text{-satisfies}$ $\SYM{G, \rho, \sigma}$ if and only if
% \[\|p(\sigma(g)x) - \rho(g)p(x)\|_2 \le \gamma \quad \forall x \in \Omega \text{ for } g \in G.\]
% %

% The polynomial vector field $p$ $\delta\text{-satisfies}$
% $\POS{ \{(P_i, N_i)\}_{i=1}^n }$,
% $\MON{ \{(P_{ij}, N_{ij}) \}_{i, j=1}^n}$, or
% $\INV {\{B_i\}_{i=1}^r}$ if and only if a constraint of the
% type \eqref{eq:valid_sdp_const} is satisfied.

% The polynomial vector field $p$ $\delta\text{-satisfies}$ $\GRAD$ (resp. $\HAMILTONIAN$) 

Observe that each of the above inequalities states that a
certain polynomial is positive over a certain closed basic
semialgebraic set whose defining polynomials satisfy the
Archimedian property by assumption. Therefore, by Putinar's
Positivstellesatz (\Cref{th:putinar}), there exists a
nonnegative integer $d$ such that each one of these
inequalities has a $\text{degree-}d$ sos-certificate (see
\cref{eq:putinar-certificate}). Therefore,
$\delta\text{-satisfiability}$ of each side information
$S_1, \ldots, S_k$ by the vector field $p$ can be proven by a sum of
squares certificate.
\end{proof}
\newcommand{\SDPD}{$\text{SDP}_d$ }
% varepsilon
\newcommand{\eeps}{{\frac{\varepsilon}{4\max\{T e^{LT}, 1+LT e^{LT}\}}}}

\bachir{An interesting corollary of \cref{thm:poly_approx} is that given a list of side information that an unknown vector field $f$ is known to satisfy, and a training dataset that is large enough, one can find  a polynomial vector field $p$, which approximately satisfies the same combination of side information and is close to $f$,
  by solving a finite sequence of semidefinite programs. 
  More formally, consider a  compact set $\Omega \subset \R^n$, time horizon
  $T > 0$, desired approximation accuracy $\varepsilon > 0$, desired side information satisfiability accuracy
  $\delta > 0$, and  a vector field $f \in \CD$ that satisfies
  any combination of the side information constraints from the
  first column of \Cref{tab:linear_L_S}. Let $$L \coloneqq \sup_{x \in \Omega} \|J_f(x)\|,$$
  where $J_f$ is the Jacobian of the vector field $f$.
  Let $D$ be any finite subset of $\Omega$ that satisfies
  \[\sup_{x \in \Omega} \min_{x_i \in D} \|x - x_i\| \le \frac{\varepsilon}{4(L+1)\max\{T e^{LT}, 1+LT e^{LT}\}}. \]
  Furthermore, assume that
  the set $\Omega$, the sets $P_i,N_i$ in the
  definition of $\POS{ \{(P_i, N_i)\}_{i=1}^n }$, the sets
  $P_{ij},N_{ij}$ in the definition of
  $\MON{ \{P_{i j}, N_{i j}\}_{i, j=1}^n}$, and the sets
  $B_i$ in the definition of $\INV{\{B_i\}_{i=1}^r}$) are all
  closed basic semialgebraic and their defining polynomials
  satisfy the Archimedian property.
 Consider the sequence of semidefinite programs \SDPD indexed by a nonnegative integer $d$:\footnote{Recall that $\mathcal P_d$ denotes the set of polynomial vector fields of degree $d$, and see \cref{eq:putinar-certificate} for the notion of a  $\text{degree-}d$  sos-certificate.}
\begin{equation*}
\begin{aligned}
\min_{p \in \mathcal P_d } \quad & \max_{x_i \in D} \|p(x_i) - f(x_i)\| \\
\textrm{s.t.} \quad & p \ \text{ has  a $\text{degree-}d$  sos-certificate } \\
&    \text{ of  the inequality } \|J_p(x)\| \le L + 1 \; \forall x \in \Omega \\
& \text{ as well as of $\delta\text{-satisfiability}$ of the side information of $f$}.
\end{aligned}
\end{equation*}
 %
%  the decision variable of \SDPD is a polynomial $p$ of degree $d$, its objective function is given by
%  \[\min_p \max_{x_i \in D} \|p(x_i) - f(x_i)\|,\]
%  and its constraints are given by the $\text{degree-}d$  sos-certificate of $\delta\text{-satisfaction}$ of the side information of $f$ (see \cref{eq:putinar-certificate}), in addition to the $\text{degree-}d$ sos certificate  of the inequality
%  \[\|J_p(x)\| \le L + 1 \; \forall x \in \Omega.\]
 We claim that for $d$ large enough, \SDPD is feasible with optimal value at most 
 \[  \varepsilon'  \coloneqq \eeps, \] 
 and that any of its feasible solutions $\tilde p$ with objective value at most $2\varepsilon'$ would $\delta\text{-satisfy}$ the same side information as $f$, and have $\metricdyn{f, \tilde p}{\Omega}{T} \le \varepsilon$.}
 \bachir{
We know from \cref{thm:poly_approx} (and its proof) that for $d$ large enough, there exists a polynomial vector field $\hat p$ of degree $d$ that $\delta\text{-satisfies}$ the same side informationas $f$ with $\text{degree-}d$ sos-certificates of $\delta\text{-satisfaction}$
and has
\begin{equation}
 \max_{\substack{x \in
\Omega\\ i,j=1,\ldots,n}} \left|\frac{\partial
  f_i}{\partial x_j}(x) - \frac{\partial  \hat p_i}{\partial
  x_j}(x)\right| \le \frac1{2n}\text{ and } \metricdyn{f, \hat p}{\Omega}{T} \le \varepsilon'.
\end{equation}
 In particular, $\hat p$ is feasible to \SDPD (because of the first inequality) and has objective value no larger than $\varepsilon'$ since $\|\hat p(x) - f(x) \| \le \metricdyn{f, \hat p}{\Omega}{T} \; \forall x \in \Omega$.}

 \bachir{Let $\tilde p$ be a feasible solution to \SDPD (for some $d \in \mathbb N$) with objective value at most $2\varepsilon'$. Then, $\tilde  p$ $\delta\text{-satisfies}$ the same side information as $f$ by construction. Moreover, for any $x \in \Omega$ and $x_i \in D$, the triangular inequality gives
 \[ \|\tilde p(x) - f(x)\| \le \|\tilde p(x) - \tilde p(x_i)\| + \|\tilde p(x_i) - f(x_i)\| + \|f(x_i) - f(x)\|.\]
  Observe that 
   \(\|\tilde p(x_i) - f(x_i)\| \le 2\varepsilon'.\)
   Moreover, since $\tilde p$ and $f$ are $L+1$ Lipschitz, we have 
  \[\|\tilde p(x) - \tilde p(x_i)\| \le (L+1) \| x - x_i\| \text{ and } \|f(x) - f(x_i)\| \le (L+1) \| x - x_i\|. \]
  Therefore,
   \[\|\tilde p(x) - f(x)\| \le 2  (L + 1) \| x - x_i\| + 2\varepsilon' \; \forall x_i \in D,\; \forall x \in \Omega.\]
   It follows that 
   \[\|\tilde p(x) - f(x)\| \le  \frac{\varepsilon}{\max\{T e^{LT}, 1+LT e^{LT}\}}  \;  \forall x \in \Omega,\]
   \cref{prop:equivalence-metric-vf} therefore gives
  $\metricdyn{f, \tilde p}{\Omega}{T}  \le \varepsilon$.}
 
  \bachir{
  We remark that the SDP construction depends on the possibly unknown Lipschitz constant $L$ of the vector field $f$. However, any upper bound on $L$ would suffice for the SDP construction and the convergence guarantee established above. 
Note that without an upper bound on $L$ to be imposed on the Lipschitz constant of our candidate vector fields, no algorithm could learn the vector field $f$ based on the values that it takes on a finite dataset $D$ alone. Indeed, for any finite training dataset $D$, there always exists another vector field $\tilde f \in C_1^\circ(\Omega)$ that is equal to $f$ on $D$, and arbitrarily far from $f$ outside of $D$.}

\section{Discussion and future research directions}
\label{sec:conclusion}

From a computational perspective, our approach to learning
dynamical systems from trajectory data while leveraging side
information relies on convex optimization. If the side
information of interest is $\NINTERPOLATION$, $\NSYM$,
$\GRAD$, or $\NHAMILTONIAN$, then our approach leads to a
least-squares problem, and thus can be implemented at large
scale. For side information constraints of $\NPOS$, $\NMON$,
or $\NINV$, our approach requires solutions to semidefinite
programs. Classical interior-point methods for SDP come with
polynomial-time solvability guarantees (see
e.g.~\cite{vandenberghe1996semidefinite}), and in practice
scale to problems of moderate sizes. In the field of
dynamical systems, many applications of interest involve a
limited number of state variables, and therefore our approach
to learning such systems leads to semidefinite programs that
off-the-shelf interior-point method solvers can readily
handle. For instance, each semidefinite program that was
considered in the numerical applications of
\Cref{sec:applications} was solved in under a second on a
standard personal machine by the solver MOSEK~\cite{mosek}.
% A recurrent criticism leveled
% against semidefinite programming is its lack of
% scalability.
An active and exciting area of research is focused on
developing algorithms for large-scale semidefinite programs
(see e.g. \cite{sdp_survey_2019,opt_sdp_cornell_2019}), and
we believe that this effort can extend our learning approach
to large-scale dynamical systems.

The size of our semidefinite programs is also affected by the
degree of our candidate polynomial vector field and the
degrees of the sos multipliers in
\eqref{eq:putinar-certificate} that result from the
application of Putinar's Positivstellesatz. In practice,
these degrees can be chosen using a statistical model
validation technique, such as cross validation. \bachir{For example,
one can split the available data into training and testing, use the training data for learning a vector field of degree $d \in \{1, 2, \ldots \}$, and choose the degree that achieves the lowest generalization error on the test data.}

These
techniques take into account the fact that lower degrees can
sometimes have a model regularization effect and lead to
better generalization on unobserved parts of the state space.

We end by mentioning some questions that
are left for future research.

\setlist{leftmargin=4.5mm}
\begin{itemize}
\item While the framework presented in this paper deals with
  continuous-time dynamical systems, we believe that most of
  the ideas could be extended to the discrete-time
  setting. It would be interesting to see how the definitions
  of side information, the approximation results, and the
  computational aspects contrast with the continuous-time
  case. Extending our framework to the problems of learning
  partial differential equations and stochastic differential
  equations with side information would also be
  interesting research directions.
\item We have shown that for any $\delta > 0$, polynomial
  vector fields can approximate to arbitrary accuracy any
  vector field $f \in \CD$ while $\delta\text{-satisfying}$
  any list of side information that $f$ is known to satisfy.
  Even though from a practical standpoint,
  $\delta\text{-satisfiability}$ is sufficient (when $\delta$
  is small), it is an interesting mathematical question in
  approximation theory to see which combinations of side
  information can be imposed exactly on polynomial vector
  fields while preserving an arbitrarily tight approximation
  guarantee to functions in $\CD$.
\item We have presented a list of six types of
  side information that arise naturally in many
  applications and that lead to a convex formulation (meaning that a convex combination of two vector fields that satisfy any one of the six side information constraints will also satisfy the same side information constraint).
  There are of course other
  interesting side information constraints that do not lead to a
  convex formulation. Examples include the knowledge that an equilibirum point is locally
  or globally stable/stabilizable, and the knowledge that trajectories of the system starting in a set $A\subseteq \R^n$
  avoid/reach another set $B\subseteq \R^n$. It is an interesting research
  direction to extend our approximation results and our
  sos-based approach to handle some of these nonconvex side
  information constraints.
\item Finally, from a statistical and information-theoretic
  point of view, it is an interesting question to quantify
  the benefit of a particular side information constraint in reducing the number of trajectory observations needed to learn a good approximation of the unknown vector field.
  % A related question in the context of learning
%   dynamical systems from trajectory observation and side
%   information is as follows: how should one decide if given
%   the choice between the knowledge of a new side information
%   constraint and the observation of a new trajectory?
\end{itemize}

\vspace{.4cm}

{\bf Aknowledgments:} The authors are grateful to two anonymous referees, Charles Fefferman, Georgina Hall, Frederick Leve, Clancey Rowley, Vikas Sindhwani, and Ufuk Topcu for insightful questions and
comments.

\newpage
\appendix
%\color{blue}
\section{Additional numerical experiments}
\subsection*{Description of the experiments}
% use notation N_S = 20, N_T = 
For the first and last learning experiments of this paper (\Cref{sec:diffusion_contagious_disease,sec:lorenz}), we conduct additional numerical tests here to show the effects of varying the degree of the learned polynomial vector field, the noise level, and the number of trajectories used in learning.

For the first (resp. last) learning application, we simulate $N_t$ trajectories starting from initial conditions picked uniformly at random from the box $B = [0, 1]^2$ (resp. $B = [-10, 10]^3$) up to time $T = 30$, and collect $20$ noisy samples 
from each trajectory at times evenly spaced on the interval $[0, T]$. Each sample has the form
$(x, f(x) + \varepsilon)$, where $\varepsilon$ is a $2 \times 1$ (resp. $3 \times 1$) Gaussian variable with mean $0$ and covariance $\sigma^2 I$, with $\sigma= 10^{-3}, 10^{-2}, 10^{-1}$. The noise terms added to samples are independent from each other.
We experiment with learning polynomial vector fields $p(x)$ with degree $d = 2,3,4$ either with no side information, or with the three side information constraints presented in \Cref{sec:diffusion_contagious_disease,sec:lorenz}.

For each combination of $N_t, \sigma, d$, we report the following metrics related to the learned vector field $p$:
\begin{itemize}
    \item {\bf The training error:} The average value over all training sample points $(x, f(x)+\varepsilon)$ of the quantity $\|p(x) - (f(x)+\varepsilon)\|$.
    \item {\bf The vector field test error:}
    The average value of  $\|p(x) - f(x)\|$, where $x$ runs over a regular discretization of the box $B$, where each dimension has size $10$. Therefore, the test set includes 100 points in our two-dimensional example and 1000 points in our three-dimensional one. 
    \item {\bf The trajectory test error:}
    The average value of  $\|x_f(t_i; x_0) -  x_p(t_i; x_0)\|$, where $x_f(t; x_0)$ (resp. $x_p(t; x_0)$) is the trajectory of $f$ (resp. $p$) starting from $x_0$. Here, the average is over $10$ random choices of $x_0$  from the box $B$, and $100$ scalars $t_i$ that form a regular subdivision of the interval $[0, 30]$.
\end{itemize}
The value $\infty$ in our table entries indicates that the trajectory of the learned vector field (as simulated by our ODE solver) is diverging.

\section*{Observations}
We make some observations on the patterns that arise in the tables below. As expected,
the training error is always higher in presence of the side information because the underlying optimization problems involve more constraints. However, 
 the test error (both in the vector field and the trajectory sense) is  lower in more than $90\%$ of the experiments.
 %
 % (percentage + median)
 %
 %For example, the last row of the first table shows that as the number of side information constraints goes from $0$  to $4$, the training error goes from $0.098$ to $0.1$, but the vector-field test error decreases from $0.63$ to $0.229$, and the trajectory test error decreases from $\infty$ (meaning that trajectories of the vector field learned without any side information constraints diverges) to $0.009$. 
 %
This indicates that side information helps with generalization and overfitting to the noise. On average, the improvement in test error is more pronounced when model complexity is high (i.e., for higher values of $d$), the noise level is high, and the number of trajectories is low. For example, the vector-field test error of the polynomial vector field of degree $2$ learned from $3$ trajectories with noise level $\sigma = 10^{-3}$ does not improve when one considers the side information constraints. However, the vector field test error of the polynomial vector field of degree $4$ learned from  $1$ trajectory with noise level $\sigma = 10^{-1}$  improves by $3$ orders of magnitudes when the side information constraints are added.

% which is an improvement of $3$ orders of magnitudes.

%      Objective value increases as we add side info, but generalization decreases. This indicates that we fight over-fitting.
%      examples: If you look at col X row Z, ...
%      Sometimes it goes up, it's rare and the amount is small (it stays in the same order of magnitude).
%     - examples: col X row Y, increase by 10\%
%  The effect is more pronounced when noise is higher
% - The effect is more pronounced when the number of trajectories is lower.

\scriptsize

\section*{Diffusion of a contagious disease, degree = $2$}\\
\begin{tabular}{ll|rr|rr|rr}
\toprule
      & {} & \multicolumn{2}{l}{training error} & \multicolumn{2}{l}{test error (vector field)} & \multicolumn{2}{l}{test error (trajectories)} \\
      & side information used? $\rightarrow$ &         no &     yes &        no &     yes &     no &     yes \\
noise & \# trajectories $\downarrow$ &           &       &          &       &       &       \\
\midrule
\multirow{3}{*}{0.001} & 1 &     0.001 & 0.001 &      189 & 0.068 &   $\infty$ & 0.008 \\
      & 2 &     0.001 & 0.002 &    0.044 & 0.081 & 0.004 & 0.005 \\
      & 3 &     0.001 & 0.002 &    0.018 & 0.083 & 0.002 & 0.005 \\
\cline{1-8}
\multirow{3}{*}{0.010} & 1 &     0.008 & 0.009 & 1.89e+03 & 0.071 &   $\infty$ & 0.012 \\
      & 2 &     0.008 & 0.009 &    0.444 & 0.061 & 0.031 & 0.004 \\
      & 3 &     0.009 & 0.009 &    0.179 & 0.074 &   $\infty$ & 0.006 \\
\cline{1-8}
\multirow{3}{*}{0.100} & 1 &      0.08 & 0.087 & 1.89e+04 & 0.075 &   $\infty$ &  0.02 \\
      & 2 &     0.085 & 0.088 &     4.44 & 0.068 &   $\infty$ & 0.019 \\
      & 3 &     0.086 & 0.088 &     1.78 &  0.11 &   $\infty$ & 0.016 \\
\bottomrule
\end{tabular}

\section*{Diffusion of a contagious disease, degree = $3$}\\
\begin{tabular}{ll|rr|rr|rr}
\toprule
      & {} & \multicolumn{2}{l}{training error} & \multicolumn{2}{l}{test error (vector field)} & \multicolumn{2}{l}{test error (trajectories)} \\
      & side information used? $\rightarrow$ &         no &     yes &        no &     yes &     no &     yes \\
noise & \# trajectories $\downarrow$ &           &       &          &       &       &       \\
\midrule
\multirow{3}{*}{0.001} & 1 &     0.001 & 0.001 &      229 &  4.22 &  $\infty$ & 0.065 \\
      & 2 &     0.001 & 0.001 &     21.8 & 0.248 &  $\infty$ & 0.005 \\
      & 3 &     0.001 & 0.001 &     1.65 & 0.151 &  $\infty$ & 0.004 \\
\cline{1-8}
\multirow{3}{*}{0.010} & 1 &     0.008 & 0.008 & 2.29e+03 &  43.4 &  $\infty$ & 0.581 \\
      & 2 &     0.008 & 0.009 &      218 &  0.47 &  $\infty$ & 0.013 \\
      & 3 &     0.008 & 0.009 &     16.5 & 0.292 &  $\infty$ & 0.007 \\
\cline{1-8}
\multirow{3}{*}{0.100} & 1 &     0.079 & 0.084 & 2.29e+04 &   435 &  $\infty$ &   5.8 \\
      & 2 &     0.082 & 0.087 & 2.18e+03 &  2.64 &  $\infty$ & 0.079 \\
      & 3 &     0.084 & 0.087 &      165 & 0.982 &  $\infty$ & 0.023 \\
\bottomrule
\end{tabular}

\section*{Diffusion of a contagious disease, degree = $4$}\\
\begin{tabular}{ll|rr|rr|rr}
\toprule
      & {} & \multicolumn{2}{l}{training error} & \multicolumn{2}{l}{test error (vector field)} & \multicolumn{2}{l}{test error (trajectories)} \\
      & side information used? $\rightarrow$ &         no &     yes &        no &     yes &     no &     yes \\
noise & \# trajectories $\downarrow$ &           &       &          &       &       &       \\
\midrule
\multirow{3}{*}{0.001} & 1 &     0.001 & 0.001 & 1.63e+04 &     21.7 &  $\infty$ & 0.092 \\
      & 2 &     0.001 & 0.001 & 2.58e+03 &    0.726 &  $\infty$ & 0.004 \\
      & 3 &     0.001 & 0.001 &      421 &    0.575 &  $\infty$ & 0.002 \\
\cline{1-8}
\multirow{3}{*}{0.010} & 1 &     0.008 & 0.008 & 1.62e+05 &      191 &  $\infty$ & 0.822 \\
      & 2 &     0.008 & 0.009 & 2.58e+04 &     1.83 &  $\infty$ & 0.046 \\
      & 3 &     0.008 & 0.009 & 4.21e+03 &     1.53 &  $\infty$ & 0.031 \\
\cline{1-8}
\multirow{3}{*}{0.100} & 1 &     0.076 & 0.083 & 1.57e+06 & 1.92e+03 &  $\infty$ &  8.22 \\
      & 2 &     0.079 & 0.087 & 2.58e+05 &     8.01 &  $\infty$ & 0.018 \\
      & 3 &     0.081 & 0.087 & 4.21e+04 &     7.99 &  $\infty$ & 0.017 \\
\bottomrule
\end{tabular}

\section*{Lorenz, degree = $2$}\\
\begin{tabular}{ll|rr|rr|rr}
\toprule
      & {} & \multicolumn{2}{l}{training error} & \multicolumn{2}{l}{test error (vector field)} & \multicolumn{2}{l}{test error (trajectories)} \\
      & side information used? $\rightarrow$ &         no &     yes &        no &     yes &     no &     yes \\
noise & \# trajectories $\downarrow$ &           &       &          &       &       &       \\
\midrule
\multirow{3}{*}{0.001} & 1 &     0.001 & 0.001 & 0.049 & 0.001 & 0.028 & 0.002 \\
      & 2 &     0.001 & 0.001 & 0.002 & 0.001 & 0.002 & 0.001 \\
      & 3 &     0.001 & 0.001 & 0.002 & 0.001 & 0.003 & 0.001 \\
\cline{1-8}
\multirow{3}{*}{0.010} & 1 &     0.007 & 0.008 & 0.491 & 0.009 & 0.285 & 0.018 \\
      & 2 &     0.008 & 0.009 &  0.02 & 0.007 & 0.018 & 0.008 \\
      & 3 &     0.008 & 0.009 & 0.018 & 0.007 & 0.025 & 0.005 \\
\cline{1-8}
\multirow{3}{*}{0.100} & 1 &     0.068 & 0.083 &  4.91 &  0.09 &  3.01 &  0.17 \\
      & 2 &     0.079 & 0.086 & 0.202 & 0.074 & 0.185 & 0.064 \\
      & 3 &     0.083 & 0.088 & 0.185 & 0.071 & 0.251 & 0.052 \\
\bottomrule
\end{tabular}

\section*{Lorenz, degree = $3$}\\
\begin{tabular}{ll|rr|rr|rr}
\toprule
      & {} & \multicolumn{2}{l}{training error} & \multicolumn{2}{l}{test error (vector field)} & \multicolumn{2}{l}{test error (trajectories)} \\
      & side information used? $\rightarrow$ &         no &     yes &        no &     yes &     no &     yes \\
noise & \# trajectories $\downarrow$ &           &       &          &       &       &       \\
\midrule
\multirow{3}{*}{0.001} & 1 &         0 &     0 &      110 &  10.9 &   $\infty$ &  13.7 \\
      & 2 &     0.001 & 0.001 &    0.046 & 0.008 & 0.016 & 0.005 \\
      & 3 &     0.001 & 0.001 &    0.014 & 0.044 & 0.006 & 0.034 \\
\cline{1-8}
\multirow{3}{*}{0.010} & 1 &     0.003 & 0.005 &      192 &  10.9 &  57.9 &  14.2 \\
      & 2 &     0.006 & 0.007 &    0.461 & 0.063 & 0.157 &  0.05 \\
      & 3 &     0.008 & 0.008 &    0.138 & 0.035 & 0.058 & 0.019 \\
\cline{1-8}
\multirow{3}{*}{0.100} & 1 &     0.032 & 0.057 & 1.47e+03 &  54.9 &   237 &  23.6 \\
      & 2 &     0.064 & 0.071 &     4.61 & 0.963 &  1.59 & 0.461 \\
      & 3 &     0.075 & 0.078 &     1.38 & 0.352 & 0.576 & 0.279 \\
\bottomrule
\end{tabular}
\newpage
\section*{Lorenz, degree = $4$}\\
\begin{tabular}{ll|rr|rr|rr}
\toprule
      & {} & \multicolumn{2}{l}{training error} & \multicolumn{2}{l}{test error (vector field)} & \multicolumn{2}{l}{test error (trajectories)} \\
      & side information used? $\rightarrow$ &         no &     yes &        no &     yes &     no &     yes \\
noise & \# trajectories $\downarrow$ &           &       &          &       &       &       \\
\midrule
\multirow{3}{*}{0.001} & 1 &         0 &     0 &      167 &      532 &  82.7 &  94.4 \\
      & 2 &         0 &     0 &     46.1 &    0.479 &  51.8 & 0.832 \\
      & 3 &     0.001 &     0 &    0.283 &    0.205 & 0.047 & 0.421 \\
\cline{1-8}
\multirow{3}{*}{0.010} & 1 &         0 & 0.001 &      858 & 6.22e+03 &   $\infty$ &   286 \\
      & 2 &     0.003 & 0.004 &      461 &     1.71 &   $\infty$ &  1.97 \\
      & 3 &     0.006 & 0.005 &     2.83 &    0.407 & 0.468 & 0.179 \\
\cline{1-8}
\multirow{3}{*}{0.100} & 1 &         0 & 0.031 & 8.92e+03 & 2.55e+03 &   $\infty$ &   120 \\
      & 2 &     0.028 & 0.057 & 4.61e+03 &     6.21 &   $\infty$ &  5.15 \\
      & 3 &     0.061 & 0.069 &     28.3 &     2.47 &  4.84 & 0.939 \\
\bottomrule
\end{tabular}

%\color{black}
\small
\bibliographystyle{custom_plain}
\bibliography{citations}

% todo: how to know when to stop adding side info?

\end{document}